
\input amstex

\define\scrO{\Cal O}
\define\Pee{{\Bbb P}}
\define\Zee{{\Bbb Z}}
\define\Cee{{\Bbb C}}
\define\Ar{{\Bbb R}}

\define\ddual{\spcheck{}\spcheck}

\define\Ker{\operatorname{Ker}}

\define\Hom{\operatorname{Hom}}

\define\Id{\operatorname{Id}}

\define\ch{\operatorname{ch}}
\define\Todd{\operatorname{Todd}}
\define\Ext{\operatorname{Ext}}
\define\Spec{\operatorname{Spec}}
\define\Hilb{\operatorname{Hilb}}
\define\rank{\operatorname{rank}}

\define\proof{\demo{Proof}}
\define\endproof{\qed\enddemo}
\define\endstatement{\endproclaim}
\define\theorem#1{\proclaim{Theorem #1}}
\define\lemma#1{\proclaim{Lemma #1}}
\define\proposition#1{\proclaim{Proposition #1}}
\define\corollary#1{\proclaim{Corollary #1}}
\define\claim#1{\proclaim{Claim #1}}

\define\section#1{\specialhead #1 \endspecialhead}
\define\ssection#1{\medskip\noindent{\bf #1}}
\loadbold

\documentstyle{amsppt}
\leftheadtext{Robert Friedman and Zhenbo Qin}
\rightheadtext{Flips and transition formulas}

\topmatter
\title
Flips of moduli spaces and transition formulas for Donaldson polynomial
invariants of rational surfaces
\endtitle
\author {Robert Friedman  and Zhenbo Qin}
\endauthor
\address Department of Mathematics, Columbia University, New York,
NY 10027, USA\endaddress
\email rf\@math.columbia.edu  \endemail
\address Department of Mathematics, Oklahoma State University,
Stillwater, OK 74078, USA
\endaddress
\email  qz\@math.okstate.edu \endemail
\thanks The first author was partially supported by NSF grants DMS-90-06116 and
DMS-92-03940. The second author was partially supported by a grant from the
ORAU
Junior Faculty Enhancement Award Program and by NSF grants DMS-91-00383 and
DMS-94-00729.
\endthanks
\endtopmatter

\NoBlackBoxes

\document
\section{1. Introduction.}

In \cite{7}, Donaldson has defined polynomial invariants for
smooth simply connected 4-manifolds with $b_2^+\geq 3$.
These invariants have also been defined for 4-manifolds with $b_2^+=1$
in \cite{24, 17, 18}, along lines suggested by the work of
Donaldson in \cite{5}.  In this case, however,
they depend on an additional piece of information, namely a chamber
defined on the positive cone of $H^2(X; \Ar)$ by
a certain locally finite set of walls. Explicitly,
let $X$ be a simply connected, oriented, and closed smooth
$4$-manifold with $b_2^+ = 1$ where $b_2^+$ is the number
of positive eigenvalues of the quadratic form $q_X$
when diagonalized over $\Ar$. Let
$$\Omega_X = \{\, x \in H^2(X, \Ar)\mid x^2 > 0 \,\}$$
be the positive cone. Fix a class $\Delta$ in $H^2(X, \Zee)$ and
an integer $c$ such that $d = 4c - \Delta ^2 - 3$ is nonnegative.
A {\sl wall of type $(\Delta, c)$\/}  is a nonempty
hyperplane:
$$W^\zeta = \{\, x \in \Omega_X\mid x \cdot \zeta = 0 \,\}$$
in $\Omega_X$ for some class $\zeta \in H^2(X, \Zee)$
with $\zeta \equiv \Delta \pmod 2$
and $\Delta ^2 - 4c \le \zeta^2 < 0$.
The connected components of the complement in $\Omega_X$ of the walls of  type
$(\Delta, c)$ are the {\sl chambers of type $(\Delta, c)$}.
Then the Donaldson polynomial invariants of $X$ associated to
$\Delta$ and $c$ are defined with respect to chambers of type $(\Delta, c)$.
The invariants only depend on the class $w=\Delta \bmod 2\in H^2(X;
\Zee/2\Zee)$
and the integer $p =  \Delta ^2 - 4c$, and we shall often refer to walls and
chambers of type $(w,p)$ as well. We shall write $D^X _{w,p}(\Cal C)$
for the Donaldson polynomial corresponding to the $SO(3)$ bundle $P$ with
invariants $w_2(P) = w$ and $p_1(P) = p$, depending on the chamber $\Cal C$.

A basic question is then the following: Suppose that $\Cal C_+$ and $\Cal C_-$
are separated by a single wall $W^\zeta$. Here there may be more than one class
$\zeta$ of type $(\Delta, c)$ defining $W^\zeta$. Then find a formula for the
difference
$$\delta ^X_{w,p}(\Cal C_+, \Cal C_-) = D^X _{w,p}(\Cal C_+) - D^X _{w,p}(\Cal
C_-).$$
We shall refer to such a difference as a {\sl transition formula}.

There has been considerable interest in the above problem. The first
result in this direction is due to Donaldson in \cite{5},
who gave a formula in case $\Delta = 0$ and $c = 1$.
Kotschick \cite{17} showed that, on the part of the symmetric algebra generated
by
$2$-dimensional classes,
$\delta ^X_{w,p}(\Cal C_+, \Cal C_-) = \pm \zeta^d$ for
$\zeta^2 = -(4c - \Delta^2)=p$, and that $\delta ^X_{w,p}(\Cal C_+, \Cal C_-)$
is in fact always divisible by $\zeta$, except when $p=-5$ and $\zeta ^2 = -1$
(cf\. also Mong \cite{24} for some partial results along these lines).
For a rational ruled surface
$X$, all the transition formulas for $\Delta  = 0$ and $2 \le c \le 4$
have been determined in \cite{24, 33, 22}. Using a gauge-theoretic approach,
Yang \cite{35} settled the problem for $\Delta = 0$ and $c = 2$, and computed
the degree $5$ Donaldson polynomials for rational surfaces. The known
examples and the work of Kotschick and Morgan \cite{18} raise
the following rather natural conjecture:

\medskip\noindent
{\bf Conjecture.} {\it The transition formula $\delta ^X_{w,p}(\Cal C_+, \Cal
C_-)$ is a homotopy invariant of the pair $(X, \zeta)$; more precisely,
if $\phi$ is an oriented homotopy equivalence from $X'$ to $X$, then
$$\delta^{X'}_{\phi ^*w,p}(\phi ^*(\Cal C_+), \phi ^*(\Cal
C_-))=\phi ^*\delta ^{X}_{w,p}(\Cal C_+, \Cal C_-).$$}
\medskip

We remark that this conjecture is essentially equivalent to the following
statement: the transition formula $\delta ^X_{w,p}(\Cal C_+, \Cal C_-)$ is a
polynomial  in $\zeta$ and the quadratic form $q_X$ with coefficients
involving only $\zeta^2$, homotopy invariants of $X$ (i.e\. $b^-_2(X)$), and
universal constants.

Our goal in this paper is to study the corresponding problem in
algebraic geometry. More precisely, let $X$ be an algebraic surface (not
necessarily with $b_2^+(X)= 1$) and let $L$ be an ample line bundle on $X$. We
can then identify the moduli space of $L$-stable rank two bundles $V$ on $X$
with
$c_1(V)=\Delta$ and $c_2(V)=c$ with the moduli space of equivalence classes of
ASD connections on $X$ with respect to a Hodge metric on $X$ corresponding to
$L$. Let $\frak M_L(\Delta, c)$ be the Gieseker compactification of this moduli
space. It is known that $\frak M_L(\Delta, c)$ changes as we change $L$, and
that $\frak M_L(\Delta, c)$ is constant on a set of chambers for the ample cone
of $X$ which are defined in a way analogous to the definition of chambers for
$\Omega _X$ given above. Using the recent result of Morgan \cite{25} and
Li \cite{21} that the Donaldson polynomial of an algebraic surface
can be evaluated using the Gieseker compactification ${\frak M}_L(\Delta, c)$
of the moduli space of stable bundles, we shall work on
${\frak M}_L(\Delta,  c)$ for  suitable choices of $L$ and in particular
analyze the change in ${\frak M}_L(\Delta, c)$ for $L\in \Cal C_+$ or
$L\in \Cal C_-$, where $\Cal C_\pm$ are two adjacent chambers.
It turns out that we can obtain $\frak M _{L_+}(\Delta, c)$ from
$\frak M _{L_-}(\Delta, c)$ by a series of blowups and blowdowns
(flips). Our results are thus very similar to those of
Thaddeus in \cite{31}. Thaddeus \cite{32} and also Dolgachev-Hu \cite{3}
have developed a general picture for the variation of GIT quotients
after a change of polarization, and although our methods are somewhat different
it seems quite possible that they fit into
their general framework. We have also found it convenient to borrow
some of Thaddeus' notation.

Next we shall apply our results on the change in the moduli spaces to determine
the transition formula for Donaldson polynomials in case
$X$ is a rational surface with
$-K_X$ effective. We shall give explicit formulas for
$\delta ^X_{w,p}(\Cal C_+, \Cal C_-)$ in case the nonnegative integer
$\ell _\zeta = (\zeta ^2 - p)/4 \leq 2$. These formulas are in agreement
with the above conjecture, in the sense that the transition formula is
indeed a polynomial in $\zeta$ and  $q_X$ with coefficients involving only
$\zeta^2$, $K_X^2$, and universal constants. We shall also give a
formula in principle for $\delta ^X_{w,p}(\Cal C_+, \Cal C_-)$ in general
(see Theorem 5.4),
but to make this formula explicit involves more knowledge of the enumerative
geometry of $\Hilb^nX$ than seems to be available at present.
In case $-K_X$ is effective, the moduli spaces are (essentially) smooth and
the centers of the blowup are smooth as well; in fact they are
$\Pee ^k$-bundles over $\Hilb ^{n_1}X\times \Hilb ^{n_2}X$ for appropriate $k$,
$n_1$ and $n_2$. In this way, we obtain general formulas which can be made
explicit for low values of $n$. For instance, we show the following
(see Theorem 6.4 for details):

\theorem{} Assume that the wall $W^\zeta$ is defined only by
$\pm \zeta$ with $\ell_\zeta = 1$ and that $\Cal C_\pm$ lies
on the $\pm$-side of $W^\zeta$. Then, on the subspace of the symmetric algebra
generated by $H_2(X)$,
$\delta^X_{w, p}(\Cal C_-, \Cal C_+)$ is equal to
$$(-1)^{{(\Delta \cdot K_X + \Delta^2) + (\zeta \cdot K_X - \zeta^2)} \over 2}
\cdot \left \{ d(d - 1) \cdot \left(\zeta \over 2 \right)^{d - 2} \cdot q_X
+ (2K_X^2 + 2d + 6) \cdot \left ({\zeta \over 2} \right )^d \right \}.$$
\endstatement

Along the direction of the work of Kronheimer and Mrowka \cite{19, 20},
we also consider the difference of Donaldson polynomial invariants
involving the natural generator $x \in H_0(X; \Zee)$. More precisely,
let $\nu$ be the corresponding $4$-dimensional class in
the instanton moduli space. For $\alpha \in H_2(X; \Zee)$,
we give a formula for the difference
$\delta^X_{w, p}(\Cal C_-, \Cal C_+)(\alpha^{d - 2}, \nu)$
in Theorem 5.5. It is worth to point out that the similarity between
Theorem 5.4 and Theorem 5.5 may indicate that there exists
a deep relation between $\delta^X_{w, p}(\Cal C_-, \Cal C_+)(\alpha^d)$ and
$\delta^X_{w, p}(\Cal C_-, \Cal C_+)(\alpha^{d - 2}, \nu)$,
and suggest a way to generalize the notion of {\sl simple type} in
\cite{19, 20} from the case of $b_2^+ > 1$ to the case of $b_2^+ = 1$.
For instance, modulo some lower degree terms,
$\delta^X_{w, p}(\Cal C_-, \Cal C_+)(\alpha^{d - 2}, \nu)$
can be obtained from $(-1/4) \cdot \delta^X_{w, p}(\Cal C_-, \Cal
C_+)(\alpha^d)$ by replacing $d$ by $(d - 2)$
(see Theorem 5.13 and Theorem 5.14).  In fact, based on some heuristic
arguments, it seems reasonable to conjecture that
$\delta^X_{w, p}(\Cal C_-, \Cal C_+)(\alpha^{d - 2}, \nu)$
is a combination of $\delta^X_{w, p'}(\Cal C_-, \Cal C_+)(\alpha^{d - 4k})$
for various nonnegative integers $k$ if the degrees are properly arranged.
We hope to return to this issue in future.

Our paper is organized as follows. In section 2, we study rank two
torsion free sheaves which are semistable with respect to ample
divisors in $\Cal C_-$ but not semistable with respect to ample
divisors in $\Cal C_+$. When the surface $X$
is rational with $-K_X$ effective, these sheaves are parametrized by an open
subset of a union of projective  bundles over the product of two Hilbert
schemes
of points in
$X$. More precisely, if $\zeta$ defines the wall separating $\Cal C_-$ from
$\Cal
C_+$, define
$E_\zeta^{n_1, n_2}$ be the set of  all isomorphism classes of nonsplit
extensions of the form
$$0 \to \scrO _X(F )\otimes I_{Z_1} \to V \to \scrO _X(\Delta -F ) \otimes
I_{Z_2} \to 0,$$
where $F$ is a divisor class such that $2F-\Delta \equiv \zeta$ and $Z_1$ and
$Z_2$ are two zero-dimensional subschemes of $X$ with $\ell (Z_i) = n_i$
such that $n_1 + n_2 = \ell _\zeta$. In case $X$ is rational,
$E_\zeta^{n_1, n_2}$ is a $\Pee ^N$ bundle over
$\Hilb ^{n_1}X\times \Hilb ^{n_2}X$, and the set of points of $E_\zeta^{n_1,
n_2}$  lying in $\frak M_{L_-}(\Delta, c)$ but not in $\frak M_{L_+}(\Delta,
c)$
is a Zariski open subset of $E_\zeta^{n_1, n_2}$. The main technical difficulty
is that it is hard to control the rational map from $E_\zeta^{n_1, n_2}$ to
$\frak M_{L_-}(\Delta, c)$, and in particular this map is not a morphism. The
general picture that we establish is the following: first, the map $E_\zeta^{0,
\ell_\zeta}\dasharrow\frak M_{L_-}(\Delta, c)$ is a morphism, and it is
possible
to make an elementary transformation, or {\sl flip\/}, along its image. The
result is a new space for which the rational map $E_\zeta^{1,
\ell_\zeta - 1}\dasharrow\frak M_{L_-}(\Delta, c)$ becomes a morphism, and it
is
possible to make a  flip along {\sl its\/} image. We continue in this way until
we reach $\frak M_{L_+}(\Delta, c)$.

It seems rather difficult to see that the above picture holds directly. Instead
we shall proceed as follows. We define abstractly a sequence of moduli spaces,
indexed by an integer $k$ with $0\leq k \leq \ell _\zeta+1$, such that the
moduli space for $k=0$ is $\frak M_{L_-}(\Delta, c)$, the moduli space for
$k=\ell _\zeta+1$ is $\frak M_{L_+}(\Delta, c)$, and moreover the
$k^{\text{th}}$ moduli space contains an embedded copy of $E_\zeta^{k, \ell
_\zeta -k}$ such that the flip along this copy yields the $(k+1)^{\text{st}}$
moduli space. Thus the picture is very similar to that developed independently
by Thaddeus in \cite{31}. To define our sequence of moduli spaces,
we define $(L_0, \boldsymbol \zeta, \bold k)$-semistability in section 3 for
rank two torsion free sheaves, where $L_0$ is any ample divisor contained  in
the common face of
$\Cal C_+$ and $\Cal C_-$,
$\boldsymbol
\zeta$ is the set of classes of type $(\Delta, c)$ defining the common wall
of $\Cal C_+$ and $\Cal C_-$, and $\bold k$ is a set of integers.
We show that $\frak M _{L_-}(\Delta, c)$ and $\frak M _{L_+}(\Delta, c)$
are linked by the moduli spaces $\frak M_0 ^{(\boldsymbol \zeta, \bold k)}$
where the data $\bold k$ is allowed to vary. When the surface $X$ is
rational with $-K_X$ effective, we can obtain $\frak M _{L_+}(\Delta, c)$
from $\frak M _{L_-}(\Delta, c)$ by a series of flips.
The fact that all $(L_0, \boldsymbol \zeta, \bold k)$-semistable
rank two torsion free sheaves do form a moduli space
$\frak M_0 ^{(\boldsymbol \zeta, \bold k)}$ in the usual sense
is proved in section 4 where we introduce an equivalent notion of stability
called {\sl mixed stability}. Our method follows Gieseker's GIT argument in
\cite{13}. Roughly speaking, the goal of mixed stability is to define stability
for a sheaf of the form $V\otimes \Xi$, where $V$ is a torsion free sheaf but
$\Xi$ is just a $\Bbb Q$-divisor. To make this idea precise, given
actual divisors
$H_1$ and $H_2$ and positive weights $a_1$ and $a_2$, we shall define a notion
of stability which ``mixes" stability for $V\otimes H_1$ with stability for
$V\otimes H_2$, together with weightings of the stability condition for
$V\otimes H_i$. The effect of this definition will be formally the same as if
we had defined stability of $V\otimes \Xi$, where $\Xi$ is the $\Bbb
Q$-divisor
$$\frac{a_1}{a_1+a_2}H_1 + \frac{a_2}{a_1+a_2}H_2.$$

In section 5, using our results on flips of moduli spaces,
we give a formula for the transition formula of Donaldson polynomials
when $X$ is rational with $-K_X$ effective,
and compute the leading term in the transition formula.
In section 6, we obtain explicit transition formulas when $\ell_\zeta \le 2$.

Some of the material in our section 2 has been worked out independently
by Hu and Li \cite{16} and G\"ottsche \cite{14}.
Moreover Ellingsrud and G\"ottsche \cite{8} have recently studied the change
in the moduli space by similar methods and have obtained results
very similar to ours. Using very different methods, the results in
Section 4 have also been obtained  by Matsuki and Wentworth \cite{23},
who also consider the case of higher rank. They use branched covers of
the surface $X$ to study the change in the moduli space.
We expect that a minor modification of the arguments in Section 4 of
this paper will also handle the case of higher rank.

\section{Conventions and notations}

We fix some conventions and notations for the rest of this paper.
Let $X$ be a smooth algebraic surface. We shall be primarily interested
in the case where $X$ is simply connected and $-K_X$ is effective and
nonzero. Thus necessarily $X$ is a rational surface. However much of the
discussion
in sections 1--4 will also apply to the general case.
Stability and semistability with respect to an ample line bundle $L$
will always be understood to mean Gieseker
stability or semistability unless otherwise noted. We shall not mention the
choice of $L$ explicitly if it is clear from the context.
Recall that a torsion free sheaf $V$ of rank two is Gieseker $L$-stable
if and only if, for every rank one subsheaf $W$ of $V$,
either $\mu _L(W) < \mu _L(V)$ or $\mu _L(W) = \mu _L(V)$ and
$2\chi(W) < \chi(V)$, where $\mu _L$ is the normalized degree
with respect to $L$. Semistability is similarly defined, where the second
inequality is also allowed to be an equality.
For a torsion free sheaf $V$, we use $V\ddual$ to stand for its double dual.
For two divisors $D_1$ and $D_2$ on $X$, the notation $D_1 \equiv D_2$
means that $D_1$ and $D_2$ are numerically equivalent, that is,
$D_1 \cdot D = D_2 \cdot D$ for any divisor $D$.
For a locally free sheaf (or equivalently a vector bundle) $\Cal E$
over a smooth variety $Y$, we use $\Pee(\Cal E)$ to denote
the associated projective space bundle,
that is, $\Pee(\Cal E)$ is the {\bf Proj} of $\oplus_{d \ge 0} S^d(\Cal E)$.

Fix a divisor $\Delta$ and an integer $c$. Let $\Cal C_-$ and $\Cal C_+$
be two adjacent chambers of type $(\Delta, c)$ separated by the wall
$W^\zeta$. We assume that $\zeta \cdot \Cal C_- < 0<\zeta \cdot \Cal C_+$.
Let $L_\pm\in \Cal C_\pm$ be an ample line bundle,
so that $L_- \cdot \zeta < 0 < L_+\cdot \zeta$, and denote by
$\frak M_\pm$ the moduli space $\frak M _{L_\pm}(\Delta, c)$
of rank two Gieseker semistable torsion free sheaves $V$ with
$c_1(V) = \Delta$ and $c_2(V) = c$. Let $L_0$ be any ample divisor
contained in the interior of the intersection of $W^\zeta$ and
the closures of $\Cal C_\pm$. Let $\zeta = \zeta _1,
\dots, \zeta _n$ be all the positive rational multiples of $\zeta$
such that $\zeta_i$ is an integral class of type $(w,p)$ which also
defines the wall $W^\zeta$. In sections 5--6, we will assume that
$n = 1$ for notational simplicity.

Finally, we point out that our $\mu$-map is half of the $\mu$-map
used in \cite{17, 18} (see (viii) and (ix) in Notation 5.1).
Thus our transition formula differs from the one defined in \cite{18}
by a universal constant.

\medskip\noindent
{\it Acknowledgements.} We would like to thank Hong-Jie Yang for invaluable
access to his calculations, which helped to keep us on the right track.
The second author would like to thank Wei-ping Li and Yun-Gang Ye for
helpful discussions, and the Institute for Advanced Study at
Princeton for its hospitality and financial support through NSF grant
DMS-9100383 during the academic year 1992--1993 when part of this work was
done.

\section{2. Preliminaries on the moduli space.}

In this section, we study rank two torsion free sheaves which are
related to walls. These sheaves arise naturally from the comparison of
$L_-$-semistability and $L_+$-semistability. We will show that
when the surface $X$ is rational with $-K_X$ effective, the moduli spaces
$\frak M_\pm$ are smooth at the points corresponding to these sheaves.
We start with the following lemma, which for simplicity is just stated for
$L_-$-stability.

\lemma{2.1} Let $V$ be a rank two torsion free sheaf on $X$ with
$c_1(V) = \Delta$ and $c_2(V) = c$. If $V$ is $L_-$-semistable,
then exactly one of the following holds:
\roster
\item"{(i)}" Both $V$ and $V\ddual$ are $L_-$-stable and Mumford
$L_-$-stable.
\item"{(ii)}" $V$ sits in an exact sequence
$$0 \to \scrO_X(F_1)\otimes I_{Z_1} \to V \to \scrO_X(F_2)
\otimes I_{Z_2} \to 0$$
where $2F_1 \equiv \Delta \equiv 2F_2$, and $Z_1$ and $Z_2$ are
zero-dimensional subschemes of $X$ such that
$\ell (Z_1) \geq \ell(Z_2)$. Moreover in this case $V$ is $L$-semistable for
every choice of an ample line bundle $L$ and $V$ is strictly
$L_\pm$-semistable if and only if $\ell(Z_1) = \ell (Z_2)$.
\endroster
\endstatement
\proof Suppose that $V$ is (Gieseker) $L_-$-semistable. The vector bundle
$V\ddual$ satisfies $c_1(V\ddual ) = \Delta$ and $c_2(V\ddual) \leq c$.
Standard arguments \cite{10} show that $V\ddual$ is Mumford
$L_-$-semistable. If $V\ddual$ is strictly Mumford
$L_-$-semistable, then by \cite{10, 30},
either $L_-$ must lie on a wall of type $(\Delta, c)$ or
if $\scrO_X(F_1)$ is a destabilizing sub-line bundle then $\Delta \equiv 2F_1$.
Since by assumption $L_-$ does not lie on a wall of type $(\Delta, c)$,
either $V\ddual$ is Mumford $L_-$-stable or there is an exact sequence
$$0 \to \scrO_X(F_1) \to V\ddual \to \scrO_X(F_2) \otimes I_Z \to 0,$$
where $F_2 = \Delta - F_1 \equiv F_1$ and $Z$ is a zero-dimensional subscheme
of $X$. If $V\ddual$ is Mumford $L_-$-stable, then $V$ is Mumford
$L_-$-stable and therefore $L_-$-stable. Thus case (i) holds. Otherwise
$\scrO_X(F_1) \cap V$ is of the form $\scrO_X(F_1)\otimes I_{Z_1}$ for
some $Z_1$ and $V/\scrO_X(F_1)\otimes I_{Z_1}$ is a subsheaf of
$\scrO_X(F_2) \otimes I_Z$ and thus of the form
$\scrO_X(F_2) \otimes I_{Z_2}$ for some $Z_2$. Thus we are in case (ii) of the
lemma. Since $\mu(\scrO_X(F_1)\otimes I_{Z_1}) = \mu (V)$
and $V$ is semistable, we have
$$2\chi(\scrO_X(F_1)\otimes I_{Z_1}) \leq \chi (V)=
\chi(\scrO_X(F_1) \otimes I_{Z_1})+ \chi(\scrO_X(F_2)\otimes I_{Z_2}).$$
Hence $\chi(\scrO_X(F_2)\otimes I_{Z_2}) -
\chi(\scrO_X(F_1)\otimes I_{Z_1}) \geq 0$.
As $F_1 \equiv F_2$ and $\chi(\scrO_X(F_i)\otimes I_{Z_i})
= \chi(\scrO_X(F_i)) - \ell(Z_i)$, we must then have
$\ell(Z_1) - \ell (Z_2) \geq 0$. The last sentence of (ii) is a
straightforward argument left to the reader.
\endproof

If $V$ satisfies the conclusions of (2.1)(ii), we shall call $V$ {\sl
universally semistable}.

Next we shall compare stability for $L_-$ and $L_+$.

\lemma{2.2} Let $V$ be a torsion free rank two sheaf with $c_1(V) = \Delta$ and
$c_2(V) = c$.
\roster
\item"{(i)}" If $V$ is $L_-$-stable but $L_+$-unstable,
then there exist a divisor class $F$ and two zero-dimensional subschemes
$Z_-$ and $Z_+$ of $X$ and an exact sequence
$$0 \to \scrO _X(F )\otimes I_{Z_-} \to V \to \scrO _X(\Delta -F ) \otimes
I_{Z_+} \to 0,$$ where $L_-\cdot (2F-\Delta) < 0 < L_+\cdot (2F-\Delta) $.
Moreover the divisor $F$, the schemes $Z_-$ and $Z_+$, and the map $F \otimes
I_{Z_-} \to V$ are unique mod scalars, and $\zeta = 2F - \Delta$ defines
a wall of type $(\Delta, c)$.

\item"{(ii)}" Conversely, suppose that there is a nonsplit exact sequence as
above. Then $V$ is simple. Moreover, $V$ is not $L_-$-stable if and
only if it is $L_-$-unstable if and only if there exist subschemes
$Z'$ and $Z''$ and an exact sequence
$$0 \to \scrO _X(\Delta - F )\otimes I_{Z'} \to V \to \scrO _X(F )
\otimes I_{Z''} \to 0,$$
if and only if $V\ddual$ is a direct sum $\scrO _X(F) \oplus
\scrO _X(\Delta -F )$. In this case the scheme $Z'$ strictly contains
the scheme $Z_+$, $\ell (Z') > \ell (Z_+)$
and $\ell (Z')+\ell (Z'') = \ell (Z_-)+\ell (Z_+)$.
Finally, if $Z_-=\emptyset$ then $V$ is always $L_-$-stable.
\endroster
\endstatement
\proof We first show (i). Suppose that $V$ is $L_-$-stable but $L_+$-unstable.
Then by (2.1) $V\ddual$ is also $L_-$-stable and $L_+$-unstable.
By \cite{30}, there is a uniquely determined line bundle $\scrO_X(F)$ and
a map $\scrO _X(F) \to V\ddual$ with torsion free quotient such that
$L_-\cdot (2F-\Delta) < 0 < L_+\cdot (2F-\Delta) $.
Moreover $\zeta = 2F - \Delta$ defines a wall of type $(\Delta, c)$.
The subsheaf $\scrO_X(F) \cap V$ of $V\ddual$ is
a subsheaf of $\scrO_X(F)$ and agrees with it away from finitely many points.
Thus $\scrO_X(F) \cap V = \scrO _X(F )\otimes I_{Z_-}$
for some well-defined subscheme $Z_-$. Moreover the quotient
$V\Big/[\scrO _X(F )\otimes I_{Z_-}]$ is a subsheaf
of $\scrO _X(\Delta -F ) \otimes I_Z$ for some zero-dimensional subscheme
$Z$, and agrees with $\scrO _X(\Delta -F ) $ away from finitely many points.
Thus the quotient is of the form $\scrO _X(\Delta -F ) \otimes I_{Z_+}$
for some zero-dimensional subscheme $Z_+$. The uniqueness is clear.

To see (ii), suppose that $V$ is given as a nonsplit exact sequence
$$0 \to \scrO _X(F )\otimes I_{Z_-} \to V \to \scrO _X(\Delta -F ) \otimes
I_{Z_+} \to 0$$ as above, where $L_-\cdot (2F-\Delta) < 0 < L_+\cdot
(2F-\Delta)
$. Again by (2.1), $V$ is $L_-$-semistable if and only if it is $L_-$-stable if
and only if
$V\ddual$ is
$L_-$-stable. Now taking double duals of the above exact sequence, there is an
exact sequence
$$0 \to \scrO _X(F ) \to V\ddual \to \scrO _X(\Delta -F ) \otimes I_Z \to 0$$
for some zero-dimensional scheme $Z$. Moreover, by \cite{30}, $V\ddual$ is
$L_-$-unstable if and only if the above exact sequence splits, and in
particular
if and only if $Z=\emptyset$ and $V\ddual = \scrO _X(F)\oplus \scrO _X(\Delta
-F
)$. In this case, the map $\scrO _X(\Delta -F )
\to V\ddual$ induces a map $\scrO _X(\Delta -F )\otimes I_{Z'} \to V$ for some
ideal sheaf $I_{Z'}$. We may clearly assume that the quotient is torsion free,
in
which case it is necessarily of the form $\scrO _X(F )
\otimes I_{Z''}$ with $\ell (Z')+\ell (Z'') = \ell (Z_-)+\ell (Z_+)$. Using the
nonzero map
$\scrO _X(\Delta -F )\otimes I_{Z'} \to
\scrO _X(\Delta -F )\otimes I_{Z_+}$, we see that there is an inclusion $I_{Z'}
\subseteq I_{Z_+}$; moreover this inclusion must be strict since the defining
exact sequence for $V$ is nonsplit. Thus $Z'$ strictly contains $Z_+$ and in
particular $\ell (Z') > \ell (Z_+)$. Conversely, if there exists a nonzero map
$\scrO _X(\Delta - F )\otimes I_{Z'} \to V$, then there is a nonzero map $\scrO
_X(\Delta - F ) \to V\ddual$ and thus $V\ddual $ is the split extension.

We next show that $V$ is simple. If $V$ is stable then it is simple. If $V$ is
not stable, then $V\ddual = \scrO _X(F)\oplus \scrO _X(\Delta -F )$. There is
an
inclusion $\Hom (V, V) \subseteq
\Hom (V\ddual, V\ddual)$. If $V\ddual$ is split, then $\Hom (V\ddual, V\ddual)
=
\Cee \oplus \Cee$. In this case, using  a nonscalar endomorphism of $V$, it is
easy to see that we can split the exact sequence defining $V$.

Finally suppose that $Z_- = \emptyset$ in the notation of (2.2). If $V$ is
$L_-$-unstable, then we can find $Z'$ with $\ell (Z') > \ell (Z_+)$ and a
subscheme $Z''$ such that $\ell (Z') + \ell (Z'') = \ell (Z_+)$. Thus $\ell
(Z')
\leq \ell (Z_+)$, a contradiction. It follows that $V$ is $L_-$-stable.
\endproof

For the rest of this section, we shall assume that $-K_X$ is effective and
nonzero and that $q(X)=0$. Thus $X$ is a rational surface.

\lemma{2.3} Suppose that $\frak M_\pm$ is nonempty. Suppose that $(w,p)\neq
(0,0)$, or equivalently that $\frak M_\pm$ does not consist of a single point
corresponding to a twist of the trivial vector bundle. Then the open subset of
$\frak M_\pm$ corresponding to Mumford stable rank two vector bundles is
nonempty
and  dense. Every component of $\frak
M_\pm$ has dimension
$4c-\Delta ^2 -3=-p-3$. The points of $\frak M_\pm$ corresponding to
$L_\pm$-stable sheaves $V$ are smooth points.
\endstatement
\proof Suppose that $\frak M_\pm$ is nonempty, and let
$V$ correspond to a point of $\frak M_\pm$. Then by general theory (e.g\.
Chapter 7 of \cite{10}), $\frak M_\pm$ is smooth of dimension $4c-\Delta ^2
-3=-p-3$
at $V$ if $V$ is stable and $\Ext^2(V, V) = 0$, since $h^2(X; \scrO
_X) = 0$. Moreover, setting $W = V\ddual$, there is a surjection from $H^2(X;
Hom(W, W))$ to  $\Ext^2(V, V)$. Thus to show that $\Ext^2(V, V) =0$ it suffices
to show that $H^2(X;  Hom(W, W))=0$. Now $H^2(X; Hom(W,
W))$ is dual to $H^0(X; Hom (W,W) \otimes K_X)$. Since $-K_X$ is
effective, there is an inclusion of $H^0(X; Hom (W,W) \otimes K_X)$ in $H^0(X;
Hom (W,W))$. If $W$ is stable, then $H^0(X; Hom (W,W)) \cong \Cee$ and $H^0(X;
Hom (W,W) \otimes K_X) = 0$. Thus
$\frak M_\pm$ is smooth at $V$. Standard theory \cite{1, 10} also shows that
every
torsion free sheaf $V$ for which $V\ddual$ is stable is smoothable. Thus the
set
of locally free sheaves is nonempty and dense in the component containing $V$
in
this case.

Now consider a $V$ such that $W=V\ddual$ is not stable.
Using the exact sequence
$$0 \to \scrO_X(F) \to W \to \scrO_X(F) \otimes I_Z \to 0$$
for $W$ which was given in the course of the proof of (2.1),
it is easy to check that there is an exact sequence
$$0 \to \Hom (I_Z, W\otimes \scrO_X(-F)\otimes K_X) \to \Hom (W,W \otimes
K_X) \to H^0(W\otimes \scrO_X(-F)\otimes K_X).$$
Since $-K_X$ is effective and nonzero, $H^0(W\otimes \scrO_X(-F)\otimes
K_X)=\Hom (I_Z, W\otimes \scrO_X(-F)\otimes K_X) =0$. Thus $\Hom (W,W
\otimes K_X)=0$ as well. Once again $V$ is smoothable.

Now we claim that a general smoothing $V'$ of $V$ is Mumford stable.
For otherwise by the proof of (2.1) there is an exact sequence
$$0 \to \scrO_X(F) \to V' \to \scrO_X(F) \otimes I_Z \to 0$$
as above, with $\ell (Z) \leq \ell (\emptyset) =0$. In this case $V'$
is an extension of $\scrO_X(F)$ by $\scrO_X(F)$, forcing $w=p=0$ and
(since $h^1(\scrO_X)= 0$) $V' = \scrO_X(F) \oplus \scrO_X(F) $.
\endproof

It is natural
to make the following conjecture, which is true for geometrically ruled $X$ by
\cite{29} and is verified in certain other cases by \cite{34}.

\proclaim{ Conjecture 2.4} If $X$ is a rational surface with
$-K_X$ effective, then for every choice of
$L$, $\Delta$ and $c$, ${\frak M}_L(\Delta, c)$ is
either empty or irreducible.
\endstatement

Let us fix some notations for the rest of this paper.

\definition{Definition 2.5} Let $X$ be an algebraic surface (not necessarily
rational), and let $\zeta$ be a fixed numerical equivalence class defining a
wall of type $(\Delta, c)$. Set
$\ell _\zeta = (4c-\Delta ^2 +\zeta ^2)/4 = (\zeta ^2-p)/4$.
Choose two nonnegative integers $n_-$ and $n_+$ with
$n_-+ n_+ = \ell _\zeta$, and let $E_\zeta^{n_-, n_+}$ be the set of
all isomorphism classes of nonsplit extensions of the form
$$0 \to \scrO _X(F )\otimes I_{Z_-} \to V \to \scrO _X(\Delta -F ) \otimes
I_{Z_+} \to 0$$
with $\zeta \equiv 2F -\Delta$ and $\ell (Z_\pm) = n_\pm$.
\enddefinition
\medskip

We remark that since $\zeta \equiv \Delta \pmod 2$
and $\Delta ^2 - 4c \le \zeta^2 < 0$, $\ell _\zeta$ is a nonnegative integer.
If $V$ corresponds to a point of $E_\zeta^{n_-, n_+}$, then $V$ is
$L_+$-unstable since $L_+\cdot \zeta > 0$. By (2.2)(ii), $V$ is simple, and
if it is $L_-$-semistable then it is actually stable. By (2.3), if
$X$ is  a rational surface with
$-K_X$ effective, then $\frak M_-$ is smooth in a neighborhood of a point
corresponding to a sheaf $V$ lying in $E_\zeta^{n_-, n_+}$ for some $\zeta,
n_-,
n_+$. We shall now study $E_\zeta^{n_-, n_+}$
in  more detail for rational surfaces.

\lemma{2.6} Suppose that $-K_X$ is effective and that $q(X) =0$. For $Z_-$ and
$Z_+$ two fixed zero-dimensional subschemes of
$X$ of lengths $n_-$ and $n_+$ respectively,
$$\dim \Ext^1(\scrO _X(\Delta -F ) \otimes
I_{Z_+} , \scrO _X(F )\otimes I_{Z_-}) =n_-+ n_+ + h(\zeta)= \ell _\zeta +
h(\zeta),$$ where
$$h(\zeta ) = h^1(X; \scrO _X(2F-\Delta)) = \frac{(\zeta \cdot K_X)}{2} -
\frac{\zeta ^2}{2} -1.$$
\endstatement
\proof Note that $\Hom(\scrO _X(\Delta -F ) \otimes
I_{Z_+} , \scrO _X(F )\otimes I_{Z_-})\subseteq H^0(\scrO_X(2F-\Delta)) = 0$,
since $L_- \cdot (2F-\Delta ) < 0$. Likewise $\Ext^2(\scrO _X(\Delta -F )
\otimes
I_{Z_+} , \scrO _X(F )\otimes I_{Z_-})$ is Serre dual to $\Hom (\scrO _X(F
)\otimes I_{Z_-} , \scrO _X(\Delta -F ) \otimes I_{Z_+}\otimes K_X) \subseteq
H^0(\scrO_X(\Delta -2F) \otimes K_X) \subseteq  H^0(\scrO_X(\Delta -2F))$,
since $-K_X$ is effective. Thus as $L_+ \cdot (\Delta - 2F) < 0$, $\Ext^2(\scrO
_X(\Delta -F ) \otimes I_{Z_+} , \scrO _X(F )\otimes I_{Z_-})=0$ as well. If we
set $\chi (\scrO _X(\Delta -F ) \otimes
I_{Z_+} , \scrO _X(F )\otimes I_{Z_-}) = \sum _i (-1)^i\dim \Ext^i(\scrO
_X(\Delta -F ) \otimes I_{Z_+} , \scrO _X(F )\otimes I_{Z_-})$, then $\chi
(\scrO _X(\Delta -F ) \otimes I_{Z_+} , \scrO _X(F )\otimes I_{Z_-}) =
 -\dim \Ext^1(\scrO _X(\Delta -F ) \otimes I_{Z_+} , \scrO _X(F )\otimes
I_{Z_-})$. Now a standard argument \cite{27} shows that
$$\gather
\chi (\scrO_X(\Delta -F ) \otimes I_{Z_+} , \scrO _X(F )\otimes I_{Z_-})\\
= \int_X\ch(\scrO _X(\Delta -F ) \otimes I_{Z_+})\spcheck\cdot
\ch(\scrO _X(F)\otimes I_{Z_-})\cdot \Todd_X.
\endgather$$
Here given a class $a = \sum a_i \in \bigoplus _iA^i(X)$, we denote by
$a\spcheck$ the class $\sum _i (-1)^ia_i$. An easy computation gives
$$\gather
\int _X\ch(\scrO_X(\Delta -F ) \otimes I_{Z_+})\spcheck\cdot
\ch(\scrO _X(F )\otimes I_{Z_-})\cdot \Todd_X\\
=\int _X\ch(\scrO_X(\Delta -F )\spcheck\cdot \ch(\scrO
_X(F )\cdot\Todd _X - \ell (Z_-) - \ell(Z_+).
\endgather$$
Reversing the above argument, we see that
$$\align
\int _X\ch(\scrO_X(\Delta -F )\spcheck\cdot \ch(\scrO
_X(F )\cdot\Todd _X &= \chi (\scrO_X(2F-\Delta))\\
  =- h^1(X; \scrO _X(2F-\Delta)) &=
\frac{\zeta ^2}{2}-\frac{(\zeta \cdot K_X)}{2} +1 = -h(\zeta ).
\endalign$$
Putting these together we see that $\dim \Ext^1(\scrO _X(\Delta -F ) \otimes
I_{Z_+} , \scrO _X(F )\otimes I_{Z_-})$ is equal to $n_-+ n_+ + h(\zeta)$.
\endproof

Let us describe the scheme structure on $E_\zeta^{n_-, n_+}$ more carefully.
For $Z_-$ and $Z_+$ fixed, the set of extensions in $E_\zeta^{n_-, n_+}$
corresponding to $Z_-$, $Z_+$, is equal to $\Pee \Ext^1(\scrO _X(\Delta -F )
\otimes I_{Z_+} , \scrO _X(F )\otimes I_{Z_-})$. To make a universal
construction, let $H_{n_\pm} = \Hilb ^{n_\pm}X$. Let $\Cal Z_{n_\pm}$ be the
universal codimension two subscheme of $X\times H_{n_\pm}$. Let $\pi _1, \pi
_2$
be the projections of $X\times H_{n_-}\times H_{n_+}$ to $X$, $H_{n_-}\times
H_{n_+}$ respectively, and let $\pi _{1,2}$, $\pi _{1,3}$ be the projections of
$X\times H_{n_-}\times H_{n_+}$ to $X\times H_{n_-}$, $X\times H_{n_+}$
respectively. Define
$$\Cal E_\zeta ^{n_-, n_+} = Ext ^1_{\pi _2}(\pi _1^*\scrO _X(\Delta -F )
\otimes \pi _{1,3}^*I_{\Cal Z_{n_+}}, \pi _1^*\scrO _X(F )\otimes \pi
_{1,2}^*I_{\Cal Z_{n_-}}).$$
The previous lemma and standard base change results show that $\Cal E_\zeta
^{n_-, n_+}$ is locally free of rank $h(\zeta)+\ell _\zeta$ over $H_{n_-}\times
H_{n_+}$. We set $E_\zeta^{n_-, n_+} = \Pee ((\Cal E_\zeta ^{n_-,
n_+})\spcheck)$,  if $h(\zeta)+\ell _\zeta > 0$. Moreover by standard facts
about  relative Ext sheaves there is an exact sequence
$$\gather
0 \to R^1\pi _2{}_*Hom\Big(\pi _1^*\scrO _X(\Delta -F )
\otimes \pi _{1,3}^*I_{\Cal Z_{n_+}}, \pi _1^*\scrO _X(F )\otimes \pi
_{1,2}^*I_{\Cal Z_{n_-}}\Big) \to \Cal E_\zeta
^{n_-, n_+} \to \\
\to \pi _2{}_*Ext ^1\Big(\pi _1^*\scrO _X(\Delta -F )
\otimes \pi _{1,3}^*I_{\Cal Z_{n_+}}, \pi _1^*\scrO _X(F )\otimes \pi
_{1,2}^*I_{\Cal Z_{n_-}}\Big) \to 0.
\endgather$$

\corollary{2.7} With $X$ as in \rom{(2.6)}, if $h(\zeta) +\ell _\zeta= h^1(X;
\scrO _X(2F-\Delta)+\ell_\zeta
\neq 0$,
$E_\zeta^{n_-, n_+}$ is a $\Pee ^{N_\zeta}$-bundle over
$H_{n_-}\times H_{n_+}$, where $N_\zeta = \dim \Ext^1 - 1
= h(\zeta)+\ell _\zeta -1$. Thus if $h(\zeta)+\ell _\zeta\neq 0$, then $\dim
E_\zeta^{n_-, n_+} = 3\ell _\zeta + h(\zeta)-1$. Moreover in this case
$E_{-\zeta}^{n_+, n_-}$ is a $\Pee ^{N_{-\zeta}}$-bundle over
$H_{n_+}\times H_{n_-}$, and $N_\zeta + N_{-\zeta} + 2\ell _\zeta = -p-4$.  If
$h(\zeta) +\ell _\zeta =0$, then $E_\zeta^{0, 0} = \emptyset$ and
$E_{-\zeta}^{0, 0} =
\Pee ^{-p-3}$ is a component of $\frak M_+$. Finally this last case arises if
and only if
$\zeta ^2 = p$ and $\zeta \cdot K_X = \zeta ^2+2 = p+2$.
\endstatement
\proof Note that $N_\zeta \geq 0$ unless $h(\zeta) +\ell _\zeta =0$. Under this
assumption, we have
$$N_\zeta + N_{-\zeta} + 2\ell _\zeta = 4\ell _\zeta - \zeta ^2 -4= -p-4.$$
The case where $h(\zeta) +\ell _\zeta =0$ is similar.
Moreover if $h(\zeta) +\ell_\zeta =0$, then it follows from (2.2)(ii) that
all of the sheaves $V$ corresponding to points of $E_{-\zeta}^{0, 0}$
are $L_+$-stable. By (2.2)(i) the map $E_{-\zeta}^{0, 0}\to \frak M_+$
is one-to-one. Since $\frak M_+$ is of dimension $-p-3$ and smooth at
points corresponding to the sheaves in $E_{-\zeta}^{0, 0}\to \frak M_+$,
the map $E_{-\zeta}^{0, 0}\to \frak M_+$ must be
an embedding onto a component of $\frak M_+$.
The final statement follows from the formulas
$\zeta ^2 = 4\ell_\zeta +p$ and $\dsize h(\zeta)=
\frac{(\zeta \cdot K_X)}{2} -\frac{\zeta ^2}{2} -1$.
\endproof

If $h(\zeta)+\ell_\zeta\neq 0$, then by Lemma 2.2 there is a rational map from
$E_\zeta^{n_-, n_+}$ to the moduli space $\frak M_-$ which is birational onto
its image. However this map will not in general be a morphism if $n_->0$
(see \cite{16}). We shall study this more carefully in the next sections.

Let us also remark that standard theory gives a universal sheaf $\Cal V$ over
$E_\zeta^{n_-, n_+}$:

\proposition{2.8} Let $\rho \: X\times E_\zeta^{n_-, n_+} \to
X\times H_{n_-}\times H_{n_+}$ be the natural projection, and let $\pi_2\:
X\times E_\zeta^{n_-, n_+} \to E_\zeta^{n_-, n_+}$ be the projection. Then
there
is a coherent sheaf $\Cal V$ over
$X\times E_\zeta^{n_-, n_+}$ and an exact sequence
$$\gather
0 \to
\rho ^*\left(\pi _1^*\scrO _X(F )\otimes \pi
_{1,2}^*I_{\Cal Z_{n_-}}\right)\otimes \pi_2^*\scrO_{E_\zeta^{n_-, n_+}}(1)\\
\to \Cal V \to \rho ^*\left(\pi _1^*\scrO _X(\Delta -F )
\otimes \pi _{1,3}^*I_{\Cal Z_{n_+}}\right) \to  0.\qed
\endgather$$
\endstatement

\noindent {\bf Remark 2.9.} Very similar results hold in the case where $-K_X$
is effective and nonzero (corresponding to certain elliptic ruled surfaces) or
$K_X=0$ (corresponding to $K3$ or abelian surfaces). For example, in the case
of
a $K3$ surface $X$, the moduli space is smooth of dimension $-p-6$ away from
the
sheaves which are strictly semistable for every ample divisor (although there
exist components consisting entirely of non-locally free sheaves for small
values of $-p$). In this case however $h(\zeta ) = -\zeta ^2/2 -2$ and $N_\zeta
+ N_{-\zeta} + 2\ell _\zeta = -p-6$, which is equal to the dimension $d$ of the
moduli space instead of to $d-1$. For example, if $\ell _\zeta = 0$, then
$N_\zeta = N_{-\zeta} = d/2$. In this case $E_\zeta ^{0,0} \cong \Pee ^{d/2}$
is
a maximal isotropic submanifold of the symplectic manifold $\frak M_-$. In
other words, the natural holomorphic $2$-form $\omega$ on $\frak M_-$ vanishes
on $E_\zeta ^{0,0}$ and identifies the normal bundle of $E_\zeta ^{0,0}$ in
$\frak M_-$ with the cotangent bundle of $E_\zeta ^{0,0}$.

\section{3. Flips of moduli spaces.}

In this section, we begin by assuming again that $X$ is an arbitrary algebraic
surface. Let $\zeta = \zeta _1, \dots, \zeta _n$ be the positive
rational multiples of $\zeta$ such that $\zeta_i$ is an integral class
also defining the wall $W^\zeta$. Our goal in this section is to deal
with the problem that
there is only a rational map in general from $E_{\zeta_i}^{n_-, n_+}$ to
$\frak M_-$. We shall do so by finding a sequence of spaces between
$\frak M_-$ and $\frak M_+$, each one given by blowing up and down
the previous one, such that for an appropriate member of the sequence
the rational map $E_{\zeta_i}^{n_-, n_+}\dasharrow \frak M_-$
becomes a morphism (and a smooth embedding in the case of rational surfaces).
Throughout the rest of this paper, $L_0$ shall denote any ample divisor
contained in the interior of  the intersection of $W^\zeta$ and
the closures of $\Cal C_\pm$. Recall that we have defined
universal semistability after the proof of (2.1).

\definition{Definition 3.1} Let $k$ be an integer. A rank two torsion free
sheaf $V$ with $c_1(V) = \Delta$ and $\Delta ^2 - 4c_2(V) = p$ is
{\sl $(L_0, \zeta, k)$-semistable\/} if $V$ is Mumford
$L_0$-semistable and if it is strictly Mumford semistable, then either it is
universally semistable or, for all divisors
$F$ such that $2F-\Delta \equiv \zeta$, we have the following:
\roster
\item"{(i)}" If there exists an exact sequence
$$0 \to \scrO _X(F )\otimes I_{Z_1} \to V \to \scrO _X(\Delta -F ) \otimes
I_{Z_2} \to 0,$$
then $\ell (Z_2) \leq k$ and thus $\ell (Z_1) \geq \ell _\zeta - k$.
\item"{(ii)}" If there exists an exact sequence
$$0 \to \scrO _X(\Delta - F )\otimes I_{Z_1} \to V \to \scrO _X(F )
\otimes I_{Z_2} \to 0,$$
then $\ell (Z_1) \geq k+1$ and thus $\ell (Z_2) \leq \ell _\zeta - k-1$.
\endroster
Likewise, setting $\boldsymbol \zeta = (\zeta _1, \dots,\zeta _n)$ and $\bold k
= (k_1, \dots, k_n)$, we say that $V$ is {\sl $(L_0, \boldsymbol \zeta,
\bold k)$-semistable\/} if
$V$ is $(L_0, \zeta _i, k_i)$-semistable for every $i$. Let $\frak
M_0^{(\boldsymbol \zeta, \bold k)}$ denote the set of isomorphism classes of
$(L_0, \boldsymbol \zeta, \bold k)$-semistable rank two sheaves $V$ with
$c_1(V)
= \Delta$ and $\Delta ^2 - 4c_2(V) = p$.
\enddefinition

Next we give some easy properties of $(L_0, \boldsymbol \zeta, \bold
k)$-semistability.

\lemma{3.2}
\roster
\item"{(i)}" If $k_i \geq \ell _{\zeta_i}$ for all $i$, and $V$ is not
universally semistable, then $V$ is
$(L_0, \boldsymbol \zeta, \bold k)$-semistable if and only if it is
$L_-$-stable. Likewise if $k_i \leq -1$ for all $i$ and $V$ is not
universally semistable, then $V$ is
$(L_0, \boldsymbol \zeta, \bold k)$-semistable
if and only if it is$L_+$-stable.

\item"{(ii)}" If $k_i \geq \ell _{\zeta_i}$ for all $i$,
then $\frak M_0^{(\boldsymbol \zeta, \bold k)} = \frak M_-$.
Likewise if $k_i \leq -1$ for all $i$,
then $\frak M_0^{(\boldsymbol \zeta, \bold k)} = \frak M_+$.

\item"{(iii)}" For $n_2 > k_i$, $\frak M_0 ^{(\boldsymbol \zeta, \bold k)}
\cap E_{\zeta_i}^{n_1, n_2}= \emptyset$.

\item"{(iv)}" There is an injection $E_{\zeta_i}^{\ell_{\zeta_i}-k_i,
k_i} \to \frak M_0 ^{(\boldsymbol \zeta, \bold k)}$.
Likewise there is an injection $E_{-\zeta_i}^{k_i+1, \ell _{\zeta_i} - k_i-1}
\to \frak M_0 ^{(\boldsymbol \zeta, \bold k)}$. Finally,
the images of $E_{\zeta_i}^{\ell _{\zeta_i} -k_i, k_i}$ and
$E_{\zeta_j}^{\ell _{\zeta_j}-k_j, k_j}$ are disjoint if $i\neq j$.
\endroster
\endstatement
\proof If $k_i\geq \ell _{\zeta_i}$ for all $i$, then the condition that
$\ell(Z_2) \leq \ell _{\zeta_i}$ and $\ell (Z_1) \geq 0$ are trivially always
satisfied and the conditions
$\ell(Z_2) \leq -1$ and $\ell (Z_1) \geq \ell _{\zeta_i} +1$ are vacuous. A
similar argument handles the case $k_i \leq -1$ for all $i$. It is easy to see
that this implies (i).  Statement (ii) follows from (i), and
(iii) follows from the definitions. As for (iv), let
$V\in  E_{\zeta_i}^{\ell _{\zeta_i} -k_i, k_i}$. To decide if  $V$ is in $\frak
M_0 ^{(\boldsymbol \zeta, \bold k)}$, we look for potentially destabilizing
subsheaves with torsion free quotient. Similar arguments as in \cite{30}
show that the only potentially destabilizing subsheaves with
torsion free quotient must be either $\scrO _X(F )\otimes I_{Z_1}$
or $\scrO _X(\Delta - F )\otimes I_{Z}$. By hypothesis, there is a
unique subsheaf of $V$ of the form
$\scrO _X(F )\otimes I_{Z_1}$, and it is not destabilizing.
If there is a subsheaf of the form $\scrO _X(\Delta - F )\otimes I_{Z}$
with torsion free quotient, then by Lemma 2.2 we have
$\ell (Z) > \ell (Z_2) = k_i$ and so $\ell (Z) \geq  k_i + 1$.
Hence such a subsheaf is also not destabilizing.
Thus by Definition 3.1 $V$ is $(L_0, \boldsymbol \zeta, \bold k)$-semistable.
The fact that the map $E_{\zeta_i}^{\ell _{\zeta_i} -k_i, k_i} \to
\frak M_0 ^{(\boldsymbol \zeta, \bold k)}$ is one-to-one and that
$E_{\zeta_i}^{\ell _{\zeta_i} -k_i, k_i}$ and
$E_{\zeta_j}^{\ell _{\zeta_j} -k_j, k_j}$ are disjoint if $i\neq j$ also
follow from similar arguments in \cite{30}. The statement about
$E_{-\zeta_i}^{k_i+1, \ell _{\zeta_i} - k_i-1}$ is similar.
\endproof

Next suppose that we are given two integral vectors $\bold k$ and $\bold k'$
and
a subset $I$ of $\{1 \dots, n\}$ such that $k_i' = k_i$ if $i\notin I$ and
$k_i'
= k_i-1$ if $i\in I$. We investigate the change as we pass from
$\frak M_0 ^{(\boldsymbol
\zeta,
\bold k)}$ to $\frak M_0 ^{(\boldsymbol \zeta, \bold k')}$.

\lemma{3.3} The set of sheaves $V$ in $\frak M_0 ^{(\boldsymbol \zeta, \bold
k)}$
which are not $(L_0, \boldsymbol \zeta, \bold k')$-semistable is exactly the
image
of
$\bigcup _{i\in I}E_{\zeta_i}^{\ell _{\zeta_i}-k_i, k_i}$. Likewise the set of
$V\in
\frak M_0 ^{(\boldsymbol \zeta, \bold k')}$ which are not
$(L_0, \boldsymbol \zeta, \bold k)$-semistable is exactly the image of
$\bigcup _{i\in I}E_{-\zeta_i}^{k_i,\ell _{\zeta_i} -k_i}$.
\endstatement
\proof If $V$ is $(L_0, \boldsymbol \zeta, \bold k)$-semistable but not $(L_0,
\boldsymbol \zeta, \bold k')$-semistable, then $V$ must be Mumford strictly
$L_0$-semistable. Suppose that the
$(L_0, \boldsymbol \zeta, \bold k')$-destabilizing subsheaf is of the form
$\scrO
_X(F )\otimes I_{Z_1}$, where $F$ corresponds to $\zeta _i$ for some $i\in I$.
Then
$\ell (Z_2)
\leq k_i$ (since $V\in
\frak M_0 ^{(\boldsymbol \zeta, \bold k)}$) but
$\ell (Z_2) \geq k_i$ (since the subsheaf is $(L_0, \boldsymbol \zeta, \bold
k')$-destabilizing, for $k_i' = k_i-1$) so that $\ell (Z_2) = k_i$. Thus
$V\in E_\zeta^{\ell _{\zeta_i} -k_i, k_i}$. The other possibility is that the
destabilizing subsheaf is of the form $\scrO _X(\Delta - F )\otimes I_{Z_1}$.
Here we need $\ell (Z_1) \geq k_i+1$ but $\ell (Z_1) < k_i$ and there are no
such
sheaves. The statement about $\frak M_0 ^{(\boldsymbol \zeta, \bold k')}$
follows
by symmetry.
\endproof

We shall now describe a sequence of actual moduli spaces $\frak M_0
^{(\boldsymbol \zeta, \bold k)}$ for which the integral vector $\bold k$
change in the way described before the statement of (3.3).

\definition{Definition 3.4} Suppose that $\zeta _i = r_i\zeta _1$, where $r_i$
is a positive rational number. Given $t\in \Bbb Q$, let $t_i = r_it$, so that
$t_1=t$.  Suppose that $\dsize
\frac{\ell _{\zeta _i} + t_i}2$ is not an integer for any
$i$. In this case, define
$$k_i(t) = \fracwithdelims[]{\ell _{\zeta _i} + t_i}{2},$$
where $[x]$ is the greatest integer function, and define $\bold k(t)$ to be
the vector formed by the
$k_i(t)$. A rational number $t$ is {\sl $\zeta _i$-critical\/} if
$\dsize \frac{\ell_{\zeta _i} + t_i}2 \in \Zee$ and $-1 \leq \dsize \frac{\ell
_{\zeta _i} + t_i}2 \leq \ell _{\zeta _i}$. We shall also say that $t_i$ is
{\sl $\zeta _i$-critical\/}. Finally $t$ is {\sl
$\boldsymbol
\zeta$-critical\/} if it is $\zeta _i$-critical for some $i$. Note that there
are
only finitely many such $t$.
\enddefinition
\medskip

Given $t \in \Bbb Q$, let $I(t) = \{\, i: t {\text{ is $\zeta
_i$-critical}}\,\}$. Suppose that $\varepsilon$ is chosen so that, for every
$i$, either there is no $\zeta _i$-critical rational number in $[t_i -
r_i\varepsilon, t_i +
r_i\varepsilon]$ or
$t_i$ is the unique
$\zeta _i$-critical rational number in $[t_i - r_i\varepsilon, t_i +
r_i\varepsilon]$. Equivalently either there is no $\boldsymbol
\zeta$-critical number in $[t-\varepsilon, t+\varepsilon]$ or $t$ is the unique
$\boldsymbol
\zeta$-critical number in $[t-\varepsilon, t+\varepsilon]$. Then we clearly
have:
$$k_i(t - \varepsilon) = \cases k_i(t+
\varepsilon) , &\text{if
$i\notin I(t)$}\\
k_i(t + \varepsilon) -1, &\text{if $i\in I(t)$.}
\endcases$$
In particular  if there is no
$\boldsymbol
\zeta$-critical number in $[t-\varepsilon, t+\varepsilon]$, so that $I(t)=
\emptyset$, then
$k_i(t -
\varepsilon) = k_i(t + \varepsilon)$ for every $i$. Further note that if $t\gg
0$, then $k_i(t) > \ell _{\zeta _i}$ for every $i$, and if $t\ll 0$, then
$k_i(t) <-1$  for every
$i$.

We then have the following theorem, whose proof will be given in the next
section:

\theorem{3.5} For all $t\in \Bbb Q$ which are not $\boldsymbol
\zeta$-critical, there exists a natural structure of a projective scheme on
$\frak
M_0 ^{(\boldsymbol
\zeta, \bold k(t))}$ for which it is a coarse moduli space.
\endstatement
\medskip
The proof of (3.5) will also show that $\frak M_0 ^{(\boldsymbol
\zeta, \bold k(t))}$ has the usual properties of a coarse moduli space: all
sheaves corresponding to points of $\frak M_0 ^{(\boldsymbol
\zeta, \bold k(t))}$ will turn out to be simple (as they will turn out to be
stable for an appropriate notion of stability), a classical or formal
neighborhood of a point of
$\frak M_0 ^{(\boldsymbol
\zeta, \bold k(t))}$ may be identified with the universal deformation space of
the corresponding sheaf, and there exists a universal sheaf locally in the
classical or
\'etale topology around every point of $\frak M_0 ^{(\boldsymbol
\zeta, \bold k(t))}$.

For the rest of this section, we shall again restrict to the case where $X$ is
a
rational surface with $-K_X$ effective, unless otherwise noted. Let $\zeta =
\zeta _i$ for some $i$ and let $\bold k=\bold k(t)$ for some $t$ which
is not
$\boldsymbol \zeta$-critical. The first step is to make some infinitesimal
calculations concerning the differential of the map
$E_\zeta^{\ell _\zeta -k, k}\to
\frak M_0 ^{(\boldsymbol \zeta, \bold k)}$ and the normal bundle to its image.

\proposition{3.6} The map $E_\zeta^{\ell _\zeta -k, k}\to
\frak M_0 ^{(\boldsymbol \zeta, \bold k)}$ is an immersion. The normal bundle
$\Cal N_\zeta ^{\ell _\zeta -k, k}$ to $E_\zeta^{\ell _\zeta -k, k}$ in $\frak
M_0 ^{(\boldsymbol \zeta, \bold k)}$ is exactly $\rho ^*\Cal E _{-\zeta}^{k,
\ell _\zeta -k}\otimes
\scrO_{E_\zeta ^{\ell _\zeta -k, k}}(-1)$, in the notation of the previous
section.
\endstatement
\proof Since every sheaf in $\frak M_0 ^{(\boldsymbol \zeta, \bold k)}$ is
actually stable and therefore simple (which was also proved in (2.2)) we may
identify an analytic neighborhood of $V\in \frak M_0 ^{(\boldsymbol \zeta,
\bold
k)}$ with the germ of the universal deformation space for $V$, i\.e\. with
$\Ext
^1(V,V)$. Let us now calculate the tangent space to
$E_\zeta ^{\ell _\zeta -k, k}$ at $V$: suppose that $\xi \in
\Ext^1(\scrO _X(\Delta -F ) \otimes I_{Z_2} , \scrO _X(F)\otimes I_{Z_1})
=\Ext^1$ is a nonzero extension class corresponding to $V$, where $\ell (Z_1) =
\ell _\zeta -k$ and $\ell (Z_2) = k$.  Let
$H_{\ell _\zeta -k} = \Hilb ^{\ell _\zeta -k}X$ and $H_k =
\Hilb ^kX$. Then there is the following exact sequence for the tangent
space to $E_\zeta ^{\ell _\zeta -k, k}$ at $\xi$:
$$0\to \Ext ^1/\Cee \cdot \xi \to T_\xi E_\zeta ^{\ell _\zeta -k, k}\to
T_{Z_1}H_{\ell _\zeta -k}
\oplus T_{Z_2}H_k \to 0.$$
Note further that the tangent space to $\Hilb ^nX$ at $Z$ is
equal to $\Hom(I_Z,\scrO_Z)$, which we may further canonically identify with
$\Ext^1(I_Z, I_Z)$ since $X$ is rational and by a local calculation.  We then
have the following:

\proposition{3.7} For all nonzero $\xi \in \Ext^1$, the natural map from a
neighborhood of
$\xi$ in $E_\zeta ^{\ell _\zeta -k, k}$ to $\frak M_0 ^{(\boldsymbol \zeta,
\bold
k)}$ is an immersion at
$\xi$. The image of $T_\xi E_\zeta ^{\ell _\zeta -k, k}$ in $\Ext ^1(V,V)$ is
exactly the kernel of the natural map $\Ext ^1(V,V)\to \Ext ^1(\scrO
_X(F)\otimes I_{Z_1},
\scrO _X(\Delta -F ) \otimes I_{Z_2})$, and the normal space to $E_\zeta
^{\ell _\zeta -k, k}$ at $\xi$ in $\frak M_0 ^{(\boldsymbol \zeta, \bold k)}$
may
be canonically identified with
$\Ext ^1 (\scrO _X(F)\otimes I_{Z_1},  \scrO _X(\Delta -F ) \otimes I_{Z_2})$.
\endstatement
\noindent
{\it Proof.} Consider the natural map from
$\Ext ^1(V,V)$ to
$\Ext ^1 (\scrO _X(F)\otimes I_{Z_1},  \scrO _X(\Delta -F ) \otimes I_{Z_2})$.
We claim that this map is surjective and will describe its kernel in more
detail. The map factors into two maps:
$$\gather
\Ext ^1(V,V) \to \Ext^1(V, \scrO _X(\Delta -F ) \otimes I_{Z_2}) \\
\Ext^1(V, \scrO _X(\Delta -F ) \otimes I_{Z_2}) \to \Ext ^1 (\scrO
_X(F)\otimes I_{Z_1},  \scrO _X(\Delta -F ) \otimes I_{Z_2}).
\endgather$$
The cokernel of the first map is contained in $\Ext^2(V, \scrO _X(F)\otimes
I_{Z_1})$. To see that this group is zero, apply Serre duality: it suffices to
show that $\Hom(\scrO _X(F)\otimes I_{Z_1}, V\otimes K_X) =0$. From the
defining exact sequence for $V$, we have an exact sequence
$$\gather
0 \to  \Hom(\scrO _X(F)\otimes I_{Z_1},\scrO _X(F)\otimes I_{Z_1}\otimes
K_X)
\to \Hom(\scrO _X(F)\otimes I_{Z_1}, V\otimes K_X)\\
 \to \Hom(\scrO _X(F)\otimes
I_{Z_1},\scrO _X(\Delta -F ) \otimes I_{Z_2}).
\endgather$$
The first term is just $H^0(K_X) = 0$ and the third is contained in
$H^0(\scrO_X(\Delta -2F) \otimes K_X) =0$. Thus $\Hom(\scrO _X(F)\otimes
I_{Z_1}, V\otimes K_X) =0$. The vanishing of the
cokernel of the second map, namely $\Ext ^2(\scrO _X(\Delta -F ) \otimes
I_{Z_2}, \scrO _X(\Delta -F )
\otimes I_{Z_2})$, is similar. Thus $\Ext ^1(V,V)\to \Ext ^1 (\scrO
_X(F)\otimes
I_{Z_1},  \scrO _X(\Delta -F ) \otimes I_{Z_2})$ is onto. If $K$ is the kernel,
then arguments as above show that there is an exact sequence
$$0 \to \Ext^1(V, \scrO _X(F)\otimes I_{Z_1}) \to K \to \Ext^1(\scrO _X(\Delta
 -F ) \otimes I_{Z_2}, \scrO _X(\Delta -F ) \otimes I_{Z_2})\to 0.$$
Here $\Ext^1(\scrO _X(\Delta
 -F ) \otimes I_{Z_2}, \scrO _X(\Delta -F ) \otimes I_{Z_2}) = \Ext^1(I_{Z_2},
I_{Z_2})$ is the tangent space to $H_k$. Moreover, there is an exact
sequence
$$\gather
\Hom(\scrO _X(F)\otimes I_{Z_1}, \scrO _X(F)\otimes I_{Z_1}) \to
\Ext^1(\scrO _X(\Delta -F ) \otimes I_{Z_2},\scrO _X(F)\otimes I_{Z_1}) \to \\
\to \Ext^1(V, \scrO _X(F)\otimes I_{Z_1}) \to \Ext^1(\scrO _X(F)\otimes
I_{Z_1},
\scrO _X(F)\otimes I_{Z_1})\to 0.
\endgather$$
The last term is $\Ext^1(I_{Z_1}, I_{Z_1})$ which is the tangent space to
$H_{\ell _\zeta -k}$ at $Z_1$, and the first two terms combine to give
$\Ext^1/\Cee \cdot \xi$. Thus the kernel $K$ looks very much like
the tangent space to $E_\zeta ^{\ell _\zeta -k, k}$ at $\xi$ and both spaces
have the same dimension.

Let us describe
the tangent space to $E_\zeta ^{\ell _\zeta -k, k}$ at $\xi$ and the
differential of the map $E_\zeta ^{\ell _\zeta -k, k}$ to $\frak M_0
^{(\boldsymbol \zeta, \bold k)}$ in more intrinsic terms. It is easy to see
that
a $\Spec
\Cee[\epsilon]$-valued point of $E_\zeta ^{\ell _\zeta -k, k}$ which
restricts to $\xi$ defines two codimension two subschemes $\Cal
Z_1\subseteq X\times
\Spec\Cee[\epsilon]$,
$\Cal Z_2\subseteq X\times \Spec\Cee[\epsilon]$, flat over
$\Spec\Cee[\epsilon]$, restricting to $Z_i$ over $X$, and an extension
$\Cal V$ over $X\times \Spec\Cee[\epsilon]$ of the form
$$0 \to \pi _1^* \scrO _X(F)\otimes I_{\Cal Z_1} \otimes \to \Cal V \to
\pi _1^*\scrO _X(\Delta -F ) \otimes I_{\Cal Z_2}\to 0.$$
Conversely such a choice of $\Cal Z_1$, $\Cal Z_2$ and $\Cal V$ define
a $\Spec \Cee[\epsilon]$-valued point of $E_\zeta ^{\ell _\zeta -k, k}$.
Thus there is a commutative diagram with exact rows and columns:
$$\CD
@. 0 @. 0 @. 0 @.\\
@. @VVV @VVV @VVV @. \\
0 @>>> \scrO _X(F)\otimes I_{Z_1} @>>> V @>>> \scrO _X(\Delta -F ) \otimes
I_{Z_2} @>>> 0\\
@. @VVV @VVV @VVV @. \\
0 @>>> \pi _1^* \scrO _X(F)\otimes I_{\Cal Z_1} @>>> \Cal V @>>> \pi _1^*\scrO
_X(\Delta -F ) \otimes I_{\Cal Z_2}@>>> 0\\
@. @VVV @VVV @VVV @. \\
0 @>>> \scrO _X(F)\otimes I_{Z_1} @>>> V @>>> \scrO _X(\Delta -F ) \otimes
I_{Z_2} @>>> 0\\
@. @VVV @VVV @VVV @. \\
@. 0 @. 0 @. 0 @.
\endCD$$
Here the extension $\Cal V$ of $V$ by $V$, viewed as a point of $\Ext
^1(V,V)$, corresponds to the Kodaira-Spencer map of the deformation $\Cal V$
of $V$. Likewise the left and right hand columns give classes in $\Ext
^1(I_{Z_1}, I_{Z_1})$ and $\Ext ^1(I_{Z_2}, I_{Z_2})$ corresponding to $\Cal
Z_1$ and $\Cal Z_2$. A straightforward diagram chase shows that if $\Cal V$
fits into this commutative diagram then the image of the extension class $\xi
\in \Ext ^1(V,V)$ corresponding to $\Cal V$ in $\Ext ^1 (\scrO
_X(F)\otimes I_{Z_1},  \scrO _X(\Delta -F ) \otimes I_{Z_2})$ is zero. To see
the converse, that every element in the kernel $K$ of the map $\Ext ^1(V,V) \to
\Ext ^1 (\scrO _X(F)\otimes I_{Z_1},  \scrO _X(\Delta -F ) \otimes I_{Z_2})$ is
the image of a tangent vector to $E_\zeta ^{\ell _\zeta -k, k}$ at $\xi$, use
the arguments above which show that there is a surjection from $K$ to
$$\Ext^1(\scrO _X(\Delta -F ) \otimes I_{Z_2}, \scrO _X(\Delta -F )
\otimes I_{Z_2}) = \Ext ^1(I_{Z_2}, I_{Z_2}).$$
Thus there is an induced extension of $\scrO _X(\Delta -F ) \otimes I_{Z_2}$
by $\scrO _X(\Delta -F ) \otimes I_{Z_2}$, necessarily of
the form $\scrO _X(\Delta -F ) \otimes I_{\Cal Z_2}$,
and a map from $\Cal V$ to
$\scrO _X(\Delta -F ) \otimes I_{\Cal Z_2}$, necessarily a surjection. The
kernel of this surjection then defines an extension $\scrO _X(F)\otimes
I_{\Cal Z_1}$ of $\scrO _X(F)\otimes I_{Z_1}$ by $\scrO _X(F)\otimes I_{Z_1}$.
It follows that $K$ is in the image of the tangent space to
$E_\zeta ^{\ell _\zeta -k, k}$ at $\xi$. By counting dimensions the map
on tangent spaces from
$T_\xi E_\zeta ^{\ell _\zeta -k, k}$ to $\Ext ^1(V,V)$ is injective, showing
that the map from $E_\zeta ^{\ell _\zeta -k, k}$ to $\frak M_0^{(\boldsymbol
\zeta, \bold k)}$ is an immersion and identifying the normal space at $\xi$.
\qed

Let us continue the proof of Proposition 3.6.
To give a global description of the normal bundle to
$E_\zeta ^{\ell _\zeta -k, k}$ in $\frak M_0^{(\boldsymbol \zeta, \bold
k)}$, recall by standard deformation
theory \cite{10} that the pullback of the tangent bundle of $\frak
M_0^{(\boldsymbol \zeta, \bold k)}$ to $E_\zeta ^{\ell _\zeta -k, k}$ is
just $Ext ^1_{\pi_2}(\Cal V, \Cal V)$,
where $\Cal V$ is the universal sheaf over
$X\times E_\zeta ^{\ell _\zeta -k, k}$ described in (2.8) and $\pi_2\: X \times
E_\zeta ^{\ell _\zeta -k, k} \to E_\zeta ^{\ell _\zeta -k, k}$ is
the second projection. Moreover the
calculations above globalize to show that the normal bundle is exactly
$$Ext ^1_{\pi_2}(\rho ^*\left(\pi _1^*\scrO _X(F )\otimes \pi
_{1,2}^*I_{\Cal Z_1}\right)\otimes \pi_2^*\scrO_{E_\zeta^{\ell _\zeta -k,
k}}(1),
\rho ^*\left(\pi _1^*\scrO _X(\Delta -F )
\otimes \pi _{1,3}^*I_{\Cal Z_2}\right)),$$
where $\rho \: X\times E_\zeta^{\ell _\zeta -k, k} \to
X \times H_{\ell _\zeta -k}\times H_k$ is the
natural projection. Using standard base change results and the projection
formula, we see that this sheaf is equal to
$$\rho ^*Ext ^1_{\pi _2}(\pi _1^*\scrO _X(F )\otimes \pi
_{1,2}^*I_{\Cal Z_1}, \pi _1^*\scrO _X(\Delta -F )
\otimes \pi _{1,3}^*I_{\Cal Z_2})\otimes \scrO_{E_\zeta^{\ell _\zeta -k,
k}}(-1),$$ which is the same as $\rho ^*\Cal E _{-\zeta}^{k,\ell _\zeta
 -k}\otimes
\scrO_{E_\zeta ^{\ell _\zeta -k, k}}(-1)$.
\endproof

Finally, to compare the moduli space $\frak M_0 ^{(\boldsymbol \zeta, \bold
k(t+\varepsilon))}$ with $\frak M_0 ^{(\boldsymbol \zeta, \bold
k(t-\varepsilon))}$, where $t$ is the unique
$\boldsymbol \zeta$-critical point in $[t-\varepsilon, t +\varepsilon]$, we
shall
need the following result which is a straightforward generalization of (A.2) of
\cite{11}.

\proposition{3.8} Let
$X$ be a smooth projective scheme or compact complex manifold, and let $T$ be
smooth. Suppose that $\Cal
V$ is a rank two reflexive sheaf over $X\times T$, flat over $T$. Let $D$ be
a reduced divisor on $T$, not necessarily smooth and let $i\: D \to T$ be the
inclusion. Suppose that $L$ is a line bundle on $X$ and that $\Cal Z$ is a
codimension two subscheme of $X\times D$, flat over $D$. Suppose further that
$\Cal V \to i_*\pi _1^*L\otimes I_{\Cal Z}$ is a surjection, and let $\Cal V'$
be its kernel:  $$0 \to \Cal V' \to \Cal V \to i_*\pi _1^*L\otimes I_{\Cal Z}
\to 0.$$ Then there is a line bundle $M$ on $X$ and a subscheme
$\Cal Z'$ of $X\times D$ codimension at least two, flat over $D$, with the
following properties:
\roster
\item"{(i)}" $\Cal V'$ is reflexive and flat over $T$.
\item"{(ii)}" There are exact sequences
$$\align
0 \to \pi _1^*M \otimes I_{\Cal Z'} \to &\Cal V|D \to \pi _1^*L\otimes
I_{\Cal Z} \to 0 ;\\
0\to \pi _1^*L\otimes I_{\Cal Z}\otimes \scrO_D(-D) \to &\Cal V'|D \to \pi
_1^*M
\otimes I_{\Cal Z'}
\to 0,
\endalign$$
which restrict for each $t\in D$ to give
exact sequences
$$\align 0 \to M\otimes I_{Z'} \to &V_t \to L\otimes I_Z \to 0;\\ 0\to
L\otimes I_Z \to &(V_t)' \to M\otimes I_{Z'} \to 0.
\endalign$$ Here $Z$ is the subscheme of $X$ defined by $\Cal Z$ for the
slice $X\times
\{t\}$ and $Z'_t$ is likewise defined by $\Cal Z'$.
\item"{(iii)}" If $D$ is smooth, then the extension class corresponding to
$(V_t)'$ in $\Ext ^1(M\otimes I_W,  L\otimes  I_Z)$ is defined by the image of
the normal vector to $D$ at $t$ under the composition of the Kodaira-Spencer
map from the tangent space of $T$ at $t$ to
$\Ext ^1(V_t, V_t)$, followed by the natural map $\Ext ^1(V_t, V_t) \to \Ext
^1(M\otimes I_{Z'},  L\otimes I_Z)$. \endroster
\endstatement
\medskip

Here $\Cal V'$ is called the {\sl elementary modification\/} of $\Cal V$ along
$D$. This construction has the following symmetry: if we make the elementary
modification of $\Cal V'$ along $D$ corresponding to the surjection $\Cal V'
\to i_*\bigl(\pi _1^*M \otimes I_{\Cal Z'}\bigr)$, then the result is $\Cal V
\otimes \scrO_{X\times T}(-(X\times D))$.

Here is the typical way that we will apply the above: given $X$, let $M$ be a
smooth manifold and $Y$ a submanifold of $M$. Let $T$ be the
blowup of $M$ along $Y$ and let $D$ be the exceptional divisor.
Let $\pi \: T \to M$ be the natural map. Then, given $\xi \in D$,
the image in the normal space to $\pi(\xi)$ of the normal direction at $\xi$ to
$D$ under $\pi _*$ may be identified with the line in the normal space
corresponding to $\xi$.

We can now state the main result as follows:

\theorem{3.9} Suppose that $t$ is the unique
$\boldsymbol \zeta$-critical point in $[t-\varepsilon, t +\varepsilon]$. If
$h(\pm\zeta_i) + \ell _{\pm\zeta_i} \neq 0$ for every $i$, then
the rational map
$\frak M_0 ^{(\boldsymbol \zeta, \bold k(t+\varepsilon))}
\dasharrow
\frak M_0 ^{(\boldsymbol \zeta, \bold k(t-\varepsilon))}$ is
obtained as follows. For every $i$, fixing $\zeta _i =\zeta$ and
$k_i(t+\varepsilon) = k$, blow up
$E_\zeta^{\ell _\zeta -k, k}$ in
$\frak M_0 ^{(\boldsymbol \zeta, \bold k(t+\varepsilon))}$. Then
the exceptional divisor $D$ is a $\Pee ^{N_\zeta}\times \Pee
^{N_{-\zeta}}$-bundle over $\Hilb ^{\ell _\zeta -k}X\times
\Hilb ^kX$. Moreover this divisor can be contracted in two different ways.
Contracting the $\Pee ^{N_{-\zeta}}$ fibers for all possible $\zeta$ gives
$\frak
M_0 ^{(\boldsymbol
\zeta, \bold k(t+\varepsilon))}$. Contracting the $\Pee
^{N_{\zeta}}$ fibers for all possible $\zeta$ gives $\frak M_0 ^{(\boldsymbol
\zeta, \bold k(t-\varepsilon))}$. Moreover the morphism from the
blowup to
$\frak M_0 ^{(\boldsymbol \zeta, \bold k(t-\varepsilon))}$ is
induced by an elementary modification as in \rom{(3.8)}, and the image of the
the component of the exceptional divisor which is the blowup of $E_\zeta^{\ell
_\zeta -k, k}$ is
$E_{-\zeta}^{k, \ell _\zeta -k}$. Finally the construction is symmetric.

Similar statements hold if $h(\pm\zeta_i) + \ell _{\pm\zeta_i}=0$ for some $i$,
where we must also add in or delete an extra component coming from $\pm\zeta
_i$.
\endstatement
\proof Begin by blowing up $E_\zeta^{\ell _\zeta -k, k}$ in $\frak M_0
^{(\boldsymbol \zeta, \bold k(t+\varepsilon))}$ for all possible
$\zeta$. For simplicity we shall just write down the argument in case there is
only one
$\zeta$; the general case is just additional notation. Let
$\widetilde{\frak M}_0^{(\boldsymbol \zeta, \bold k(\bold
t+\varepsilon))}$ denote the blowup and $D$ the exceptional divisor.
Note that the normal bundle
$\Cal N_\zeta ^{\ell _\zeta -k, k}$ to $E_\zeta^{\ell _\zeta -k, k}$ in $\frak
M_0 ^{(\boldsymbol \zeta, \bold k(\bold
t+\varepsilon))}$ is $\rho ^*\Cal E _{-\zeta}^{k, \ell _\zeta -k}\otimes
\scrO_{E_\zeta ^{\ell _\zeta -k, k}}(-1)$, where $\rho\: E_\zeta ^{\ell _\zeta
 -k, k} \to \Hilb ^{\ell _\zeta -k}X\times \Hilb ^kX$ is the projection. In
particular $\Cal N_\zeta ^{\ell _\zeta
 -k, k}$ restricts to each fiber $\Pee ^{N_\zeta}$ to a bundle of the form
$\left[\scrO_{\Pee ^{N_\zeta}}^N\right]\otimes \scrO_{\Pee ^{N_\zeta}}(-1)$,
and an easy calculation using (2.7) shows that $N= N_{-\zeta}+1$. It follows
that
the fibers of the induced map from
$D$ to $\Hilb ^{\ell _\zeta -k}X\times \Hilb ^kX$ are naturally $\Pee
^{N_\zeta}\times \Pee ^{N_{-\zeta}}$. Moreover it is easy to see that
$\scrO(D)|\Pee ^{N_\zeta}= \scrO_{\Pee ^{N_\zeta}}(-1)$, using for example the
fact that $\scrO(D)|\Pee ^{N_\zeta}\times \Pee ^{N_{-\zeta}}= \scrO(a,-1)$ for
some $a$ and the fact that
$$\Cal N_\zeta ^{\ell _\zeta
 -k, k}|\Pee ^{N_\zeta} = R^0\pi _1{}_*[\scrO(-D)|\Pee ^{N_\zeta}\times \Pee
^{N_{-\zeta}}]= \left[\scrO_{\Pee ^{N_\zeta}}^{N_{-\zeta}+1}\right]\otimes
\scrO_{\Pee ^{N_\zeta}}(-a).$$

For the rest of the argument, we assume that there exists a universal family
on $X\times \frak M_0^{(\boldsymbol \zeta,
\bold k(\bold t+\varepsilon))}$. Of course, such a family will only exist
locally in the classical or \'etale topology, but this will suffice for the
argument. Let $\Cal U$ be the pullback of the universal family to
$X\times \widetilde{\frak M}_0^{(\boldsymbol \zeta,
\bold k(\bold t+\varepsilon))}$. Locally again we may assume that the
restriction of $\Cal U$ to $X\times D$ is the pullback of the universal
extension $\Cal V$ of (2.8):
$$\gather
0 \to
\rho ^*\left(\pi _1^*\scrO _X(F )\otimes \pi
_{1,2}^*I_{\Cal Z_{n_-}}\right)\otimes \pi_2^*\scrO_{E_\zeta^{n_-, n_+}}(1)\\
\to \Cal V \to \rho ^*\left(\pi _1^*\scrO _X(\Delta -F )
\otimes \pi _{1,3}^*I_{\Cal Z_{n_+}}\right) \to  0.
\endgather$$
Now consider the effect of making an elementary
transformation of  $\Cal U$ on
$X\times \widetilde{\frak M}_0^{(\boldsymbol \zeta,
\bold k(\bold t+\varepsilon))}$ along the divisor $D$, using the morphism from
$\Cal U$ to the pullback of $\rho ^*\left(\pi _1^*\scrO _X(\Delta -F )  \otimes
\pi _{1,3}^*I_{\Cal Z_k}\right)$  given by considering the pullback of the
universal extension.  Applying
(3.8) to the elementary transformation $\Cal U'$, we see that the fiber of
$\Cal U'$ at a point of the fiber $\Pee ^{N_\zeta}\times \Pee
^{N_{-\zeta}}$ lying over a point $(Z_1, Z_2)\in \Hilb ^{\ell _\zeta -k}X\times
\Hilb ^kX$ is given by a nonsplit extension of the form
$$0 \to \scrO_X(\Delta - F)\otimes I_{Z_2} \to U \to \scrO_X(F)\otimes I_{Z_1}
\to 0.$$
Moreover the extension class corresponding to $U$ is given by the projectivized
normal vector in $\Pee ^{N_{-\zeta}}$. Thus it is independent of the first
factor
$\Pee ^{N_\zeta}$ and the set of all possible such classes is parametrized by
the
second factor $\Pee
^{N_{-\zeta}}$. There is then an induced morphism from $\widetilde{\frak
M}_0^{(\boldsymbol \zeta, \bold k(t+\varepsilon))}$ to $\frak
M_0^{(\boldsymbol \zeta, \bold k(t-\varepsilon))}$ and clearly it
has the effect of contracting
$D$ along its first ruling and has the property that the image of $D$ is
exactly $E_{-\zeta}^{k, \ell _\zeta -k}$. We leave the symmetry of the
construction to the reader. This concludes the proof of (3.9).
\endproof

\noindent {\bf Remark 3.10.} In the $K3$ or abelian case, the arguments of this
section show that the rational map $\frak M_0 ^{(\boldsymbol \zeta, \bold
k(\bold
t+\varepsilon))} \dasharrow \frak M_0 ^{(\boldsymbol \zeta,  \bold
k(t-\varepsilon))}$ is a Mukai elementary transformation \cite{26, 28}.
\medskip

We can also use (3.8) to analyze the rational map from $E_\zeta ^{n_-, n_+}$
to $\frak M_-$, in the case where it is not a morphism. For simplicity we shall
only consider the case of
$E_\zeta ^{1,0}$, i\.e\.
$\ell _\zeta = 1$. In this case $Z_- = p\in X$ and
$I_{Z_-} = \frak m_p$ is the maximal ideal sheaf of $p$. Moreover $\Ext
^1(\scrO
_X(\Delta -F ), \scrO _X(F)\otimes
\frak m_p) = H^1(\scrO _X(2F-\Delta)\otimes \frak m_p)$. There is an exact
sequence
$$0 \to H^0(\Cee _p) \to H^1(\scrO _X(2F-\Delta)\otimes \frak m_p) \to
H^1(\scrO _X(2F-\Delta)) \to 0.$$
Moreover, for $p$ fixed, the extensions $V$ corresponding to a split extension
for $V\ddual$ are exactly the kernel of the map from $H^1(\scrO
_X(2F-\Delta)\otimes \frak m_p)$ to $H^1(\scrO _X(2F-\Delta))$, i\.e\. the
image of $H^0(\Cee _p)$. The normal space is thus identified with $H^1(\scrO
_X(2F-\Delta))$. Now if the extension for $V\ddual$ is split, then there is a
map $\scrO _X(\Delta -F )\otimes \frak m_p \to V$ with quotient $\scrO _X(F)$.
This way of realizing $V$ as an extension gives a surjection $\Ext ^1(V,V) \to
\Ext ^1(\scrO _X(\Delta -F )\otimes \frak m_p, \scrO _X(F))$, and we must look
at the image of the normal space $H^1(\scrO _X(2F-\Delta))$ in this extension
group. On the other hand, we have an exact sequence
$$0\to H^1(\scrO _X(2F-\Delta)) \to \Ext ^1(\scrO _X(\Delta -F )\otimes \frak
m_p, \scrO _X(F)) \to H^0(\Cee _p) \to 0$$
coming from the long exact Ext sequence, and it is an easy diagram chase to see
that the induced map $\Ext ^1(\scrO _X(\Delta -F ), \scrO _X(F)\otimes
\frak m_p) \to \Ext ^1(\scrO _X(\Delta -F )\otimes \frak m_p, \scrO _X(F))$
factors through the map $\Ext ^1(\scrO _X(\Delta -F ), \scrO _X(F)\otimes
\frak m_p) \to H^1(\scrO _X(2F-\Delta))$ and that the image is exactly the
natural subgroup $H^1(\scrO _X(2F-\Delta))$ of $\Ext ^1(\scrO _X(\Delta -F
)\otimes \frak m_p, \scrO _X(F))$.

The above has the following geometric interpretation: the locus $U$ in $E_\zeta
^{1,0}$ of $L_-$-unstable sheaves is in fact a section of $E_\zeta ^{1,0}$. If
we blow up this section and then make the elementary transformation, the result
is exactly the set of elements of $E_\zeta ^{0,1}$ corresponding to nonlocally
free sheaves. This set is already a divisor in $E_\zeta ^{0,1}$. There is thus
a morphism from the blowup of $E_\zeta ^{1,0}$ along $U$ to $\frak M_-$ which
is an embedding into $\frak M_-$. Its image $(E_\zeta ^{1,0})'$ in $\frak M_-$
meets $E_\zeta ^{0,1}$ exactly along the divisor in $E_\zeta ^{0,1}$ of
nonlocally free sheaves.

We can now give a picture of the birational map from $\frak M_-$ to $\frak M_+$
in this case. Begin with the subvariety $E_\zeta ^{0,1}$ in $\frak M_-$ and
blow it up. Let $D^{0,1}$ be the exceptional divisor, ruled in two different
ways. As $E_\zeta ^{0,1}$ meets
$(E_\zeta ^{1,0})'$ along a divisor, the proper transform of $(E_\zeta
^{1,0})'$
in the blowup is again
$(E_\zeta ^{1,0})'$. Making the elementary modification along $D^{0,1}$, we
then
blow down $D^{0,1}$ to get a new moduli space. This moduli space then contains
$E_\zeta ^{1,0}$. At this point we can then blow  up $E_\zeta ^{1,0}$ and
contract the new exceptional divisor $D^{1,0}$ to obtain $\frak M_+$ (a few
extra details need to be checked here concerning the Kodaira-Spencer class).
Note again the symmetry of the situation. In principle we could hope to carry
through this analysis to the case where $\ell _\zeta >1$ as well, but we run
into trouble with the birational geometry of $\Hilb ^nX$. Somehow the
construction of our auxiliary sequence of moduli spaces has eliminated the
necessity for understanding this birational geometry in detail.

\section{4. Mixed stability and mixed moduli spaces.}

Our goal in this section is to give a proof of Theorem 3.5 (for an arbitrary
algebraic surface $X$). By way of motivation for our construction, let us
analyze Gieseker semistability more closely. In the notation of the last
section, we suppose that $L_0$ is an ample line bundle lying on a unique wall
$W$ of type $(w,p)$, and let $\zeta _1, \dots, \zeta _n$ be the integral
classes
of type $(w,p)$ defining $W$. Let $V$ be an
$L_0$-semistable rank two sheaf. Thus either $V$ is Mumford $L_0$-stable or it
is Mumford strictly semistable. In the second case, let $\scrO_X(F)\otimes
I_{Z_1}$ be a destabilizing subsheaf and suppose that there is an exact
sequence
$$0 \to \scrO_X(F)\otimes I_{Z_1} \to V \to \scrO_X(\Delta -F)\otimes
I_{Z_2}\to
0.$$
Let $\zeta = 2F-\Delta$. We shall assume that $\zeta =\zeta _i$ for some $i$,
or equivalently that $\zeta$ is not numerically equivalent to zero (i.e\., $V$
is not universally semistable). By assumption
$\mu _{L_0}(V)
\geq
\mu _{L_0}(\scrO_X(F)\otimes I_{Z_1})$, and so $\chi(V) \geq 2\chi
(\scrO_X(F)\otimes I_{Z_1})$. Since $\chi (V) = \chi (\scrO_X(F)\otimes
I_{Z_1})
+ \chi (\scrO_X(\Delta -F)\otimes I_{Z_2})$, we may rewrite this last condition
as
$$\chi(\scrO_X(\Delta -F)\otimes I_{Z_2}) - \chi (\scrO_X(F)\otimes I_{Z_1})
\geq 0.$$
Now from the exact sequence
$$0 \to \scrO_X(F)\otimes I_{Z_1} \to \scrO_X(F) \to \scrO_{Z_1} \to 0,$$
we see that $\chi(\scrO_X(F)\otimes I_{Z_1}) = \chi(\scrO_X(F)) - \ell (Z_1)$,
and similarly $\chi(\scrO_X(\Delta -F)\otimes I_{Z_2}) = \chi(\scrO_X(\Delta
 -F)) - \ell(Z_2)$. By Riemann-Roch,
$$\align
\chi(\scrO_X(\Delta -F)) - \chi(\scrO_X(F)) &= \frac12((\Delta -F)^2 -
(\Delta -F)\cdot K_X - F^2 + F\cdot K_X)\\
&= \frac12(\Delta ^2 - 2\Delta \cdot F + \zeta \cdot K_X) \\
&= \frac12\zeta\cdot(K_X-\Delta)=t.
\endalign$$
Thus we have the following conditions on $Z_1$ and $Z_2$:
$$\align
\ell (Z_2) - \ell (Z_1) &\leq t;\\
\ell (Z_2) + \ell (Z_1) &= \ell _\zeta,
\endalign$$
and so $2\ell (Z_2) \leq \ell _\zeta +t$.
Setting $k= \dsize\fracwithdelims[]{\ell _\zeta +t}2$, we have $\ell (Z_2) \leq
k$. Applying a similar analysis to a subsheaf of the form $\scrO_X(\Delta -F)
\otimes I_{Z_1}$ shows that, if there is such a subsheaf, with a torsion free
quotient $\scrO_X(F)\otimes I_{Z_2}$, then
$$\ell (Z_2)\leq \frac{\ell _\zeta - t}2 = \ell _\zeta - \frac{\ell _\zeta
+t}2.$$ In particular, if $\dsize\frac{\ell _\zeta +t}2$ is not an integer,
then
this condition becomes $\ell (Z_2) \leq \ell _\zeta - k-1$. Thus, provided
$\dsize\frac{\ell _\zeta +t}2$ is not an integer for every $\zeta$ defining the
wall $W$ (i.e\. $t$ is not $\zeta$-critical for every $\zeta$), $V$ is $(L_0,
\zeta, k)$-semistable for
$k =\dsize\fracwithdelims[]{\ell _\zeta +t}2$ and indeed $V$ is $(L_0,
\boldsymbol\zeta, \bold k)$-semistable, where $\bold k$ is defined in the
obvious way. Conversely, assuming  that $t$ is not $\zeta$-critical for every
$\zeta$,
$V$ is Gieseker $L_0$-semistable, indeed Gieseker $L_0$-stable, if it is $(L_0,
\boldsymbol\zeta, \bold k)$-semistable for $\bold k$ as above.

We would like to produce a similar condition where $t$ is allowed to be any
rational number which is not $\zeta$-critical. One way to think of
this problem is to consider the analogous problem where we replace $\Delta$ by
$\Delta + 2 \Xi$ and make the corresponding change in $c$, so that $\Delta$ and
$p$ remain the same. This corresponds to twisting $V$ by $\scrO_X(\Xi)$, and
$t$
is replaced by
$t - \zeta \cdot \Xi$. In particular, we see that the notion of Gieseker
stability is rather sensitive to twisting by a line bundle. Moreover if $W$ is
defined by exactly one
$\zeta$ such that there exists a divisor $\Xi$ with $\zeta \cdot \Xi=1$, for
example if $\zeta$ is primitive and $p_g(X)=0$, it is easy to see that we can
construct the appropriate moduli spaces as Gieseker moduli spaces
corresponding to twists of $V$ by various multiples of
$\Xi$. In general however we will need to consider a problem which is roughly
analogous to allowing twists of $V$ by a $\Bbb Q$-divisor $\Xi$. This is the
goal of the following definition of mixed stability:

\definition{Definition 4.1} Let $X$ be an algebraic surface and let $L_0$ be an
ample line bundle on $X$. Fix line bundles $H_1$ and $H_2$ on $X$ and positive
integers $a_1$ and $a_2$. For every torsion free sheaf $V$ on $X$ of rank $r$,
define
$$p_{V; H_1, H_2, a_1, a_2}(n)= \frac{a_1}{r}\chi(V\otimes H_1 \otimes L_0^n) +
\frac{a_2}{r}\chi(V\otimes H_2 \otimes L_0^n).$$
A torsion free sheaf $V$ is {\sl$(H_1, H_2, a_1, a_2)$ $L_0$-stable\/} if,
for all subsheaves $W$ of $V$ with $0< \rank W < \rank V$ and for all $n \gg
0$,
$$p_{V; H_1, H_2, a_1, a_2}(n) > p_{W; H_1, H_2, a_1, a_2}(n).$$
$(H_1, H_2, a_1, a_2)$ $L_0$-semistable and unstable are defined similarly.
\enddefinition

The usual arguments show the following:

\lemma{4.2} If $V$ is $(H_1, H_2, a_1, a_2)$ $L_0$-stable, then it is simple.
\qed
\endstatement

In the case of rank two on a surface $X$ (which is the only case
which shall concern us), $V$ is $(H_1, H_2, a_1, a_2)$ $L_0$-stable if and only
if, for all rank one subsheaves $W$, and for all $n\gg 0$, we have
$$a_1(\chi(V\otimes H_1 \otimes L_0^n) - 2\chi(W\otimes H_1 \otimes L_0^n))+
a_2(\chi(V\otimes H_2 \otimes L_0^n) - 2\chi(W\otimes H_2 \otimes L_0^n)) >0.$$
In particular, if $V$ is $(H_1, H_2, a_1, a_2)$ $L_0$-stable then either
$V\otimes
H_1$ or $V\otimes H_2$ is stable, and a similar statement holds for
semistability. A short calculation shows that the coefficient of
$n$ in the above expression (which is a degree two polynomial in $n$) is
$(a_1+a_2)(L_0\cdot (c_1(V) - 2F))$ and that the constant term is
$$(a_1+a_2)(\chi(V) - 2\chi(W)) + a_1H_1\cdot (c_1(V) - 2F) + a_2H_2\cdot
(c_1(V)
 - 2F).$$ Thus $V$ is $(H_1, H_2, a_1, a_2)$ $L_0$-stable (resp\. semistable)
if
and only if it is either Mumford $L_0$-stable or Mumford strictly semistable
and
the above constant term is positive (resp\. nonnegative). It is easy to see,
comparing this with the discussion at the beginning of this section, that
formally
this is the same as requiring that $V\otimes \Xi$ is (Gieseker)
$L_0$-stable or semistable, where $\Xi$ is the $\Bbb Q$-divisor
$$\frac{a_1}{a_1+a_2}H_1 + \frac{a_2}{a_1+a_2}H_2.$$ Thus for example taking
$H_2=0$ and replacing $H_1$ by a positive integer multiple we see that we can
take for $\Xi$ an arbitrary $\Bbb Q$-divisor.

Let us explicitly relate mixed stability to our previous notion of $(L_0,
\boldsymbol \zeta, \bold k)$-semistability:

\lemma{4.3} Given $\Delta$ and $c$ and the corresponding $w$ and $p$, let $L_0$
be
an ample divisor lying on a unique wall of type $(w,p)$ and let $V$ be a rank
two torsion free sheaf with $c_1(V) = \Delta$ and $c_2(V)=c$. Let $\Xi$ be the
$\Bbb Q$-divisor $\dsize \frac{a_1}{a_1+a_2}H_1 + \frac{a_2}{a_1+a_2}H_2$ and
suppose that the rational number $t_i = \frac12\zeta_i\cdot(K_X-\Delta) -
\zeta_i
\cdot \Xi$ is not $\zeta_i$-critical for every $\zeta_i$ of type $(w,p)$
defining
$W$. Then, with $t=t_1$,
$V$ is $(L_0, \boldsymbol \zeta, \bold k(t))$-semistable
if and only if it is $(H_1, H_2, a_1, a_2)$
$L_0$-semistable if and only if it is $(H_1, H_2, a_1, a_2)$
$L_0$-stable.
\endstatement
\proof Using the additivity of the polynomials $p_{V; H_1, H_2, a_1, a_2}$ over
exact sequences, it is easy to check that $V$ is $(H_1, H_2, a_1, a_2)$
$L_0$-semistable if and only if it is Mumford $L_0$-semistable, and for every
Mumford destabilizing subsheaf of the form $\scrO_X(F)\otimes I_{Z_1}$, either
$V$ is universally semistable or we have
$$\chi (V) - 2\chi (\scrO_X(F)\otimes I_{Z_1}) - \zeta _i\cdot \Xi >0,$$
where $\zeta_i = 2F-\Delta$. Using our calculations above, this works out to
$$\ell (Z_2) - \ell (Z_1) \leq \frac12\zeta_i\cdot(K_X-\Delta) - \zeta_i
\cdot \Xi =t_i.$$
Equivalently since $\ell (Z_1) + \ell (Z_2) = \ell _{\zeta _i}$, this becomes
$\ell (Z_2) \leq \dsize \fracwithdelims[]{\ell _{\zeta _i}+t_i}2$. Thus $V$ is
$(H_1, H_2, a_1, a_2)$
$L_0$-semistable if and only if it is  $(L_0,
\boldsymbol \zeta,
\bold k(t))$-semistable. Moreover, since $t$ is not $\zeta_i$-critical, the
inequalities are automatically strict, so that $V$ is also $(H_1, H_2, a_1,
a_2)$
$L_0$-stable.
\endproof

Now choosing a $\Xi_0$ such that $\zeta _1\cdot \Xi _0 \neq 0$, every
rational number $t$ is of the form $\frac12\zeta_1\cdot(K_X-\Delta) - \zeta_1
\cdot r\Xi_0$ for some rational number $r$. Thus Theorem 3.5 will follow
from Lemma 4.3 and from the more general result below:

\theorem{4.4} Let $X$ be an algebraic surface $X$ and let  $L_0$ be an ample
line
bundle on $X$. Given a divisor $\Delta$ and an integer $c$,  line bundles $H_1$
and $H_2$ on
$X$ and positive integers $a_1$ and $a_2$, suppose that every rank
two torsion free sheaf
$V$ with $c_1(V) = \Delta$, $c_2(V) = c$ which is
$(H_1, H_2, a_1, a_2)$ $L_0$-semistable is actually $(H_1, H_2, a_1, a_2)$
$L_0$-stable. Then there exists a projective coarse moduli space
$\frak M _{L_0}(\Delta, c;H_1, H_2, a_1, a_2)$ of isomorphism classes of rank
two torsion free sheaves
$V$ with $c_1(V) = \Delta$, $c_2(V) = c$, which are
$(H_1, H_2, a_1, a_2)$ $L_0$-semistable.
\endstatement
\proof  The argument will follow the arguments in \cite{13} as closely as
possible,
and we shall assume a familiarity with that paper.

Suppose that $V$ is $(H_1, H_2, a_1, a_2)$ $L_0$-semistable. Then either
$V\otimes H_1$ or $V\otimes H_2$ is $L_0$-semistable, and thus by \cite{13},
Lemma 1.3 the set of all such $V$ is bounded. We may thus choose an $n$ such
that, for all
$V$ which are $(H_1, H_2, a_1, a_2)$ $L_0$-semistable, $V\otimes H_i \otimes
L_0^n$ is generated by its global sections and has no higher cohomology, for
$i=1,2$. Fix such an $n$ for the moment, and let $d_i = h^0(V\otimes H_i
\otimes
L_0^n)$. Then
$d_i$ is independent of $V$ and $V$ is a quotient of $(H_i^{-1}\otimes
L_0^{-n})^{\oplus d_i}$. Let $Q_i$ be the open subset of the corresponding Quot
scheme associated to $(H_i^{-1}\otimes L_0^{-n})^{\oplus d_i}$ consisting of
quotients which are rank two torsion free sheaves $V_i$ with $c_1(V_i) =
\Delta$
and
$c_2(V_i) = c$, and such that $V_i\otimes H_i \otimes L_0^n$ is generated by
its
global sections and has no higher cohomology. We will write a point of $Q_i$ as
$V_i$, suppressing the surjection
$(H_i^{-1}\otimes L_0^{-n})^{\oplus d_i} \to V_i$. Inside
$Q_1\times Q_2$, we have the closed subscheme
$I_0$ consisting of quotients $V_1$ and $V_2$ such that $\dim \Hom (V_1, V_2)
\geq 1$. There is also the open subvariety $I_0'$ of $I_0$ consisting of $(V_1,
V_2)$ with $\dim \Hom (V_1, V_2) =1$. Using the universal sheaves $\Cal U_i$
over $X\times Q_i$, we can construct a $\Cee ^*$ bundle $I$ over $I_0'$ whose
points are $(V_1, V_2, \varphi)$, where $\varphi\: V_1 \to V_2$ is a nonzero
homomorphism, unique up to scalars.

For $i=1,2$, let
$E_i$ be a fixed vector space of dimension equal to $d_i = h^0(V\otimes H_i
\otimes L_0^n)$. Fix once and for all an isomorphism $(H_i^{-1}\otimes
L_0^{-n})^{\oplus d_i} \cong (H_i^{-1}\otimes L_0^{-n})\otimes E_i$.
A surjection $(H_i^{-1}\otimes L_0^{-n})^{\oplus d_i}\to V_i$ then gives a map
$E_i \to H^0(V_i\otimes H_i \otimes L_0^n)$ and via such a surjection  a basis
$v_1,\dots, v_{d_1}$ of $E_1$ gives $d_1$ sections of $V_1
\otimes H_1 \otimes L_0^n$ and similarly for a basis $w_1, \dots, w_{d_2}$ of
$E_2$. Moreover
$GL(d_i)$ acts on $(H_i^{-1}\otimes L_0^{-n})^{\oplus d_i}$ and on $Q_i$. By
the universal property of the Quot scheme, this action extends to a
$GL(d_i)$-linearization of the universal sheaf $\Cal U_i$ over $X\times Q_i$.
Thus there is a right action of $GL(d_1) \times GL(d_2)$ on $I$, and it is easy
to see that the elements $(\lambda\Id, \lambda\Id)$ act trivially. Let
$F_i$ be the fixed vector space
$H^0(\Delta \otimes H_i^2 \otimes L_0^{2n})$, and $F$ the fixed vector space
$H^0(\Delta \otimes H_1 \otimes H_2 \otimes L_0^{2n})$. Let
$$U= \Hom(\bigwedge ^2E_1, F_1) \oplus \Hom(\bigwedge ^2E_2, F_2)\oplus \Hom
(E_1\otimes E_2, F).$$
(The factor $\Hom (E_1\otimes E_2, F)$ is there to make sure that the
destabilizing subsheaves for $V\otimes H_1$ and $V\otimes H_2$ are
in fact the same.) Note that $GL(d_1) \times GL(d_2)$ operates on
the right on $U$ and $\Pee U$. For example, the pair $(\lambda\Id, \mu\Id)$
acts on the triple $(T_1, T_2, T) \in U$ via
$(T_1, T_2, T) \mapsto (\lambda ^2T_1, \mu ^2T_2, \lambda\mu T)$.
Thus $(A_1, A_2)$ acts trivially on $\Pee U$ if and only if
$(A_1, A_2) = (\lambda\Id, \lambda \Id)$. Given a quintuple
$\underline{V} = (V_1, V_2,\psi _1, \psi_2, \varphi)$, where
$V_i \in Q_i$, $\psi _i\: E_i \to H^0(V_i \otimes H_i \otimes L_0^n)$
is an isomorphism, and $\varphi\: V_1 \to V_2$ is a nonzero map,
we will define a point
$(T_1(\underline{V}), T_2(\underline{V}), T(\underline{V})) \in \Pee U$.
To do so, fix an isomorphism $\alpha _2 \: \det V_2 \to \scrO_X(\Delta)$,
and set $\alpha _1 = \alpha _2 \circ \det\varphi$. (Thus $\alpha_1 = 0$
if $\varphi$ is not an isomorphism.) Given $v,v' \in E_1$ and
$w, w' \in E_2$, identify $v,v'$ with their images in
$H^0(V_i \otimes H_1 \otimes L_0^n)$ and similarly for $w,w'$, and let
$$\align
T_1(\underline{V})(v\wedge v') &= \alpha _1(v\wedge v') =\alpha
_2\circ
\det\varphi (v\wedge v') \in H^0(\Delta \otimes H_1^2\otimes L_0^{2n});\\
T_2(\underline{V})(w\wedge w') &= \alpha _2(w\wedge w') \in H^0(\Delta \otimes
H_2^2\otimes L_0^{2n});\\
T(\underline{V})(v\otimes w) &= \alpha _2(\varphi(v)\wedge w) \in H^0(\Delta
\otimes H_1 \otimes H_2 \otimes L_0^{2n}).
\endalign$$
Changing $\alpha _2$ by a nonzero scalar $\lambda$ multiplies
$(T_1(\underline{V}), T_2(\underline{V}), T(\underline{V}))$ by $\lambda$, so
that
the induced element of $\Pee U$ is well defined. Similarly, if we replace
$\varphi$ by
$\lambda\varphi$, then $(T_1(\underline{V}), T_2(\underline{V}),
T(\underline{V}))$ is replaced by $(\lambda ^2T_1(\underline{V}),
T_2(\underline{V}),
\lambda T(\underline{V}))$. It is easy to check that the map $\underline{V}
\mapsto T(\underline{V})$ induces a morphism from
$I$ to
$\Pee U$ which is $GL(d_1) \times GL(d_2)$-equivariant. Further note that we
can
define
$(T_1(\underline{V}), T_2(\underline{V}), T(\underline{V}))$ more generally if
we are given the data $\underline{V}$ of two rank two torsion free sheaves
$V_1$
and $V_2$ with $\det V_i = \Delta$, a morphism $\varphi\: V_1 \to V_2$, and
linear maps
$\psi _i\: E_i \to H^0(V_i \otimes H_i \otimes L_0^n)$, not necessarily
isomorphisms, although it is possible for $(T_1(\underline{V}),
T_2(\underline{V}), T(\underline{V}))$ to be zero in this case.

We have not yet introduced the extra parameters $a_1$ and $a_2$. To do so,
define $G(a_1, a_2) \subset GL(d_1) \times GL(d_2)$ as follows:
$$G(a_1, a_2) = \{\, (A_1,A_2)\mid \det A_1^{a_1}\det A_2^{a_2} = \Id\,\}.$$
Thus unlike Thaddeus we don't change the polarization or the linearization but
the actual group which we use to determine stability; still our construction
could
probably be interpreted in his general framework. Fixing $a_1$ and $a_2$ for
the
rest of the discussion, we shall denote $G(a_1, a_2)$ by $G$. Since $a_1$ and
$a_2$ are positive, the matrix $(\lambda\Id, \lambda\Id)$ lies in $G$ if and
only
if
$\lambda$ is an $m^{\text{th}}$ root of unity, where $m = a_1d_1+a_2d_2$. Thus
a
quotient of
$G$ by a finite group acts faithfully on $\Pee U$. Moreover, the problem
of finding a good quotient of
$\Pee U$ (for an appropriate open subset of $\Pee U$) for $G$ is the same  as
that of finding a good quotient of $\Pee U$ for $GL(d_1) \times GL(d_2)$, since
$$G\cdot \Cee ^*(\Id, \Id) = GL(d_1) \times GL(d_2).$$
This last statement follows since $G$ clearly contains $SL(d_1) \times SL(d_2)$
and since $\Cee ^*\times \Cee ^*$ is generated by its diagonal subgroup and by
the
subgroup
$$\{\,(\lambda, \mu):
\lambda ^{a_1d_1}\mu ^{a_2d_2} = 1\,\}.$$

We may thus apply the general machinery of GIT to the group $G$ acting on $\Pee
U$. A one parameter subgroup of $G$ is given by a basis $\{v_i\}$ of $E_1$, a
basis $\{w_k\}$ of $E_2$ and weights $n_i$, $m_k\in \Zee$, such that
$v_i^\lambda
= \lambda ^{n_i}v_i$, $w_k ^\lambda = \lambda ^{m_k}w_k$, and
$$a_1\sum _in_i + a_2 \sum _km_k = 0.$$
We shall always arrange our choice of basis so that $n_1 \leq n_2 \leq \dots
\leq n_{d_1}$ and $m_1 \leq m_2 \leq \dots \leq m_{d_2}$. Given
$(T_1, T_2, T)
\in U$ and a one parameter subgroup of $G$ as above, we see that $\lim
_{\lambda
\to 0}(T_1, T_2, T)^\lambda =0$  if and only if $T_1(v_i\wedge v_j) = 0$ for
every
pair of indices
$i,j$ such that $n_i+n_j \leq 0$, $T_2(w_k\wedge w_\ell) = 0$ for every pair of
indices $k, \ell$ such that $m_k+m_\ell \leq 0$, and  $T(v_i\otimes w_j) = 0$
for
every pair $i,k$ such that $n_i+m_k\leq 0$. Likewise the condition that $\lim
_{\lambda
\to 0}(T_1, T_2, T)^\lambda$ exists is similar, replacing the $\leq$ by strict
inequality. Finally note that if $n_i+n_j \leq 0$, then $n_1+n_j \leq 0$, if
$m_k+m_\ell \leq 0$ then $m_1+m_\ell \leq 0$, and if $n_i+m_k\leq 0$ then
$n_1+m_k\leq 0$ and $n_i+m_1\leq 0$.

We then have the following:

\lemma{4.5}
\roster
\item"{(i)}" Suppose that we are given the data $\underline{V}$ of two rank two
torsion free sheaves $V_1$ and $V_2$ with $\det V_i = \Delta$, a morphism
$\varphi\: V_1 \to V_2$, and a linear map
$E_i \to H^0(V_i \otimes H_i \otimes L_0^n)$, not necessarily an isomorphism.
If $E_i \to H^0(V_i \otimes H_i \otimes L_0^n)$ is not injective for some
$i$ or if $\varphi$ is not an isomorphism, then $(T_1(\underline{V}),
T_2(\underline{V}),  T(\underline{V}))$ is either zero or $G$-unstable.
\item"{(ii)}" For $n$ sufficiently large depending only on $\Delta$ and $c$ and
for
$V$ a rank two torsion free sheaf with $\det V = \Delta$ and $c_2(V) = c$,
$V$ is $(H_1, H_2, a_1, a_2)$
$L_0$-unstable if and only if $(T_1(\underline{V}), T_2(\underline{V}),
T(\underline{V}))$ is $G$-unstable for all choices of data $\underline{V}$ such
that $E_i \to H^0(V_i \otimes H_i \otimes L_0^n)$ is injective and $\varphi\:
V_1 \to V_2\cong V$ is an isomorphism, and
$V$ is $(H_1, H_2, a_1, a_2)$
$L_0$-strictly semistable if and only if $(T_1(\underline{V}),
T_2(\underline{V}), T(\underline{V}))$ is $G$-strictly semistable for all such
$\underline{V}$. Thus
$V$ is $(H_1, H_2, a_1, a_2)$
$L_0$-stable if and only if $(T_1(\underline{V}), T_2(\underline{V}),
T(\underline{V}))$ is $G$-stable for all such $\underline{V}$.
\endroster
\endstatement
\proof First let us prove (i). We may assume that $(T_1(\underline{V}),
T_2(\underline{V}),  T(\underline{V}))\neq 0$. Suppose for example that $v_1
\in E_1
\mapsto 0
\in H^0(V_1\otimes H_1\otimes L_0^n)$. Complete $v_1$ to a basis of $E_1$ and
choose a basis $\{w_k\}$ for $E_2$. Then
$T_1(\underline{V})(v_1\wedge v_i) = 0$ for all $i$ and
$T(\underline{V})(v_1\otimes w_k) =0$ for all $k$. Define a one parameter
subgroup of $G$ as follows: let $v_1^\lambda = \lambda ^{-N}v_1$, $v_i^\lambda
=
\lambda ^av_i$ for $i>1$, and $w_k^\lambda = \lambda ^bw_k$ for all $k$.
Clearly
$\lim _{\lambda \to 0}(T_1(\underline{V}), T_2(\underline{V}),
T(\underline{V}))^\lambda =0$ provided that $a$ and $b$ are positive, so that
$(T_1(\underline{V}), T_2(\underline{V}),  T(\underline{V}))$ is $G$-unstable
provided that the one parameter subgroup so constructed lies in $G$, or on
other
words provided that
$$a_1(-N + a(d_1-1))+ a_2bd_2 =0.$$
It thus suffices to take $a$ an arbitrary positive integer, $b=a_1$, and $N=
a(d_1-1) + a_2d_2$.

The argument in case $\varphi$ has a kernel is similar: in this case let $v_1
\in \Ker \varphi$. Then $T_1(\underline{V})=0$ and
$T(\underline{V})(v_1\otimes
w_k) =0$ for all $k$, so that the previous argument handles this case also.

Next we show (ii). Let $p_{V\otimes H_i}$ be the usual normalized Hilbert
polynomial of $V\otimes H_i$, and similarly for $p_{W\otimes H_i}$, where $W$
is a rank one subsheaf of $V$. Thus $p_{V\otimes H_i}$ and $p_{W\otimes
H_i}$ have the same leading term. Given a polynomial
$p$, let $\Delta p$ denote the difference polynomial. In our case, all of the
polynomials $p$ that occur are quadratic polynomials with the same fixed degree
two term. Thus if $p_1$ and
$p_2$ are two such polynomials, then
$p_1(n) > p_2(n)$ for all $n\gg 0$ if and only if the linear term of $p_1$ is
greater than or equal to the linear term of $p_2$, and if the linear terms are
equal then the constant term of $p_1$ is greater than the constant term of
$p_2$.
In this last case, where the linear terms are also equal, we see that $p_1(n) >
p_2(n)$ for all
$n\gg 0$ if and only if $p_1(n) > p_2(n)$ for some $n$. Finally the linear
term of $p_1$ is greater than or equal to the linear term of $p_2$ if and only
if
$\Delta p_1(n) \geq \Delta p_2(n)$ for all $n$, which we shall write as $\Delta
p_1\geq \Delta p_2$. Thus if $\Delta p_1 \geq \Delta p_2$ and $p_1(n) > p_2(n)$
for some $n$, then $p_1(n) > p_2(n)$ for all $n\gg 0$. If $\Delta p_1 =
\Delta p_2$, then $p_1(n) > p_2(n)$
for some $n$ if and only if $p_1(n) > p_2(n)$ for all $n$.

We shall show that, for sufficiently large $n$, if $V$ is $(H_1, H_2, a_1,
a_2)$
$L_0$-semistable and $\underline{V}$ corresponds to data where $E_i \to
H^0(V\otimes H_i\otimes L_0^n)$ is injective and $\varphi$ is an isomorphism,
then
$(T_1(\underline{V}), T_2(\underline{V}), T(\underline{V}))$  is
$G$-semistable.
Note that $V$ is Mumford semistable. First we may choose
$n$ so that
$V\otimes H_i$ is generated by its global sections and has no higher
cohomology,
and so $\chi (V\otimes H_i \otimes L_0^n) = h^0(V\otimes H_i \otimes L_0^n) =
d_i$. Hence, since $E_i \to H^0(V\otimes H_i\otimes L_0^n)$ is injective, it
is an isomorphism. Let
$W$ be a rank one subsheaf of $V$. Since $V\otimes H_i$ is Mumford semistable,
$\Delta p_{W\otimes H_i} \leq \Delta p_{V\otimes H_i}$. Now the proof of (3) of
Lemma 1.2 in \cite{13}  shows that
there exists an $N$ so that, for all $n\geq N$, with $d_i$ as above, if
$W$ is a rank one subsheaf of
$V$ and such that  $h^0(W\otimes H_i \otimes L_0^n )
\geq d_i/2$ for at least one
$i$ ($i=1,2$), then in fact $\Delta p_{V\otimes H_i} = \Delta
p_{W\otimes H_i}$ for all such $W$, and thus $\mu _{L_0}(V) = \mu _{L_0}(W)$.
It is then easy to see that there is a twist $W\otimes H_i\otimes L_0^{-k}$,
depending only on $L_0$ and $\Delta$, such that $h^0((V/W)\otimes H_i\otimes
L_0^{-k}) = 0$. The proof of Proposition 3.1 in \cite{13} shows that in this
case
$h^1(W\otimes H_i \otimes  L_0^{-k})$ is bounded by $Q$, where $Q$ is some
universal bound for the numbers
$h^1(V\otimes H_i\otimes L_0^{-k})$ as $V\otimes H_i$ ranges over the
appropriate set of
$L_0$-semistable sheaves  Thus by (4) of Lemma 1.2 in \cite{13}, the $W$
satisfying the condition that $h^0(W\otimes H_i \otimes L_0^n )
\geq d_i/2$ for at least one $i$ form a bounded family, and we may
choose $n$ so large, depending only on $L_0$, $\Delta$, $c$, such that
$h^j(W\otimes H_i \otimes L_0^n)=0$ for $j\geq 1$ and $i=1,2$.

Now suppose that
$(T_1(\underline{V}), T_2(\underline{V}), T(\underline{V}))$ is $G$-unstable.
Then there exists a  one parameter subgroup of $G$ as above such that
$\lim_{\lambda \to 0} (T_1(\underline{V}), T_2(\underline{V}),
T(\underline{V}))^\lambda =0$. Let
$$\align
s_1&= \#\{\, j: T_1(\underline{V})(v_1\wedge v_j) =0\,\} \geq
\max \{\, j: n_1 + n_j \leq 0\,\};\\
s_2 &= \#\{\, j: T_2(\underline{V})(w_1\wedge
w_\ell) =0\,\} \geq
\max \{\, \ell: m_1 + m_\ell \leq 0\,\}.
\endalign$$
Since $a_1\sum _in_i + a_2 \sum _km_k =0$, at least one of $n_1, m_1$ is
negative. By symmetry we may assume that $n_1$ is negative, and that $n_1 \leq
m_1$.  Since for
$j\leq s_1$, $v_1\wedge v_j$ is zero as a section of
$\det (V\otimes H_1\otimes L_0^n)$, the sections corresponding to
$v_j$, $1\leq j \leq s_1$, are all sections of a rank one subsheaf $W_1$ of
$V$. Likewise the sections $w_\ell$, $1\leq \ell \leq s_2$, if there are any
such, are all sections of a rank one subsheaf $W_2$ of $V$.
The condition that $T(\underline{V})(v_1\otimes w_1) =0$
insures that $W_1$ and $W_2$ are contained in a  saturated
rank one subsheaf $W$, if $s_2\neq 0$, otherwise we shall just take for $W$ the
saturated rank one subsheaf containing $W_1$. Moreover
$h^0(W\otimes H_1\otimes L_0^n)
\geq s_1$  and $h^0(W\otimes H_2\otimes L_0^n) \geq s_2$. Suppose that we show
that
$$a_1(d_1 - 2s_1) + a_2(d_2 - 2s_2) <0.$$
Thus in particular $s_i\geq d_i/2$ for at least one $i$. By our choice of $n$
and  the previous paragraph, if $s_i\geq d_i/2$ for at least one $i$,
then $h^0(W\otimes H_i\otimes L_0^n)  = \chi (W\otimes H_i\otimes L_0^n)$
and furthermore $\mu _{L_0}(V) = \mu _{L_0}(W)$.
Thus
$$h^0(W\otimes H_i \otimes L_0^n)
= \chi (W\otimes H_i \otimes L_0^n) \geq s_i$$
for $i=1,2$ and so
$p_{V; H_1, H_2, a_1, a_2}(n) < p_{W; H_1, H_2, a_1, a_2}(n)$.
On the other hand, $p_{V; H_1, H_2, a_1, a_2}$ and $p_{W; H_1, H_2, a_1, a_2}$
are two quadratic polynomials with the same linear and quadratic terms
(since $\mu _{L_0}(V) = \mu _{L_0}(W)$),
and $p_{V; H_1, H_2, a_1, a_2}(n) < p_{W; H_1, H_2, a_1, a_2}(n)$
for one value of $n$. Thus the constant term of $p_{W; H_1, H_2, a_1, a_2}$
must be larger than that of $p_{V; H_1, H_2, a_1, a_2}$.
This contradicts the $(H_1, H_2, a_1, a_2)$ $L_0$-semistability of $V$.

To see that $a_1(d_1 - 2s_1) + a_2(d_2 - 2s_2) <0$, let
$$t_1 = \#\{\, j: n_j +m_1 \leq 0\,\} \leq s_1.$$
Here $t_1 \leq s_1$ since $T(\underline{V})(v_j \otimes w_1) = 0$
implies that $v_j$ and $w_1$ are contained in a rank one subsheaf of $V$,
necessarily $W$, and thus that $v_1\wedge v_j = 0$. Let
$$t_2 = \#\{\, \ell: n_1 +m_\ell \leq 0\,\} \leq s_2.$$ We
have assumed that $n_1 \leq m_1$. Then consider the expression
$$a_1\sum _j (n_1+ n_j) + a_2\sum _\ell (n_1 + m_\ell).$$
On the one hand from the definition of the one parameter subgroup we have
$$a_1\sum _j (n_1+ n_j) + a_2\sum _\ell (n_1 + m_\ell)
= a_1d_1n_1 + a_2d_2n_1.$$
On the other hand, to estimate $\sum _j (n_1+ n_j)$, we can ignore the positive
terms where $n_1+n_j \geq 0$ and each of the $s_1$ negative terms are at least
$n_1 + n_1 \geq 2n_1$. Thus $\sum _j (n_1+ n_j) \geq 2s_1n_1$.
Since $n_1 <0$, this term is $\geq 2s_1n_1$.
Also this inequality is strict or $n_1 + n_i \leq 0$ for
every $i$, which would say that every section of $V\otimes H_1 \otimes L_0^n$
is
really a section of $W\otimes H_1 \otimes L_0^n$ contradicting the fact that
$V\otimes H_1 \otimes L_0^n$ is generated by global sections. So $\sum _j (n_1+
n_j) < 2s_1n_1$. Likewise we claim that $\sum _\ell (n_1 + m_\ell) \geq
2s_2n_1$. Here, to estimate $\sum _\ell (n_1 + m_\ell)$, we may ignore the
terms with $n_1 + m_\ell$ positive, leaving $t_2$ terms $n_1+m_\ell$ which are
$\leq 0$, and moreover each such term is at least $n_1+m_1 \geq 2n_1$. Thus
$\sum _\ell (n_1 + m_\ell) \geq 2t_2n_1$, and since $t_2 \leq s_2$ and $n_1
<0$, we have $2t_2n_1 \geq 2s_2n_1$.

Putting this together we have
$$\align
a_1d_1n_1 + a_2d_2n_1 &= a_1\sum _j (n_1+ n_j) + a_2\sum _\ell
(n_1 + m_\ell) \\ &> a_1(2s_1n_1)+ a_2(2s_2n_1),
\endalign$$
so that
$$a_1(d_1 - 2s_1)n_1 + a_2(d_2 - 2s_2)n_1 >0.$$ As $n_1 <0$, we must have
$a_1(d_1 - 2s_1) + a_2(d_2 - 2s_2) <0$, as desired.

We have thus shown that, if $(T_1(\underline{V}), T_2(\underline{V}),
T(\underline{V}))$  is $G$-unstable, then $V$ is $(H_1, H_2, a_1, a_2)$
$L_0$-unstable. A very similar argument
handles the $G$-strictly semi\-stable case.

Now we turn to the converse statement, that if $V$ is $(H_1, H_2, a_1, a_2)$
$L_0$-unstable then $(T_1(\underline{V}), T_2(\underline{V}),
T(\underline{V}))$ is $G$-unstable. Suppose instead that
$$(T_1(\underline{V}), T_2(\underline{V}), T(\underline{V}))$$
is $G$-semistable. Let $W$ be a rank one subsheaf of $V$ such that
$p_{W; H_1, H_2, a_1, a_2}(m) > p_{V; H_1, H_2, a_1, a_2}(m)$ for all $m\gg 0$.
We may assume that the quotient $W' = V/W$ is torsion free.
Thus $p_{W'; H_1, H_2, a_1, a_2}(m) < p_{V; H_1, H_2, a_1, a_2}(m)$
for all $m\gg 0$, and so
$\Delta p_{W'\otimes H_i} \leq \Delta p_{V\otimes H_i}$.
Now we have the map $E_i \to H^0(V\otimes H_i \otimes L_0^n)$.
Consider $E_i \cap H^0(W\otimes H_i \otimes L_0^n)\subseteq E_i$.
Let $\dim E_i \cap H^0(W\otimes H_i \otimes L_0^n) = s_i$.
Suppose first that $a_1(d_1 - 2s_1) + a_2(d_2 - 2s_2) <0$.
We claim that in this case
$(T_1(\underline{V}), T_2(\underline{V}), T(\underline{V}))$  is $G$-unstable,
a contradiction. To see this, choose a basis $v_1, \dots, v_{d_1}$ for
$E_1$ such that
$$v_i \in E_1\cap H^0(W\otimes H_1 \otimes L_0^n)$$
for $i\leq s_1$, and similarly choose a basis $w_1, \dots, w_{d_2}$ for
$E_2$ such that $w_k \in E_2\cap H^0(W\otimes H_2 \otimes
L_0^n)$ for $i\leq s_2$. Thus, if $i,j\leq s_1$ then
$T_1(\underline{V})(v_i\wedge v_j) =0$; if $k, \ell \leq s_2$ then $T_2(
\underline{V})(w_k\wedge w_\ell) = 0$; if $i\leq s_1$ and $k\leq s_2$ then
$T(\underline{V})(v_i\otimes w_k) =0$.

We will try to find a one parameter subgroup of $G$ of the form
$$v_i^\lambda = \cases \lambda ^{-m}v_i, &\text{for $i \leq s_1$;} \\
\lambda ^nv_i , &\text{for $i > s_1$,}
\endcases$$ and similarly
$$w_k^\lambda = \cases \lambda ^{-m}w_k, &\text{for $i \leq s_2$;} \\
\lambda ^nw_k , &\text{for $i > s_2$.}
\endcases$$
It is easy to check that $\lim _{\lambda \to 0}(T_1(\underline{V}),
T_2(\underline{V}), T(\underline{V}))^\lambda=0$ if and only if $n>m$.
What we must arrange is the condition
$$a_1(-ms_1 + n(d_1-s_1)) + a_2(-ms_2 +n(d_2-s_2)) = 0.$$
Now consider the linear function with rational coefficients
$$f(t) = a_1(-s_1 + t(d_1-s_1)) + a_2(-s_2 +t(d_2-s_2)).$$
Since the coefficient of $t$ is strictly positive $f(t)$ is increasing, and
$$\align f(1) &= a_1(-s_1 + (d_1-s_1)) + a_2(-s_2 +(d_2-s_2))\\
&= a_1(d_1 - 2s_1) + a_2(d_2 - 2s_2) <0.
\endalign$$
Thus there is a rational number $t= n/m>1$ such that $f(t) = 0$,
and this gives the desired choice of $n$ and $m$.
Thus if $a_1(d_1 - 2s_1) + a_2(d_2 - 2s_2) <0$, then
$(T_1(\underline{V}), T_2(\underline{V}), T(\underline{V}))$ is $G$-unstable,
contradicting our hypothesis.

The other possibility is that $a_1(d_1 - 2s_1) +
a_2(d_2 - 2s_2) \geq 0$. In this case  $d_i \geq 2 s_i$ for at least one $i$.
Recalling that we have the quotient $W'$ of $V$ by $W$, it then follows that
for
such an $i$ the image of $E_i$ in $H^0(W' \otimes H_i \otimes L_0^n)$ must have
dimension at least $d_i/2$. Arguing as in Proposition 3.2 of \cite{13},
it then follows from Lemma 1.2 of \cite{13} that
$\Delta p_{W'\otimes H_i} = \Delta p_{V\otimes H_i}$ and so that $V$ is
Mumford $L_0$-semistable and $\mu _{L_0}(V) = \mu _{L_0}(W)$.
Moreover, after enlarging $n$ if necessary (independently of $V$)
we may assume that $h^j(V\otimes H_i \otimes L_0^n) = 0$ for $j>0$.
In particular, $d_i = \dim H^0(V\otimes H_i \otimes L_0^n)$ for $i=1,2$,
and $E_i \to H^0(V \otimes H_i \otimes L_0^n)$ is an isomorphism;
so $s_i = h^0(W \otimes H_i \otimes L_0^n)$.
As $\mu _{L_0}(V) = \mu _{L_0}(W)$, the polynomials
$p_{W; H_1, H_2, a_1, a_2}$ and $p_{V; H_1, H_2, a_1, a_2}$ have
the same terms in degree one and two, and thus since
$p_{W; H_1, H_2, a_1, a_2}(m) > p_{V; H_1, H_2, a_1, a_2}(m)$ for some $m$
the same is true for all $m$, in particular for $m=n$. Moreover,
for a general choice of a smooth curve $C$ in the linear system corresponding
to $L_0$, there is a fixed bound on the line bundle $W\otimes H_i|C$.
A standard argument as in the proof of (2) of Lemma 1.2 of \cite{13}
shows that, for $n$ sufficiently large but independent of $V$,
we have $H^2(W\otimes H_i\otimes L_0^n) =0$. Thus
$s_i = h^0(W \otimes H_i \otimes L_0^n) \geq p_{W\otimes H_i}(n)$.
It follows that
$$\align
a_1(d_1 - 2s_1) +  a_2(d_2 - 2s_2)&\leq a_1(d_1 - 2p_{W\otimes H_1}(n)) +
a_2(d_2 - 2p_{W\otimes H_2}(n))\\
&= 2(p_{V; H_1, H_2, a_1, a_2}(n) - p_{W; H_1, H_2, a_1, a_2}(n)) <0.
\endalign$$
This contradicts the assumption that $a_1(d_1 - 2s_1) +
a_2(d_2 - 2s_2) \geq 0$. It then follows that $(T_1(\underline{V}),
T_2(\underline{V}), T(\underline{V}))$ is $G$-unstable.

The strictly semistable case is similar.
\endproof

We may now finish the proof of Theorem 4.4. Let $\Pee U_{\text{ss}}$ be the set
of $G$-semistable points of $\Pee U$.  Let $I_{\text{ss}}$ be the inverse image
of $\Pee U_{\text{ss}}$ under the morphism $I \to \Pee U$. Since every
semistable sheaf is stable, $I_{\text{ss}}$ is a $\Cee ^*$-bundle over
its image in $Q_1\times Q_2$. Moreover the representable functor
corresponding to $I_{\text{ss}}$ is easily seen to be formally smooth
over the moduli functor. Arguments very similar to those for Lemma 4.3 and 4.5
of
\cite{13} show that the morphism
$I_{\text{ss}} \to \Pee U_{\text{ss}}$ is one-to-one and proper, and thus in
particular finite. Thus we may construct a quotient $\frak M _{L_0}(\Delta,
c;H_1, H_2, a_1, a_2)$ of $I_{\text{ss}}$ by $G$. This quotient maps in a
one-to-one and proper way to the GIT quotient of $\Pee U_{\text{ss}}$ and is
therefore projective.  By the discussion at the beginning of the proof of
Theorem
4.4 the points of $\frak M _{L_0}(\Delta, c;H_1, H_2, a_1, a_2)$ may be
identified with isomorphism classes of $(H_1, H_2, a_1, a_2)$
$L_0$-semistable rank two sheaves. Standard arguments then show that $\frak M
_{L_0}(\Delta, c;H_1, H_2, a_1, a_2)$ has the usual properties of a coarse
moduli space.
\endproof

\section{5. The transition formula for Donaldson polynomial invariants.}

 From now on, we will assume that the surface $X$ is rational with
$-K_X$ effective, and will study the transition formula of
Donaldson polynomial invariants:
$$\delta ^X_{w,p}(\Cal C_+, \Cal C_-)
= D^X _{w,p}(\Cal C_+) - D^X _{w,p}(\Cal C_-)$$
where $\Cal C_-$ and $\Cal C_+$ are two adjacent chambers separated by
a single wall $W^{\zeta}$ of type $(w, p)$
or equivalently of type $(\Delta, c)$. For simplicity, we assume that
the wall $W^{\zeta}$ is only represented by $\pm \zeta$ since
the general case just involves additional notation. We use $\frak
M_0^{(k)}$ to stand for the moduli space $\frak M_0^{(\boldsymbol \zeta, k)}$.
When $\ell_\zeta = 0$, we also assume that
$$h(\zeta) = h^1(X; \scrO _X(2F-\Delta)  \neq 0$$
(see Corollary 2.7). The special case when $\ell_\zeta = h(\zeta) = 0$
will be treated in Theorem 6.1. By Theorem 3.9 and Lemma 3.2 (ii),
we have the following diagram:
$$\matrix
& &{\widetilde{\frak M}}_0^{(\ell_\zeta)}&&&&\ldots&&&&{\widetilde{\frak
M}}_0^{(0)}&&\\ &\swarrow&&\searrow&&\swarrow&&\searrow&&\swarrow&&\searrow\\
\frak M_0^{(\ell_\zeta)} & & && \frak M_0^{(\ell_\zeta - 1)} &
&&& \frak M_0^{(0)} & &&& \frak M_0^{(-1)}\\
\Vert \quad&&&&&&&&&&&& \Vert \quad\\
\frak M_- &&&&&&&&&&&&\frak M_+\\
\endmatrix
$$
where the morphism ${\widetilde{\frak M}}_0^{(k)} \to \frak M_0^{(k)}$
is the blowup of $\frak M_0^{(k)}$ at $E_\zeta^{\ell_\zeta -k, k}$,
and the morphism ${\widetilde{\frak M}}_0^{(k)} \to \frak M_0^{(k - 1)}$
is the blowup of $\frak M_0^{(k - 1)}$ at $E_{-\zeta}^{k, \ell_\zeta -k}$.

Next, we collect and establish some notations. Recall that in section 2
we have constructed the bundle $\Cal E_\zeta^{\ell_\zeta -k, k}$ over
$ H_{\ell_\zeta - k} \times H_{k}$, where $H_k = \Hilb^k X$.

\medskip
\noindent
{\bf Notation 5.1}. Let $\zeta$ define a wall of type $(w, p)$.
\roster

\item"{(i)}" $\lambda_k$ is the tautological line bundle over
$E_\zeta^{\ell_\zeta - k, k} = \Pee((\Cal E_\zeta^{\ell_\zeta -k,
k})\spcheck)$;  for simplicity, we also use $\lambda_k$ to denote its first
Chern class;

\item"{(ii)}" $\rho_k: X \times E_\zeta^{\ell_\zeta - k, k} \to
X \times H_{\ell_\zeta - k} \times H_{k}$ is the natural projection;

\item"{(iii)}" $p_k: \widetilde{\frak M}_0^{(k)} \to \frak M_0^{(k)}$ is the
blowup of $\frak M_0^{(k)}$ at $E_\zeta^{\ell_\zeta -k, k}$;

\item"{(iv)}" $q_{k - 1}: \widetilde{\frak M}_0^{(k)} \to \frak M_0^{(k-1)}$ is
the contraction of $\widetilde{\frak M}_0^{(k)}$ to $\frak M_0^{(k-1)}$;

\item"{(v)}" $\Cal N_k$ is the normal bundle of $E_\zeta^{\ell_\zeta -k, k}$
in $\frak M_0^{(k)}$; by Proposition 3.7, we have
$$\Cal N_k = \rho_k^*\Cal E_{-\zeta}^{k, \ell_\zeta -k}
\otimes \lambda_k^{-1};$$

\item"{(vi)}" $D_k = \Pee(\Cal N_k\spcheck)$ is the exceptional divisor in
$\widetilde{\frak M}_0^{(k)}$;

\item"{(vii)}" $\xi_k = \Cal O_{\widetilde{\frak M}_0^{(k)}}(-D_k)|D_k$ is
the tautological line bundle on $D_k$; again, for simplicity, we also use
$\xi_k$ to denote its first Chern class;

\item"{(viii)}" $\mu^{(k)}(\alpha) = -{1 \over 4} p_1(\Cal U^{(k)})/\alpha$
where
$\alpha \in H_2(X; \Zee)$ and $\Cal U^{(k)}$ is a universal sheaf
over $X \times \frak M_0^{(k)}$.
Let $\mu^{(\ell_\zeta)}(\alpha) = \mu_-(\alpha)$ and
that $\mu^{(-1)}(\alpha) = \mu_+(\alpha)$.

\item"{(ix)}" $\nu^{(k)} = -{1 \over 4} p_1(\Cal U^{(k)})/x$ where
$x \in H_0(X; \Zee)$ is the natural generator. Let
$\nu^{(\ell_\zeta)}= \nu_-$ and that $\nu^{(-1)} = \nu_+$.

\endroster
Note that, in (viii) and (ix) above, the sheaf $\Cal U^{(k)}$ is only defined
locally in the classical topology. However, since it is defined on the level of
the Quot scheme  a straightforward argument shows that $p_1(\Cal U^{(k)})$ is a
well-defined element in the rational cohomology of $X \times \frak M_0^{(k)}$,
at least in the complement of the universally semistable sheaves. In case there
are universally semistable sheaves, then the work of Li \cite{21} extends the
$\mu$-map to $\frak M_0^{(k)}$, at least for the two-dimensional algebraic
classes. We can then extend the $\mu$-map to the 4-dimensional class via a
blowup
formula due to O'Grady (unpublished). Moreover, there
is a universal sheaf
$\Cal V_k$ over
$X \times E_\zeta^{\ell_\zeta - k, k}$. In what follows, we shall work as if
there were a universal sheaf $\Cal U^{(k)}$, and leave it to the reader to
check
that our final Chern class calculations can be verified directly even when no
universal sheaf exists.

In the following lemma, we study the restrictions of
$p_k^*\mu^{(k)}(\alpha)$ and $p_k^*\nu^{(k)}$ to $D_k$.

\lemma{5.2} Let $\alpha \in H_2(X; \Zee)$ and $a = (\zeta \cdot \alpha)/2$.
Let $\tau_1$ and $\tau_2$ be the projections of $E_\zeta^{\ell_\zeta - k, k}$
to $H_{\ell_\zeta - k}$ and $H_k$ respectively. Then,
$$\align
&(\operatorname{Id} \times p_{k})^*c_1(\Cal U^{(k)})|(X \times D_k)
= \pi _1^*\Delta + (p_k|D_k)^*\lambda_k\\
&p_k^*\mu^{(k)}(\alpha)|D_k =
(p_k|D_k)^*\left[\tau_1^*([{\Cal Z_{\ell_\zeta - k}}]/\alpha)
+ \tau_2^*([{\Cal Z_{k}}]/\alpha) - a \lambda_k\right]\\
&p_k^*\nu^{(k)}|D_k = {1 \over 4} (p_k|D_k)^* \left [
4 \tau_1^*([{\Cal Z_{\ell_\zeta - k}}]/x) + 4 \tau_2^*([{\Cal Z_{k}}]/x)
 - \lambda_k^2 \right ]. \\
\endalign$$
\endstatement
\proof Note that $\Cal U^{(k)}|X \times E_\zeta^{\ell_\zeta - k, k}
= \Cal V_k$, where the sheaf $\Cal V_k$ is constructed by Proposition 2.8
and sits in the exact sequence:
$$0 \to \pi _1^*\scrO _X(F )\otimes \rho_k ^*\pi
_{1,2}^*I_{\Cal Z_{\ell_\zeta - k}} \otimes \pi_2^*\lambda_k
\to \Cal V_k \to \pi _1^*\scrO _X(\Delta -F )
\otimes \rho_k ^*\pi _{1,3}^*I_{\Cal Z_{k}} \to  0.$$
Thus, $c_1(\Cal V_k) = \pi _1^*\Delta + \pi_2^*\lambda_k$ and
$(\operatorname{Id} \times p_{k})^*c_1(\Cal U^{(k)})|(X \times D_k)
= \pi _1^*\Delta + (p_k|D_k)^*\lambda_k$. Moreover,
$c_2(\Cal V_k) = \rho_k ^*\pi_{1,2}^*[{\Cal Z_{\ell_\zeta - k}}]
+ \rho_k ^*\pi_{1,3}^*[{\Cal Z_{k}}] + (\pi _1^*F + \pi_2^*\lambda_k)
\cdot \pi _1^*(\Delta - F)$.
Since $p_k^*\mu^{(k)}(\alpha)|D_k =
(p_k|D_k)^*[\mu^{(k)}(\alpha)|E_\zeta^{\ell_\zeta - k, k}]
= (p_k|D_k)^*[-{1 \over 4} p_1(\Cal V_k)/\alpha]$, we have
$$p_k^*\mu^{(k)}(\alpha)|D_k =
(p_k|D_k)^*\left[\tau_1^*([{\Cal Z_{\ell_\zeta - k}}]/\alpha)
+ \tau_2^*([{\Cal Z_{k}}]/\alpha) - a \lambda_k\right].$$
Similarly, $p_k^*\nu^{(k)}|D_k = {1 \over 4} (p_k|D_k)^* \left [
4 \tau_1^*([{\Cal Z_{\ell_\zeta - k}}]/x) + 4 \tau_2^*([{\Cal Z_{k}}]/x)
 - \lambda_k^2 \right ]$.
\endproof

It follows from the work of Morgan \cite{25} and Li \cite{21},
together with unpublished work of Morgan, that
$D^X _{w,p}(\Cal C_\pm)(\alpha^d) = \delta(\Delta) \cdot \mu_\pm(\alpha)^d$ and
$$D^X _{w,p}(\Cal C_\pm)(\alpha^{d - 2}, x)
= \delta(\Delta) \cdot \mu_\pm(\alpha)^{d - 2} \cdot \nu_\pm$$
where $d = -p - 3$,
$\delta(\Delta) = (-1)^{{{(\Delta^2 + \Delta \cdot K_X)}/2}}$
is the difference between the complex orientation and the standard orientation
on the instanton moduli space (see \cite{6}), and $x \in H_0(X; \Zee)$
is the natural generator. Strictly speaking, their methods only handle the case
of $D^X _{w,p}(\Cal C_\pm)(\alpha^d)$. To handle the case of $D^X _{w,p}(\Cal
C_\pm)(\alpha^{d - 2}, x) $, one needs a blowup formula in algebraic geometry,
which has been established by O'Grady (unpublished).  To compute the
differences
$$\mu_+(\alpha)^d - \mu_-(\alpha)^d \quad \text{and} \quad
\mu_+(\alpha)^{d - 2} \cdot \nu_+ - \mu_-(\alpha)^{d - 2} \cdot \nu_-,$$
we need to know how $\mu^{(k)}(\alpha)$ and $\mu^{(k - 1)}(\alpha)$
are related, and also how $\nu^{(k)}$ and $\nu^{(k - 1)}$ are related.
The following lemma handles this problem.

\lemma{5.3} For $\alpha \in H_2(X; \Zee)$ and the natural generator
$x \in H_0(X; \Zee)$, we have
$$\align
&q_{k - 1}^*\mu^{(k-1)}(\alpha) = p_k^*\mu^{(k)}(\alpha) - aD_k\\
&q_{k - 1}^*\nu^{(k - 1)} = p_k^*\nu^{(k)} - {1 \over 4}
[D_k^2 + 2(p_k|D_k)^*\lambda_k].\\
\endalign$$
\endstatement
\proof From the construction,
the sheaf $(\operatorname{Id} \times q_{k - 1})^*\Cal U^{(k - 1)}$
on $X \times \widetilde{\frak M}_0^{(k)}$ is
the elementary modification of $(\operatorname{Id} \times p_{k})^*\Cal U^{(k)}$
along the divisor $X \times D_k$, using the surjection from
$(\operatorname{Id} \times p_{k})^*\Cal U^{(k)}$ to the pullback of
$\rho_k^*(\pi_1^*\Cal O_X(\Delta - F) \otimes
\pi_{1, 3}^*I_{\Cal Z_k})$:
$$0 \to (\operatorname{Id} \times q_{k - 1})^*\Cal U^{(k - 1)}
\to (\operatorname{Id} \times p_{k})^*\Cal U^{(k)}$$
$$\to (\operatorname{Id} \times p_{k}|D_k)^* \rho_k^*(\pi_1^*\Cal O_X(\Delta
 - F) \otimes \pi_{1, 3}^*I_{\Cal Z_k}) \to 0$$
where $(2F - \Delta) = \zeta$ and $\pi_1$ is the natural projection
$X \times H_{\ell_\zeta - k} \times H_{k} \to X$. Note that
$(\operatorname{Id} \times p_{k}|D_k)^* \rho_k^*(\pi_1^*\Cal O_X(\Delta - F)
\otimes \pi_{1, 3}^*I_{\Cal Z_k})$ is a sheaf supported on $X \times D_k$,
and that its first and second Chern classes are equal to $(X \times D_k)$
and $(X \times D_k^2) - \pi_1^*(\Delta - F) \cdot (X \times D_k)$ respectively.
It follows that
$$\align
&(\operatorname{Id} \times q_{k - 1})^* c_1(\Cal U^{(k - 1)}) =
(\operatorname{Id} \times p_{k})^* c_1(\Cal U^{(k)}) - (X \times D_k)\\
&(\operatorname{Id} \times q_{k - 1})^* c_2(\Cal U^{(k - 1)}) =
(\operatorname{Id} \times p_{k})^* c_2(\Cal U^{(k)}) -
(\operatorname{Id} \times p_{k})^* c_1(\Cal U^{(k)}) \cdot (X \times D_k)\\
&\qquad\qquad\qquad\qquad\qquad\qquad\qquad
+ \pi_1^*(\Delta - F) \cdot (X \times D_k).\\
\endalign$$
By Lemma 5.2, $(\operatorname{Id} \times p_{k})^*c_1(\Cal U^{(k)})
\cdot (X \times D_k) = (\Delta \times D_k) + (X \times (p_k|D_k)^*\lambda_k)$.
Thus,
$$\align
(\operatorname{Id} \times q_{k - 1})^* p_1(\Cal U^{(k - 1)})
&= (\operatorname{Id} \times p_{k})^* p_1(\Cal U^{(k)}) + (X \times D_k^2)
 - 4 (\Delta - F) \times D_k\\
&\qquad\qquad
+ 2 (\operatorname{Id} \times p_{k})^*c_1(\Cal U^{(k)}) \cdot (X \times D_k)\\
&= (\operatorname{Id} \times p_{k})^* p_1(\Cal U^{(k)})
+ 2(2F - \Delta) \times D_k\\
&\qquad\qquad + X \times [D_k^2 + 2(p_k|D_k)^*\lambda_k].\\
\endalign$$
Now the conclusions follow from some straightforward calculations.
\endproof

In the next two theorems, we will give formulas for the differences
$[\mu_+(\alpha)]^d - [\mu_-(\alpha)]^d$ and
$[\mu_+(\alpha)]^{d - 2} \cdot \nu_+ - [\mu_-(\alpha)]^{d - 2} \cdot \nu_-$
in terms of the intersections in $H_{\ell_\zeta - k} \times H_k$
and the Segre classes of the vector bundles
$\Cal E_{\zeta}^{\ell_\zeta - k, k} \oplus
(\Cal E_{-\zeta}^{k, \ell_\zeta -k})\spcheck$ on $H_{\ell_\zeta - k} \times
H_k$,
where $k = 0, 1, \ldots, \ell_\zeta$.
The arguments are a little complicated, but the idea is that we are trying
to get rid of the exceptional divisors $D_k$ as well as
the Chern classes of the tautological line bundles $\xi_k$ and $\lambda_k$.

\theorem{5.4} Let $\zeta$ define a wall of type $(w, p)$,
and $d = (-p - 3)$. For $\alpha \in H_2(X; \Zee)$,
put $a = (\zeta \cdot \alpha)/2$. Then,
$[\mu_+(\alpha)]^d - [\mu_-(\alpha)]^d$ is equal to
$$\sum_{j = 0}^{2\ell_\zeta}~ {d \choose j} \cdot
(-1)^{h(\zeta) + \ell_\zeta + j} \cdot a^{d - j}
\cdot \sum_{k = 0}^{\ell_\zeta}~
([{\Cal Z_{\ell_\zeta - k}}]/\alpha + [{\Cal Z_{k}}]/\alpha)^j
\cdot s_{2\ell_\zeta - j}(\Cal E_{\zeta}^{\ell_\zeta - k, k} \oplus
(\Cal E_{-\zeta}^{k, \ell_\zeta -k})\spcheck).$$
\endstatement
\proof By Lemma 5.3, we have
$q_{k - 1}^*\mu^{(k-1)}(\alpha) = p_k^*\mu^{(k)}(\alpha) - aD_k$.
Since $p_k$ and $q_{k - 1}$ are birational morphisms,
$[p_k^*\mu^{(k)}(\alpha)]^d = [\mu^{(k)}(\alpha)]^d$ and
$[q_{k - 1}^*\mu^{(k-1)}(\alpha)]^d = [\mu^{(k-1)}(\alpha)]^d$. Thus,
$[\mu^{(k - 1)}(\alpha)]^d - [\mu^{(k)}(\alpha)]^d$ is equal to
$$\align
&\quad \sum_{i = 1}^d~ {d \choose i} \cdot
[p_k^*\mu^{(k)}(\alpha)|D_k]^{d - i} \cdot (-D_k|D_k)^{i - 1} \cdot (-a^i) \\
&= \sum_{i = 1}^d~ {d \choose i} \cdot
[p_k^*\mu^{(k)}(\alpha)|D_k]^{d - i} \cdot \xi_k^{i - 1} \cdot (-a^i). \\
\endalign$$
By Lemma 5.2, $p_k^*\mu^{(k)}(\alpha)|D_k =
(p_k|D_k)^*([{\Cal Z_{\ell_\zeta - k}}]/\alpha
+ [{\Cal Z_{k}}]/\alpha - a \lambda_k)$. So we have
$$\align
&\quad [\mu^{(k - 1)}(\alpha)]^d - [\mu^{(k)}(\alpha)]^d \\
&= \sum_{i = 1}^d~ {d \choose i} \cdot \sum_{j = 0}^{2\ell_\zeta}~
{{d - i} \choose j} \cdot
([{\Cal Z_{\ell_\zeta - k}}]/\alpha + [{\Cal Z_{k}}]/\alpha)^j \cdot
(-a\lambda_k)^{d - i - j} \cdot \xi_k^{i - 1} \cdot (-a^i)\\
&= \sum_{j = 0}^{2\ell_\zeta}~ \sum_{i = 1}^{d - j}~ {d \choose j}
\cdot {{d - j} \choose i} \cdot (-a^{d - j}) \cdot
([{\Cal Z_{\ell_\zeta - k}}]/\alpha + [{\Cal Z_{k}}]/\alpha)^j \cdot
\xi_k^{i - 1} \cdot (-\lambda_k)^{d - i - j}\\
&= \sum_{j = 0}^{2\ell_\zeta}~ {d \choose j} \cdot (-a^{d - j}) \cdot
([{\Cal Z_{\ell_\zeta - k}}]/\alpha + [{\Cal Z_{k}}]/\alpha)^j
\cdot \sum_{i = 1}^{d - j}~ {{d - j} \choose i} \cdot
\xi_k^{i - 1} \cdot (-\lambda_k)^{d - i - j}\\
&= \sum_{j = 0}^{2\ell_\zeta}~ {d \choose j} \cdot (-a^{d - j}) \cdot
([{\Cal Z_{\ell_\zeta - k}}]/\alpha + [{\Cal Z_{k}}]/\alpha)^j
\cdot \sum_{i = 0}^{d - 1 - j}~ {{d - j} \choose {i + 1}} \cdot
\xi_k^{i} \cdot (-\lambda_k)^{d - 1 - j - i}\\
\endalign$$
Now, our formula follows from the following claim by summing
$k$ from $0$ to $\ell_\zeta$.

\claim{}
$$([{\Cal Z_{\ell_\zeta - k}}]/\alpha + [{\Cal Z_{k}}]/\alpha)^j
\cdot \sum_{i = 0}^{d - 1 - j}~ {{d - j} \choose {i + 1}} \cdot
\xi_k^{i} \cdot (-\lambda_k)^{d - 1 - j - i}$$
$$= ([{\Cal Z_{\ell_\zeta - k}}]/\alpha + [{\Cal Z_{k}}]/\alpha)^j
\cdot (-1)^{h(\zeta) + \ell_\zeta + j - 1} \cdot
s_{2\ell_\zeta - j}(\Cal E_{\zeta}^{\ell_\zeta - k, k}
\oplus (\Cal E_{-\zeta}^{k, \ell_\zeta -k})\spcheck).$$
\endstatement
\par\noindent
{\it Proof.} For simplicity, on the exceptional divisor $D_k$, we put
$$\sigma_s = ([{\Cal Z_{\ell_\zeta - k}}]/\alpha + [{\Cal Z_{k}}]/\alpha)^j
\cdot \sum_{i = 0}^s {{s + 1} \choose {i + 1}} \cdot \xi_k^i \cdot
(-\lambda_k)^{s - i}.$$
So we must compute $\sigma_{d - 1 - j}$. Notice the relation
$$\sigma_s + \lambda_k \cdot \sigma_{s - 1}
= ([{\Cal Z_{\ell_\zeta - k}}]/\alpha + [{\Cal Z_{k}}]/\alpha)^j
\cdot (\xi_k - \lambda_k)^s.$$
Thus for $0 \le t \le s$, we have
$$\sigma_s = (-\lambda_k)^t \cdot \sigma_{s - t} +
([{\Cal Z_{\ell_\zeta - k}}]/\alpha + [{\Cal Z_{k}}]/\alpha)^j
\cdot \sum_{i = 0}^{t - 1} (\xi_k - \lambda_k)^{s - i} \cdot (-\lambda_k)^i.$$
Put $s = d - 1 - j$ and $t = s - {N_{-\zeta}} = d - 1 - j - {N_{-\zeta}}$,
where ${N_{-\zeta}} = \ell_{-\zeta} + h(-\zeta) - 1
= \ell_{\zeta} + h(-\zeta) - 1$ as defined in Corollary 2.7. Then,
$\sigma_{d - 1 - j}$ is equal to
$$(-\lambda_k)^{d - 1 - j - {N_{-\zeta}}} \cdot \sigma_{N_{-\zeta}} +
([{\Cal Z_{\ell_\zeta - k}}]/\alpha + [{\Cal Z_{k}}]/\alpha)^j
\cdot \sum_{i = 0}^{d - 2 - j - {N_{-\zeta}}}
(\xi_k - \lambda_k)^{(d - 1 - j) - i} \cdot (-\lambda_k)^i.$$
Since $\dim E_\zeta^{\ell_\zeta - k, k} = d - 1 - {N_{-\zeta}}$,
we see that $(-\lambda_k)^{d - 1 - j - {N_{-\zeta}}} \cdot
\sigma_{N_{-\zeta}}$ is equal to
$$\align
&\quad (-\lambda_k)^{d - 1 - j - {N_{-\zeta}}} \cdot
([{\Cal Z_{\ell_\zeta - k}}]/\alpha + [{\Cal Z_{k}}]/\alpha)^j
\cdot \sum_{i = 0}^{N_{-\zeta}} {{{N_{-\zeta}} + 1} \choose {i + 1}}
\cdot \xi_k^i \cdot (-\lambda_k)^{{N_{-\zeta}} - i}\\
&= (-\lambda_k)^{d - 1 - j - {N_{-\zeta}}} \cdot
([{\Cal Z_{\ell_\zeta - k}}]/\alpha + [{\Cal Z_{k}}]/\alpha)^j
\cdot \xi_k^{N_{-\zeta}}\\
&= ([{\Cal Z_{\ell_\zeta - k}}]/\alpha + [{\Cal Z_{k}}]/\alpha)^j \cdot
(-\lambda_k)^{d - 1 - j - {N_{-\zeta}}} \cdot (\xi_k -
\lambda_k)^{N_{-\zeta}}\\
\endalign$$
since the restriction of $\xi_k$ to a fiber of
$D_k \to  E_\zeta^{\ell_\zeta - k, k}$ is a hyperplane. Therefore,
$$\sigma_{d - 1 - j} =
([{\Cal Z_{\ell_\zeta - k}}]/\alpha + [{\Cal Z_{k}}]/\alpha)^j \cdot
\sum_{i = 0}^{d - 1 - j - {N_{-\zeta}}} (\xi_k - \lambda_k)^{(d - 1 - j) - i}
\cdot (-\lambda_k)^i.$$

Now, we shall simplify $(\xi_k - \lambda_k)^{(d - 1 - j) - i}$.
Since $\xi_k$ is the tautological line bundle on
$D_k = \Pee(\Cal N_k\spcheck)$, the line bundle $(\xi_k \otimes
\lambda_k^{-1})$
is the tautological line bundle on
$$\Pee(\Cal N_k\spcheck \otimes \lambda_k^{-1}) =
\Pee[((\rho_k|E_\zeta^{k, \ell_\zeta - k})\spcheck\Cal E_{-\zeta}^{k,
\ell_\zeta -k})\spcheck].$$
Since ${N_{-\zeta}} + 1$ is the rank of $\Cal E_{-\zeta}^{k, \ell_\zeta -k}$,
it follows that
$$(\xi_k - \lambda_k)^{1 + {N_{-\zeta}}} = -\sum_{j = 1}^{1 + {N_{-\zeta}}}
c_j(\Cal E_{-\zeta}^{k, \ell_\zeta -k})
\cdot (\xi_k - \lambda_k)^{1 + {N_{-\zeta}} - j}.$$
One verifies that in general, for $u' \ge {N_{-\zeta}}$, one has
$$(\xi_k - \lambda_k)^{u'} =
s_{u' - {N_{-\zeta}}}(\Cal E_{-\zeta}^{k, \ell_\zeta -k}) \cdot
(\xi_k - \lambda_k)^{{N_{-\zeta}}} +
O\left((\xi_k - \lambda_k)^{{N_{-\zeta}} - 1}\right)$$
where $s_i(\Cal E_{-\zeta}^{k, \ell_\zeta -k})$ is the $i^{\text{th}}$ Segre
class of $\Cal E_{-\zeta}^{k, \ell_\zeta -k}$. Therefore,
since $(d - 1 - j) - i \ge {N_{-\zeta}}$, we see that
$(\xi_k - \lambda_k)^{(d - 1 - j) - i}$ is equal to
$$s_{d - 1 - j - i - {N_{-\zeta}}}(\Cal E_{-\zeta}^{k, \ell_\zeta -k}) \cdot
(\xi_k - \lambda_k)^{{N_{-\zeta}}}
+ O\left((\xi_k - \lambda_k)^{{N_{-\zeta}} - 1}\right)$$
and that $([{\Cal Z_{\ell_\zeta - k}}]/\alpha + [{\Cal Z_{k}}]/\alpha)^j \cdot
(\xi_k - \lambda_k)^{(d - 1 - j) - i} \cdot (-\lambda_k)^i$ is equal to
$$([{\Cal Z_{\ell_\zeta - k}}]/\alpha + [{\Cal Z_{k}}]/\alpha)^j \cdot
\left[s_{d - 1 - j - i - {N_{-\zeta}}}(\Cal E_{-\zeta}^{k, \ell_\zeta -k})
\cdot (\xi_k - \lambda_k)^{{N_{-\zeta}}}\right] \cdot (-\lambda_k)^i$$
$$= ([{\Cal Z_{\ell_\zeta - k}}]/\alpha + [{\Cal Z_{k}}]/\alpha)^j \cdot
s_{d - 1 - j - i - {N_{-\zeta}}}(\Cal E_{-\zeta}^{k, \ell_\zeta -k})
\cdot (-\lambda_k)^i.$$

Next, we note that $([{\Cal Z_{\ell_\zeta - k}}]/\alpha +
[{\Cal Z_{k}}]/\alpha)^j \cdot
s_{d - 1 - j - i - {N_{-\zeta}}}(\Cal E_{-\zeta}^{k, \ell_\zeta -k})$
is a cycle on $E_{\zeta}^{\ell_\zeta -k, k}$ pulled-back from
$H_{\ell_\zeta -k} \times H_k$. So this term is zero unless
$d - 1 - i - {N_{-\zeta}} \le 2\ell_\zeta$, that is,
$i \ge d - 1 - {N_{-\zeta}} - 2\ell_\zeta$.
Note that by Corollary 2.7, $d - 1 - {N_{-\zeta}} - 2\ell_\zeta = {N_\zeta}$
and ${N_\zeta} + 1 = h(\zeta) + \ell_\zeta$ is the rank of
$\Cal E_{\zeta}^{\ell_\zeta -k, k}$.
Since $\lambda_k$ is the tautological line bundle on
$E_{\zeta}^{\ell_\zeta -k, k} = \Pee((\Cal E_{\zeta}^{\ell_\zeta -k,
k})\spcheck)$, we see as before that
$$\lambda_k^i = s_{i - {N_\zeta}}(\Cal E_{\zeta}^{\ell_\zeta -k, k}) \cdot
\lambda_k^{N_\zeta} + O\left(\lambda_k^{{N_\zeta} - 1}\right).$$

Putting all these together, we conclude that $\sigma_{d - 1 - j}$ is equal to
$$\align
&\quad ([{\Cal Z_{\ell_\zeta - k}}]/\alpha + [{\Cal Z_{k}}]/\alpha)^j \cdot
\sum_{i = {N_\zeta}}^{d - 1 - j - {N_{-\zeta}}}~
s_{d - 1 - j - i - {N_{-\zeta}}}(\Cal E_{-\zeta}^{k, \ell_\zeta -k})
\cdot (-1)^i \cdot s_{i - {N_\zeta}}(\Cal E_{\zeta}^{\ell_\zeta -k, k})\\
&= ([{\Cal Z_{\ell_\zeta - k}}]/\alpha + [{\Cal Z_{k}}]/\alpha)^j
\cdot \sum_{i = 0}^{2\ell_\zeta - j}~ (-1)^{i + {N_\zeta}} \cdot
s_{(2\ell_\zeta - j) - i}(\Cal E_{-\zeta}^{k, \ell_\zeta -k}) \cdot
s_{i}(\Cal E_{\zeta}^{\ell_\zeta -k, k})\\
&= ([{\Cal Z_{\ell_\zeta - k}}]/\alpha + [{\Cal Z_{k}}]/\alpha)^j
\cdot (-1)^{j + {N_\zeta}} \cdot\sum_{i = 0}^{2\ell_\zeta - j}~
s_{(2\ell_\zeta - j) - i}((\Cal E_{-\zeta}^{k, \ell_\zeta -k})\spcheck) \cdot
s_{i}(\Cal E_{\zeta}^{\ell_\zeta -k, k})\\
&= ([{\Cal Z_{\ell_\zeta - k}}]/\alpha + [{\Cal Z_{k}}]/\alpha)^j
\cdot (-1)^{j + {N_\zeta}} \cdot
s_{2\ell_\zeta - j}(\Cal E_{\zeta}^{\ell_\zeta -k, k} \oplus
(\Cal E_{-\zeta}^{k, \ell_\zeta -k})\spcheck) \qed\\
\endalign$$

This completes the proof of the Theorem.
\endproof

For the difference $[\mu_+(\alpha)]^{d - 2} \cdot \nu_+
 - [\mu_-(\alpha)]^{d - 2} \cdot \nu_-$, we have the following.

\theorem{5.5} Let $\zeta$ define a wall of type $(w, p)$,
and $d = -p - 3$. For $\alpha \in H_2(X; \Zee)$,
put $a = (\zeta \cdot \alpha)/2$. Then, $[\mu_+(\alpha)]^{d - 2} \cdot \nu_+
 - [\mu_-(\alpha)]^{d - 2} \cdot \nu_-$ is equal to
$${1 \over 4} \cdot \sum_{j = 0}^{2 \ell_\zeta} {{d - 2} \choose j} \cdot
(-1)^{h(\zeta) + \ell_\zeta - 1 + j} \cdot a^{d - 2 - j} \cdot$$
$$\sum_{k = 0}^{\ell_\zeta}
([{\Cal Z_{\ell_\zeta - k}}]/\alpha + [{\Cal Z_{k}}]/\alpha)^j \cdot
\left [ s_{2 \ell_\zeta - j} -
4 ([{\Cal Z_{\ell_\zeta - k}}] + [{\Cal Z_{k}}])/x \cdot
s_{2 \ell_\zeta - 2 - j} \right ]$$
where $s_i$ stands for the $i^{\text{th}}$ Segre class of
$\Cal E_{\zeta}^{\ell_\zeta - k, k}
\oplus (\Cal E_{-\zeta}^{k, \ell_\zeta -k})\spcheck$.
\endstatement
\proof By Lemma 5.3, we have $q_{k - 1}^*\mu^{(k-1)}(\alpha)
= p_k^*\mu^{(k)}(\alpha) - aD_k$ and
$$q_{k - 1}^*\nu^{(k - 1)}
= p_k^*\nu^{(k)} - {1 \over 4} [D_k^2 + 2 (p_k|D_k)^*\lambda_k].$$
It follows that $[\mu^{(k-1)}(\alpha)]^{d - 2} \cdot \nu^{(k-1)} -
[\mu^{(k)}(\alpha)]^{d - 2} \cdot \nu^{(k)} = I_1 + I_2$ where
$$\align
&I_1 = [\mu^{(k)}(\alpha) - aD_k]^{d - 2} \cdot
{1 \over 4} [-D_k^2 - 2 (p_k|D_k)^*\lambda_k] \\
&\quad = [\mu^{(k)}(\alpha)|D_k + a\xi_k]^{d - 2} \cdot
{1 \over 4} (\xi_k - 2 \lambda_k) \\
&I_2 = \sum_{i = 1}^{d - 2} {{d - 2} \choose i}
\cdot \mu^{(k)}(\alpha)^{d - 2 - i} \cdot (-aD_k)^i \cdot \nu^{(k)} \\
&\quad = \sum_{i = 1}^{d - 2}
{{d - 2} \choose i} \cdot [\mu^{(k)}(\alpha)|D_k]^{d - 2 - i}
\cdot \xi_k^{i - 1} \cdot (-a^i) \cdot (\nu^{(k)}|D_k).\\
\endalign$$

First of all, since $\mu^{(k)}(\alpha)|D_k =
([{\Cal Z_{\ell_\zeta - k}}]/\alpha
+ [{\Cal Z_{k}}]/\alpha - a \lambda_k)$, we see that
$$\align
I_1 &= \left [([{\Cal Z_{\ell_\zeta - k}}]/\alpha
+ [{\Cal Z_{k}}]/\alpha) + a (\xi_k - \lambda_k) \right]^{d - 2}
\cdot {1 \over 4} (\xi_k - 2 \lambda_k)\\
&= {1 \over 4} \sum_{j = 0}^{2 \ell_\zeta} {{d - 2} \choose j} \cdot
a^{d - 2 - j} \cdot
([{\Cal Z_{\ell_\zeta - k}}]/\alpha + [{\Cal Z_{k}}]/\alpha)^j
\cdot (\xi_k - \lambda_k)^{d - 2 - j} \cdot (\xi_k - 2 \lambda_k)\\
&= {1 \over 4} \sum_{j = 0}^{2 \ell_\zeta} {{d - 2} \choose j} \cdot
a^{d - 2 - j} \cdot
([{\Cal Z_{\ell_\zeta - k}}]/\alpha + [{\Cal Z_{k}}]/\alpha)^j \cdot \\
&\quad\quad\quad\quad \cdot \left [ (\xi_k - \lambda_k)^{d - 1 - j} -
\lambda_k \cdot (\xi_k - \lambda_k)^{d - 2 - j} \right ]\\
&= {1 \over 4} \sum_{j = 0}^{2 \ell_\zeta} {{d - 2} \choose j} \cdot
a^{d - 2 - j} \cdot
([{\Cal Z_{\ell_\zeta - k}}]/\alpha + [{\Cal Z_{k}}]/\alpha)^j \cdot \\
&\quad\quad\quad\quad \cdot \left [
s_{d - 1 - j - {N_{-\zeta}}}(\Cal E_{-\zeta}^{k, \ell_\zeta -k}) -
\lambda_k \cdot
s_{d - 2 - j - {N_{-\zeta}}}(\Cal E_{-\zeta}^{k, \ell_\zeta -k}) \right ].\\
\endalign$$

Next, by Lemma 5.2, we have $\nu^{(k)}|D_k = {1 \over 4}
\left [ 4 [{\Cal Z_{\ell_\zeta - k}}]/x +
4 [{\Cal Z_{k}}]/x - \lambda_k^2 \right ]$. Thus,
as in the proof of Theorem 5.4, we can verify that $I_2$ is equal to
$$\align
&\quad {1 \over 4} \left [ 4 [{\Cal Z_{\ell_\zeta - k}}]/x +
4 [{\Cal Z_{k}}]/x - \lambda_k^2 \right ] \cdot
\sum_{i = 1}^{d - 2} {{d - 2} \choose i} \cdot
[\mu^{(k)}(\alpha)|D_k]^{d - 2 - i} \cdot \xi_k^{i - 1} \cdot (-a^i) \\
&= {1 \over 4} \left [ 4 [{\Cal Z_{\ell_\zeta - k}}]/x +
4 [{\Cal Z_{k}}]/x - \lambda_k^2 \right ] \cdot
\sum_{j = 0}^{2 \ell_\zeta} {{d - 2} \choose j} \cdot (-a^{d - 2 - j}) \cdot \\
&\quad\quad\quad \cdot
([{\Cal Z_{\ell_\zeta - k}}]/\alpha + [{\Cal Z_{k}}]/\alpha)^j \cdot
\sum_{i = 0}^{2 \ell_\zeta + {N_\zeta} - 2 - j}
s_{d - 3 - j - i - {N_{-\zeta}}}(\Cal E_{-\zeta}^{k, \ell_\zeta -k}) \cdot
(- \lambda_k)^i\\
&= {1 \over 4} \sum_{j = 0}^{2 \ell_\zeta} {{d - 2} \choose j} \cdot
(-a^{d - 2 - j}) \cdot
([{\Cal Z_{\ell_\zeta - k}}]/\alpha + [{\Cal Z_{k}}]/\alpha)^j \cdot \\
&\cdot \left [ 4 ([{\Cal Z_{\ell_\zeta - k}}] + [{\Cal Z_{k}}])/x \cdot
(-1)^{j + {N_\zeta}} \cdot s'
 - \sum_{i = 0}^{2 \ell_\zeta + {N_\zeta} - 2 - j}
s_{d - 3 - j - i - {N_{-\zeta}}}(\Cal E_{-\zeta}^{k, \ell_\zeta -k}) \cdot
(- \lambda_k)^{i + 2} \right ] \\
&= {1 \over 4} \sum_{j = 0}^{2 \ell_\zeta} {{d - 2} \choose j} \cdot
a^{d - 2 - j} \cdot
([{\Cal Z_{\ell_\zeta - k}}]/\alpha + [{\Cal Z_{k}}]/\alpha)^j \cdot \\
&\cdot \left [ \sum_{i = 0}^{2 \ell_\zeta + {N_\zeta} - 2 - j}
s_{d - 3 - j - i - {N_{-\zeta}}}(\Cal E_{-\zeta}^{k, \ell_\zeta -k}) \cdot
(- \lambda_k)^{i + 2} - 4 ([{\Cal Z_{\ell_\zeta - k}}] +
[{\Cal Z_{k}}])/x \cdot (-1)^{j + {N_\zeta}} \cdot s' \right ] \\
\endalign$$
where $s'$ stands for
$s_{2 \ell_\zeta - 2 - j}(\Cal E_{\zeta}^{\ell_\zeta - k, k}
\oplus (\Cal E_{-\zeta}^{k, \ell_\zeta -k})\spcheck)$. Thus,
$I_1 + I_2$ is equal to
$$\align
&\quad {1 \over 4} \sum_{j = 0}^{2 \ell_\zeta} {{d - 2} \choose j} \cdot
a^{d - 2 - j} \cdot
([{\Cal Z_{\ell_\zeta - k}}]/\alpha + [{\Cal Z_{k}}]/\alpha)^j \cdot \\
&\cdot \left [ \sum_{i = -2}^{2 \ell_\zeta + {N_\zeta} - 2 - j}
s_{d - 3 - j - i - {N_{-\zeta}}}(\Cal E_{-\zeta}^{k, \ell_\zeta -k}) \cdot
(- \lambda_k)^{i + 2} - 4 ([{\Cal Z_{\ell_\zeta - k}}] +
[{\Cal Z_{k}}])/x \cdot (-1)^{j + {N_\zeta}} \cdot s' \right ] \\
&= {1 \over 4} \sum_{j = 0}^{2 \ell_\zeta} {{d - 2} \choose j} \cdot
a^{d - 2 - j} \cdot
([{\Cal Z_{\ell_\zeta - k}}]/\alpha + [{\Cal Z_{k}}]/\alpha)^j \cdot \\
&\quad \cdot \left [ (-1)^{j + {N_\zeta}} \cdot s'' -
4 ([{\Cal Z_{\ell_\zeta - k}}] + [{\Cal Z_{k}}])/x \cdot
(-1)^{j + {N_\zeta}} \cdot s' \right ] \\
&= {1 \over 4} \sum_{j = 0}^{2 \ell_\zeta} {{d - 2} \choose j} \cdot
(-1)^{h(\zeta) + \ell_\zeta - 1 + j} \cdot a^{d - 2 - j} \cdot
([{\Cal Z_{\ell_\zeta - k}}]/\alpha + [{\Cal Z_{k}}]/\alpha)^j \cdot \\
&\quad \cdot
\left [ s'' - 4 ([{\Cal Z_{\ell_\zeta - k}}] + [{\Cal Z_{k}}])/x \cdot
s' \right ] \\
\endalign$$
since ${N_\zeta} = h(\zeta) + \ell_\zeta - 1$, where $s''$ stands for
$s_{2 \ell_\zeta - j}(\Cal E_{\zeta}^{\ell_\zeta - k, k}
\oplus (\Cal E_{-\zeta}^{k, \ell_\zeta -k})\spcheck)$. Letting $k$ run from
$0$ to $\ell_\zeta$, we obtain the desired formula.
\endproof

\par\noindent
{\bf Remark 5.6.} For the sake of convenience, we record here the following
relation among the Chern classes and the Segre classes of a vector bundle:
$$s_n = -c_1 \cdot s_{n - 1} - c_2 \cdot s_{n - 2} - \ldots - c_n$$
with the convention that $s_0 = 1$. We refer to \cite{12} for details.
\medskip

In the next section, using Theorem 5.4 and Theorem 5.5,
we shall compute $[\mu_+(\alpha)]^d - [\mu_-(\alpha)]^d$
and $[\mu_+(\alpha)]^{d - 2} \cdot \nu_+
 - [\mu_-(\alpha)]^{d - 2} \cdot \nu_-$
explicitly when $0 \le \ell_\zeta \le 2$.
In principle, Theorem 5.4 and Theorem 5.5 give formulas for
these differences in terms of certain intersections
in $H_{\ell_\zeta -k} \times H_k$. However, it is difficult
to evaluate these intersection numbers in general. In the following,
we shall compute the term
$$S_j = \sum_{k = 0}^{\ell_\zeta}~
([{\Cal Z_{\ell_\zeta - k}}]/\alpha + [{\Cal Z_{k}}]/\alpha)^j
\cdot s_{2\ell_\zeta - j}(\Cal E_{\zeta}^{\ell_\zeta - k, k} \oplus
(\Cal E_{-\zeta}^{k, \ell_\zeta -k})\spcheck) \eqno (5.7)$$
in the special cases when
$j = 2\ell_\zeta$ and $2\ell_\zeta - 1$. We start with a simple lemma.

\lemma{5.8} Let $\alpha, \beta \in H_2(X; \Zee)$. Then
$$\align
&([\Cal Z_k]/\alpha)^{2k} =
{{(2k)!} \over {2^k \cdot k!}} \cdot (\alpha^2)^k\\
&([\Cal Z_k]/\alpha)^{2k - 1} \cdot ([\Cal Z_k]/\beta) =
{{(2k)!} \over {2^k \cdot k!}} \cdot (\alpha^2)^{k - 1}
\cdot (\alpha \cdot \beta)\\
&([\Cal Z_k]/\alpha)^{2k - 2} \cdot ([\Cal Z_k]/\beta)^2 = \\
&{{(2k - 2)!} \over {2^{k - 1} \cdot (k - 1)!}} \cdot (\alpha^2)^{k - 1}
\cdot \beta^2
+ {{(2k - 2)!} \over {2^{k - 2} \cdot (k - 2)!}} \cdot (\alpha^2)^{k - 2}
\cdot (\alpha \cdot \beta)^2.\\
\endalign$$
\endstatement
\proof The first equality is well-known (see \cite{28} for instance).
The other statements follow from the first one by considering
$$ ([\Cal Z_k]/\alpha + [\Cal Z_k]/\beta)^{2k}= {{(2k)!} \over {2^k \cdot k!}}
\cdot ((\alpha + \beta)^2)^k,$$
and formally equating the terms involving $(2k-1)$ copies of $\alpha$
and one $\beta$ or $(2k-2)$ copies of $\alpha$ and two copies of $\beta$.
\endproof

The next result computes the term (5.7) when $j = 2\ell_\zeta$.

\proposition{5.9} Let $\zeta$ define a wall of type $(w, p)$,
and $\alpha \in H_2(X; \Zee)$. Then,
$$S_{2\ell_\zeta} = \sum_{k = 0}^{\ell_\zeta}~
([{\Cal Z_{\ell_\zeta - k}}]/\alpha + [{\Cal Z_{k}}]/\alpha)^{2\ell_\zeta} =
{{(2\ell_\zeta)!} \over {\ell_\zeta !}} \cdot (\alpha^2)^{\ell_\zeta}.$$
\endstatement
\par\noindent
{\it Proof.} This follows in a straightforward way from Lemma 5.8 (i):
$$\align
&\quad \sum_{k = 0}^{\ell_\zeta}~
([{\Cal Z_{\ell_\zeta - k}}]/\alpha + [{\Cal Z_{k}}]/\alpha)^{2\ell_\zeta}\\
&= \sum_{k = 0}^{\ell_\zeta}~{{2\ell_\zeta} \choose {2k}} \cdot
([{\Cal Z_{\ell_\zeta - k}}]/\alpha)^{2(\ell_\zeta - k)} \cdot
([{\Cal Z_{k}}]/\alpha)^{2k} \\
&= \sum_{k = 0}^{\ell_\zeta}~{{2\ell_\zeta} \choose {2k}} \cdot
\left [{{(2\ell_\zeta - 2k)!} \over {2^{\ell_\zeta - k} \cdot
(\ell_\zeta - k)!}} \cdot (\alpha^2)^{\ell_\zeta - k} \right ] \cdot
\left [{{(2k)!} \over {2^k \cdot k!}} \cdot (\alpha^2)^k \right ]\\
&= \sum_{k = 0}^{\ell_\zeta}~ {\ell_\zeta \choose k} \cdot
{{(2\ell_\zeta)!} \over {2^{\ell_\zeta} \cdot \ell_\zeta !}}
\cdot (\alpha^2)^{\ell_\zeta} \\
&= {{(2\ell_\zeta)!} \over {\ell_\zeta !}} \cdot (\alpha^2)^{\ell_\zeta}\qed\\
\endalign$$

To compute the term (5.7) when $j = (2\ell_\zeta - 1)$,
we study $\Cal E_{-\zeta}^{k, \ell_\zeta -k}$ and
$\Cal E_{\zeta}^{\ell_\zeta -k, k}$,
and evaluate their first Chern classes. We begin with a general lemma.

\lemma{5.10} Let $Z, W$ be codimension $2$ cycles in
a smooth variety $Y$.
\roster

\item"{(i)}" If $Z \subseteq W$, then $Hom(I_W, I_Z) = \Cal O_Y$;

\item"{(ii)}" If $(Z - Z \cap W)$ is open and dense in $Z$,
then $Hom(I_W, I_Z) = I_Z$;

\item"{(iii)}" If $Z$ and $W$ are local complete intersections meeting
properly,
then there is an exact sequence:
$$0 \to  Ext ^1(I_W, I_Z) \to \Cal O _W \otimes \det
N_W \to   \Cal O _{W \cap Z} \otimes \det N_W \to 0$$
where $N_W$ is the normal bundle of $W$ in $Y$;

\item"{(iv)}" Assume that $Z \cap W$ is nowhere dense in $W$ and that
$W$ is smooth at a generic point. Then, as a sheaf on $W$,
$Ext^1(I_W, I_Z)$ is of rank $1$; thus,
$$c_0(Ext^1(I_W, I_Z)) = c_1(Ext^1(I_W, I_Z)) = 0,
\quad c_2(Ext^1(I_W, I_Z)) = -[W].$$
\endroster
\endstatement
\proof (i) Applying the functor $Hom(I_W, \cdot)$ to
the exact sequence
$$0 \to I_Z \to \Cal O_Y \to \Cal O_Z \to 0,$$
we obtain $0 \to Hom(I_W, I_Z) \to  Hom(I_W, \Cal O_Y) = \Cal O_Y$.
Thus, $Hom(I_W, I_Z) = I_U$ for some closed subscheme $U$ of $Y$.
On the other hand, since $Z \subseteq W$,
$$H^0(Y; Hom(I_W, I_Z)) = \Hom(I_W, I_Z) \ne 0.$$
Thus, $U$ must be empty, and $Hom(I_W, I_Z) = \Cal O_Y$.

(ii) As in the proof of (i), $Hom(I_W, I_Z) = I_U$ for
some closed subscheme $U$ of $Y$. Applying the functor
$Hom(\cdot, I_Z)$ to the exact sequence
$$0 \to I_W \to \Cal O_Y \to \Cal O_W \to 0,$$
we get $0 \to I_Z \to  Hom(I_W, I_Z) = I_U \to  Ext^1(\Cal O_W, I_Z)$.
Thus, $U \subseteq Z$; moreover, since $Ext^1(\Cal O_W, I_Z) = 0$
on $(X - W)$, we have $(Z - Z \cap W) = (U - U \cap W)$. So
$$(Z - Z \cap W) \subseteq U \subseteq Z.$$
Since $(Z - Z \cap W)$ is open and dense in $Z$, it follows that $U = Z$.

(iii) We begin with the local identification:
let $R$ be a regular local ring, and let $Z$ and $W$ be two codimension $2$
local complete intersection subschemes of $R$ meeting properly.
Applying the functor $Hom_R(\cdot, I_Z)$ to the Koszul resolution of $W$
$$0 \to R \to R \oplus R \to I_W \to 0$$
gives $I_Z \oplus I_Z \to I_Z \to Ext_R^1(I_W, I_Z) \to 0$.
It follows that $Ext_R^1(I_W, I_Z) = I_Z/(I_Z \cdot I_W)$.
Since $Z$ and $W$ are codimension $2$ local complete intersections
meeting properly, we have $I_Z \cdot I_W = I_Z \cap I_W$.
Thus, $Ext^1_R(I_W, I_Z) \cong I_Z/(I_Z \cap I_W)$,
and we can fit it into an exact sequence
$$0 \to Ext^1_R(I_W, I_Z) \to R/I_W \to R/(I_W + I_Z) \to 0.$$
Here $(I_W + I_Z)$ corresponds to the intersection $W \cap Z$.
The identification of $Ext_R^1(I_W, I_Z)$ and $I_Z/(I_Z \cap I_W)$
is not canonical. Globally we must correct by $\det N_W$.
Thus globally we have an exact sequence:
$$0 \to  Ext ^1(I_W, I_Z) \to \Cal O _W \otimes \det
N_W \to   \Cal O _{W \cap Z} \otimes \det N_W \to 0.$$

(iv) It is clear that $Ext^1(I_W, I_Z)$ is a sheaf supported on $W$.
To show that it has rank $1$ as a sheaf on $W$, it suffices to verify
that it has rank $1$ at a generic point $w$ of $W$.
Since $Z \cap W$ is nowhere dense in $W$ and
$W$ is smooth at a generic point, we may assume that
$w \not \in Z$ and that $w$ is a smooth point of $W$.
Then it follows from (iii) that $Ext^1(I_W, I_Z)$ is of rank $1$ at $w$.
\endproof

\lemma{5.11} Let $Hom =  Hom(I_{\Cal Z_k},
I_{\Cal Z_{\ell_\zeta - k}})$,
$Ext^1 = Ext^1(I_{\Cal Z_k}, I_{\Cal Z_{\ell_\zeta - k}})$,
$\pi_1$ and $\pi_2$ be the projections from
$X \times (H_{\ell_\zeta - k} \times H_k)$ to $X$ and
$(H_{\ell_\zeta - k} \times H_k)$ respectively.
\roster

\item"{(i)}" There exist a row exact sequence and a column exact sequence:
$$\matrix
          &0&\\
          &\downarrow&\\
          &\pi_{2*}(\pi_1^*\Cal O_X(\zeta) \otimes
           \Cal O_{\Cal Z_{\ell_\zeta - k}}) &\\
          &\downarrow&\\
0 \to & R^1\pi_{2*} \left (\pi_1^*\Cal O_X(\zeta) \otimes  Hom \right )
& \to \Cal E_{\zeta}^{\ell_\zeta -k, k} \to
\pi_{2*} \left (\pi_1^*\Cal O_X(\zeta) \otimes  Ext^1 \right ) \to 0; \\
          &\downarrow&\\
          &[\Cal O_{H_{\ell_\zeta - k} \times H_k}]^{\oplus~ h(\zeta)} &\\
          &\downarrow&\\
          &0&\\
\endmatrix$$

\item"{(ii)}" $c_1 \left(R^1\pi_{2*} (\pi_1^*\Cal O_X(\zeta)
\otimes  Hom)\right) = [\Cal Z_{\ell_\zeta - k}]/(\zeta - K_X/2)
+ \pi_{2*}[c_3(\Cal O_{\Cal Z_{\ell_\zeta - k}})]/2$;

\item"{(iii)}" $c_1 \left(\pi_{2*} (\pi_1^*\Cal O_X(\zeta) \otimes
Ext^1 )\right ) = [\Cal Z_k]/(\zeta - K_X/2)
+ \pi_{2*}[c_3(Ext^1)]/2$.

\endroster
\endstatement
\proof (i) Note that the bundle $\Cal E_\zeta ^{\ell_\zeta -k, k}$
is defined as
$$Ext^1_{\pi _2}(\pi_1^*\scrO _X(\Delta - F) \otimes I_{\Cal Z_k},
\pi _1^*\scrO _X(F)\otimes I_{\Cal Z_{\ell_\zeta -k}})
= Ext^1_{\pi _2}(I_{\Cal Z_k},
\pi _1^*\scrO _X(\zeta) \otimes I_{\Cal Z_{\ell_\zeta -k}}).$$
Since $R^2\pi_{2*}(\pi_1^*\Cal O_X(\zeta) \otimes  Hom) = 0$,
the row exact sequence follows from standard facts about relative Ext sheaves.
To see the column exact sequence, we use Lemma 5.10 (ii)
and apply the functor $\pi_{2*}$ to the exact sequence
$$0 \to \pi _1^*\scrO _X(\zeta) \otimes I_{\Cal Z_{\ell_\zeta -k}} \to
\pi _1^*\scrO _X(\zeta) \to \pi _1^*\scrO _X(\zeta) \otimes
\Cal O_{\Cal Z_{\ell_\zeta -k}} \to 0.$$

(ii) Note that $Hom = I_{\Cal Z_{\ell_\zeta -k}}$ and that
$R^i\pi_{2*}(\pi_1^*\Cal O_X(\zeta) \otimes  Hom) = 0$
for $i = 0, 2$. By the Grothendieck-Riemann-Roch Theorem, we have
$$\align
&\quad -\ch \left (R^1\pi_{2*}(\pi_1^*\Cal O_X(\zeta)
\otimes  Hom)\right )\\
&= \pi_{2*}\left (\ch    (\pi_1^*\Cal O_X(\zeta)
\otimes I_{\Cal Z_{\ell_\zeta -k}}) \cdot \pi_1^*\Todd (T_X) \right )\\
&= \pi_{2*}\left (\pi_1^*\ch    (\Cal O_X(\zeta)) \cdot
\ch    (I_{\Cal Z_{\ell_\zeta -k}})
\cdot \pi_1^*\Todd (T_X) \right ).\\
\endalign$$
Now, the conclusion follows by comparing the degree $1$ terms and
by the fact that
$$\ch    (I_{\Cal Z_{\ell_\zeta -k}}) =
1 - \ch    (\Cal O_{\Cal Z_{\ell_\zeta - k}})
= 1 - [\Cal Z_{\ell_\zeta - k}] -
{c_3(\Cal O_{\Cal Z_{\ell_\zeta - k}}) \over 2}
+ (\text{terms with degree} \ge 4).$$

(iii) We have $R^i\pi_{2*}(\pi_1^*\Cal O_X(\zeta) \otimes  Ext^1) = 0$
for $i = 1, 2$. By Lemma 5.10 (iv),
$$\ch    (Ext^1) = [\Cal Z_k] +
{c_3(Ext^1) \over 2} + (\text{terms with degree} \ge 4).$$
Again, using the Grothendieck-Riemann-Roch Theorem, we obtain
$$\align
&\quad \ch    \left (\pi_{2*}(\pi_1^*\Cal O_X(\zeta)
\otimes  Ext^1)\right )\\
&= \pi_{2*}\left (\ch    (\pi_1^*\Cal O_X(\zeta)
\otimes  Ext^1) \cdot \pi_1^*\Todd (T_X) \right )\\
&= \pi_{2*}\left (\pi_1^*\ch    (\Cal O_X(\zeta)) \cdot
\ch    (Ext^1) \cdot \pi_1^*\Todd (T_X) \right ).\\
\endalign$$
Then, our conclusion follows by comparing the degree $1$ terms.
\endproof

Now, we can compute the term (5.7) for $j = 2\ell_\zeta - 1$.

\proposition{5.12} Let $\alpha \in H_2(X; \Zee)$ and
$a = (\zeta \cdot \alpha)/2$. Then,
$$\align
S_{2\ell_\zeta - 1} &= \sum_{k = 0}^{\ell_\zeta}~
([{\Cal Z_{\ell_\zeta - k}}]/\alpha + [{\Cal Z_{k}}]/\alpha)^{2\ell_\zeta - 1}
\cdot s_1(\Cal E_{\zeta}^{\ell_\zeta - k, k} \oplus
(\Cal E_{-\zeta}^{k, \ell_\zeta -k})\spcheck) \\
&= (-4) \cdot {{(2\ell_\zeta)!} \over {\ell_\zeta !}} \cdot
(\alpha^2)^{\ell_\zeta - 1} \cdot a.\\
\endalign$$
\endstatement
\par\noindent
{\it Proof.} By the symmetry between $k$ and $(\ell_\zeta - k)$,
we see that $S_{2\ell_\zeta - 1}$ is equal to
$$\sum_{k = 0}^{\ell_\zeta}~
([{\Cal Z_{\ell_\zeta - k}}]/\alpha + [{\Cal Z_{k}}]/\alpha)^{2\ell_\zeta - 1}
\cdot {{s_1(\Cal E_{\zeta}^{\ell_\zeta - k, k} \oplus
(\Cal E_{-\zeta}^{k, \ell_\zeta -k})\spcheck) +
s_1(\Cal E_{\zeta}^{k, \ell_\zeta -k} \oplus
(\Cal E_{-\zeta}^{\ell_\zeta - k, k})\spcheck)} \over 2}.$$
 From Lemma 5.11, we conclude that $c_1(\Cal E_{\zeta}^{\ell_\zeta - k, k})$
is equal to
$$([\Cal Z_{\ell_\zeta - k}] + [\Cal Z_k])/(\zeta - K_X/2) +
{{\pi_{2*}[c_3(\Cal O_{\Cal Z_{\ell_\zeta - k}}) +
c_3(Ext^1(I_{\Cal Z_k}, I_{\Cal Z_{\ell_\zeta - k}}))]} \over 2}.$$
Since $s_1(\Cal E_{\zeta}^{\ell_\zeta - k, k} \oplus
(\Cal E_{-\zeta}^{k, \ell_\zeta -k})\spcheck) =
c_1(\Cal E_{-\zeta}^{k, \ell_\zeta -k}) -
c_1(\Cal E_{\zeta}^{\ell_\zeta - k, k})$, we see that
$${{s_1(\Cal E_{\zeta}^{\ell_\zeta - k, k} \oplus
(\Cal E_{-\zeta}^{k, \ell_\zeta -k})\spcheck) +
s_1(\Cal E_{\zeta}^{k, \ell_\zeta -k} \oplus
(\Cal E_{-\zeta}^{\ell_\zeta - k, k})\spcheck)} \over 2} =
(-2) \cdot ([\Cal Z_{\ell_\zeta - k}] + [\Cal Z_k])/\zeta$$
where the $c_3$'s are cancelled out. Therefore, by Lemma 5.8,
$$\align
S_{2\ell_\zeta - 1}&= \sum_{k = 0}^{\ell_\zeta}~
([{\Cal Z_{\ell_\zeta - k}}]/\alpha + [{\Cal Z_{k}}]/\alpha)^{2\ell_\zeta - 1}
\cdot (-2) \cdot ([\Cal Z_{\ell_\zeta - k}]/\zeta + [\Cal Z_k]/\zeta)\\
&= (-2) \cdot \sum_{k = 0}^{\ell_\zeta}~ [
{{2\ell_\zeta - 1} \choose {2k}} \cdot
([{\Cal Z_{\ell_\zeta - k}}]/\alpha)^{2\ell_\zeta - 2k - 1}
\cdot ([{\Cal Z_{k}}]/\alpha)^{2k} \cdot [\Cal Z_{\ell_\zeta - k}]/\zeta\\
&\quad\quad\quad\quad + {{2\ell_\zeta - 1} \choose {2k - 1}} \cdot
([{\Cal Z_{\ell_\zeta - k}}]/\alpha)^{2\ell_\zeta - 2k}
\cdot ([{\Cal Z_{k}}]/\alpha)^{2k - 1} \cdot [\Cal Z_k]/\zeta]\\
&= (-4) \cdot \sum_{k = 1}^{\ell_\zeta}~
{{2\ell_\zeta - 1} \choose {2k - 1}} \cdot
([{\Cal Z_{\ell_\zeta - k}}]/\alpha)^{2\ell_\zeta - 2k}
\cdot ([{\Cal Z_{k}}]/\alpha)^{2k - 1} \cdot [\Cal Z_k]/\zeta\\
&= (-4) \cdot {{(2\ell_\zeta)!} \over {\ell_\zeta !}} \cdot
(\alpha^2)^{\ell_\zeta - 1} \cdot a \qed\\
\endalign$$

It is possible, but far more complicated, to compute (5.7) for
$j = 2\ell_\zeta - 2$.

Next, we shall draw some consequences from our previous computations.
Recall that $q_X$ denotes the intersection form of $X$,
and that
$$\delta(\Delta) = (-1)^{{{\Delta^2 + \Delta \cdot K_X} \over 2}}$$
is the difference between the complex orientation and the standard orientation
on the instanton moduli space (see \cite{6}). Theorem 5.13 below has
already been obtained by Kotschick and Morgan \cite{18} for
any smooth $4$-manifold with $b_2^+ = 1$.

\theorem{5.13} Let $\zeta$ define a wall of type $(w, p)$,
and $d = -p - 3$. Then,
$$[\mu_+(\alpha)]^d - [\mu_-(\alpha)]^d \equiv
(-1)^{h(\zeta) + \ell_\zeta} \cdot
{{d!} \over {\ell_\zeta! \cdot (d - 2\ell_\zeta)!}}
\cdot a^{d - 2\ell_\zeta} \cdot (\alpha^2)^{\ell_\zeta}
\pmod {a^{d - 2\ell_\zeta + 2}}$$
for $\alpha \in H_2(X; \Zee)$, where $a = (\zeta \cdot \alpha)/2$.
In other words,
$$\delta^X_{w, p}(\Cal C_-, \Cal C_+) \equiv
\delta(\Delta) \cdot (-1)^{h(\zeta) + \ell_\zeta} \cdot
{{d!} \over {\ell_\zeta! \cdot (d - 2\ell_\zeta)!}}
\cdot \left ( {\zeta\over 2} \right )^{d - 2\ell_\zeta} \cdot q_X^{\ell_\zeta}
\pmod {\zeta^{d - 2\ell_\zeta + 2}}.$$
\endstatement
\par\noindent
{\it Proof.} By Theorem 5.4 and our notation (5.7), we have
$$[\mu_+(\alpha)]^d - [\mu_-(\alpha)]^d \equiv
\sum_{j = 2\ell_\zeta - 1}^{2\ell_\zeta}~ {d \choose j} \cdot
(-1)^{h(\zeta) + \ell_\zeta + j} \cdot a^{d - j} \cdot S_j
\pmod {a^{d - 2\ell_\zeta + 2}}.$$
By Proposition 5.12, $S_{2\ell_\zeta - 1}$ is divisible by $a$. Therefore,
$$[\mu_+(\alpha)]^d - [\mu_-(\alpha)]^d \equiv
{d \choose {2\ell_\zeta}} \cdot
(-1)^{h(\zeta) + \ell_\zeta} \cdot a^{d - 2\ell_\zeta} \cdot S_{2\ell_\zeta}
\pmod {a^{d - 2\ell_\zeta + 2}}.$$
Now, our conclusion follows from Proposition 5.9 and the fact that
$$\gamma_{\pm}(\alpha^d) = \delta(\Delta) \cdot \mu_{\pm}(\alpha)^d.\qed$$

The following is proved by using a similar method.

\theorem{5.14} Let $\zeta$ define a wall of type $(w, p)$.
For $\alpha \in H_2(X; \Zee)$, let $a = (\zeta \cdot \alpha)/2$.
Then, modulo $a^{d - 2\ell_\zeta}$, $[\mu_+(\alpha)]^{d - 2} \cdot \nu_+ -
[\mu_-(\alpha)]^{d - 2} \cdot \nu_-$ is equal to
$${1 \over 4} \cdot (-1)^{h(\zeta) + \ell_\zeta - 1} \cdot
{{(d - 2)!} \over {\ell_\zeta! \cdot (d - 2 - 2\ell_\zeta)!}}
\cdot a^{d - 2 - 2\ell_\zeta} \cdot (\alpha^2)^{\ell_\zeta}.$$
\endstatement
\par\noindent
{\it Proof.} By Theorem 5.5, $[\mu_+(\alpha)]^{d - 2} \cdot \nu_+ -
[\mu_-(\alpha)]^{d - 2} \cdot \nu_-$ is equal to
$${1 \over 4} \cdot \sum_{j = 2 \ell_\zeta - 1}^{2 \ell_\zeta}
{{d - 2} \choose j} \cdot (-1)^{h(\zeta) + \ell_\zeta - 1 + j}
\cdot a^{d - 2 - j} \cdot S_j$$
modulo $a^{d - 2\ell_\zeta}$, where $S_j$ is the notation introduced in (5.7).
By Proposition 5.12, $S_{2 \ell_\zeta - 1}$ is divisible by $a$;
by Proposition 5.9, we have
$$S_{2\ell_\zeta} = {{(2\ell_\zeta)!} \over {\ell_\zeta !}}
\cdot (\alpha^2)^{\ell_\zeta}.$$
Therefore, modulo $a^{d - 2\ell_\zeta}$,
$[\mu_+(\alpha)]^{d - 2} \cdot \nu_+ -
[\mu_-(\alpha)]^{d - 2} \cdot \nu_-$ is equal to
$${1 \over 4} \cdot (-1)^{h(\zeta) + \ell_\zeta - 1} \cdot
{{(d - 2)!} \over {\ell_\zeta! \cdot (d - 2 - 2\ell_\zeta)!}}
\cdot a^{d - 2 - 2\ell_\zeta} \cdot (\alpha^2)^{\ell_\zeta}. \qed$$

\section{6. The formulas when $\ell_\zeta = 0, 1, 2$.}

In this section, we shall compute
$[\mu_+(\alpha)]^d - [\mu_-(\alpha)]^d$ and
$$[\mu_+(\alpha)]^{d - 2} \cdot \nu_+
 - [\mu_-(\alpha)]^{d - 2} \cdot \nu_-$$
by assuming that $\ell_\zeta = 0, 1, 2$.
Our first result, Theorem 6.1 below, was first obtained by Mong and
Kotschick \cite{17}.

\theorem{6.1} Let $\zeta$ define a wall of type $(w, p)$ with
$\ell_\zeta = 0$. Then,
$$[\mu_+(\alpha)]^d - [\mu_-(\alpha)]^d = (-1)^{h(\zeta)}
\cdot \left ({{\zeta \cdot \alpha} \over 2} \right )^d$$
for $\alpha \in H_2(X; \Zee)$. In other words,
$\delta^X_{w, p}(\Cal C_-, \Cal C_+) =
\delta(\Delta) \cdot (-1)^{h(\zeta)} \cdot ({\zeta}/2)^d$.
\endstatement
\par\noindent
{\it Proof.} There are two cases: $h(\zeta)> 0$ and $h(\zeta)= 0$.
In the first case when $h(\zeta)> 0$,
the formula follows immediately from Theorem 5.4. In the second case
when $h(\zeta)= 0$, we must have $\zeta^2 = p$
and $\zeta \cdot K_X = \zeta^2 + 2 = p + 2$ by Corollary 2.7.
Then $\frak M_+$ consists of $\frak M_-$ and an additional connected component
$E_{-\zeta}^{0, 0} \cong \Pee^{-p - 3}$. We have constructed
a universal sheaf $\Cal U$ over $X \times E_{-\zeta}^{0, 0}$:
$$0 \to \pi_1^*\Cal O_X(\Delta - F) \otimes \pi_2^*\lambda \to
\Cal U \to \pi_1^*\Cal O_X(F) \to 0$$
where $F$ is the unique divisor satisfying $(2F - \Delta) = \zeta$,
$\lambda$ is the line bundle corresponding to a hyperplane in
$E_{-\zeta}^{0, 0} \cong \Pee^{-p - 3}$,
and $\pi_1$ and $\pi_2$ are the natural projections of
$X \times E_{-\zeta}^{0, 0}$. Thus for $\alpha \in H_2(X; \Zee)$, we have
$$\mu_+(\alpha) = \mu_-(\alpha) - {1 \over 4} \cdot p_1(\Cal U)/\alpha
= \mu_-(\alpha) + a \lambda$$
where $a = ({\zeta \cdot \alpha})/2$. Since $h(\zeta) = 0$, we conclude that
$$\mu_+(\alpha)^d = \mu_-(\alpha)^d +
\left ({{\zeta \cdot \alpha} \over 2} \right )^d
= \mu_-(\alpha)^d + (-1)^{h(\zeta)} \cdot
\left ({{\zeta \cdot \alpha} \over 2} \right )^d. \qed$$

The proof of the next result is similar to the proof of Theorem 6.1.

\theorem{6.2} Let $\zeta$ define a wall of type $(w, p)$ with
$\ell_\zeta = 0$, let $d = -p - 3$. Then, for $\alpha \in H_2(X; \Zee)$,
we have
$$[\mu_+(\alpha)]^{d - 2} \cdot \nu_+ - [\mu_-(\alpha)]^{d - 2} \cdot \nu_-
= {1 \over 4} \cdot (-1)^{h(\zeta) - 1} \cdot
\left ({{\zeta \cdot \alpha} \over 2} \right)^{d - 2}. \qed$$
\endstatement

Next, we shall study the difference $\delta^X_{w, p}(\Cal C_-, \Cal C_+)$
when $\ell_\zeta = 1$. In this case, we have to know (5.7)
for $j = 2, 1, 0$. In view of Propositions 5.9 and 5.12,
it suffices to calculate (5.7) for $j = 0$. The following lemma
deals with this.

\lemma{6.3} Let $\zeta$ define a wall of type $(w, p)$ with
$\ell_\zeta = 1$. Then
$$S_0 = \sum_{k = 0}^1~ s_2(\Cal E_{\zeta}^{1 - k, k} \oplus
(\Cal E_{-\zeta}^{k, 1 - k})\spcheck) = (6 \zeta^2 + 2K_X^2).$$
\endstatement
\par\noindent
{\it Proof.} First, we compute the Chern classes of $\Cal E_{\zeta}^{1, 0}$.
Let notations be  as in Lemma 5.11,
and set $\ell_\zeta = 1$ and $k = 0$ in Lemma 5.11.
Then $Ext^1 = 0$. Since $(H_{\ell_\zeta - k} \times H_k) = X$,
the codimension $2$ cycle $\Cal Z_1$ is exactly the diagonal in
$X \times (H_{\ell_\zeta - k} \times H_k) = X \times X$. Thus,
$\pi_{2*}(\pi_1^*\Cal O_X(\zeta) \otimes \Cal O_{\Cal Z_{\ell_\zeta - k}})
= \Cal O_X(\zeta)$. By Lemma 5.11 (i), the bundle $\Cal E_{\zeta}^{1, 0}$
sits in an exact sequence:
$$0 \to \Cal O_X(\zeta) \to \Cal E_{\zeta}^{1, 0} \cong
R^1\pi_{2*} \left (\pi_1^*\Cal O_X(\zeta) \otimes  Hom \right )
\to \Cal O_X^{\oplus~ h(\zeta)} \to 0.$$
Thus, $c_1(\Cal E_{\zeta}^{1, 0}) = \zeta$ and
$c_2(\Cal E_{\zeta}^{1, 0}) = 0$.

Next, we compute the Chern classes of $\Cal E_{\zeta}^{0, 1}$.
Let $\ell_\zeta = 1$ and $k = 1$ in Lemma 5.11. Then,
$Ext^1 = \det (N)$ where $N$ is the normal bundle
of $\Cal Z_1$ in $X \times X$. Thus,
$$\pi_{2*} \left (\pi_1^*\Cal O_X(\zeta) \otimes  Ext^1 \right )
= \Cal O_X(\zeta - K_X).$$
By Lemma 5.11 (i), the bundle $\Cal E_{\zeta}^{0, 1}$
sits in an exact sequence:
$$0 \to \Cal O_X^{\oplus~ h(\zeta)} \to \Cal E_{\zeta}^{0, 1}
\to \Cal O_X(\zeta - K_X) \to 0.$$
Thus, $c_1(\Cal E_{\zeta}^{0, 1}) = \zeta - K_X$ and
$c_2(\Cal E_{\zeta}^{0, 1}) = 0$. Replacing $\zeta$ by $-\zeta$ gives
$c_1(\Cal E_{-\zeta}^{0, 1}) = -\zeta - K_X$
and $c_2(\Cal E_{-\zeta}^{0, 1}) = 0$.
It follows that $c_1(\Cal E_{\zeta}^{1, 0} \oplus
(\Cal E_{-\zeta}^{0, 1})\spcheck) = 2\zeta + K_X$ and that
$$c_2(\Cal E_{\zeta}^{1, 0} \oplus (\Cal E_{-\zeta}^{0, 1})\spcheck)
= \zeta \cdot (\zeta + K_X) = \zeta^2 + \zeta \cdot K_X.$$
So we conclude that the Segre class
$s_2(\Cal E_{\zeta}^{1, 0} \oplus (\Cal E_{-\zeta}^{0, 1})\spcheck)$ is equal
to
$$c_1(\Cal E_{\zeta}^{1, 0} \oplus (\Cal E_{-\zeta}^{0, 1})\spcheck)^2
 - c_2(\Cal E_{\zeta}^{1, 0} \oplus (\Cal E_{-\zeta}^{0, 1})\spcheck)
= 3\zeta^2 + 3 \zeta \cdot K_X + K_X^2.$$
Replacing $\zeta$ by $-\zeta$ gives
$s_2(\Cal E_{-\zeta}^{1, 0} \oplus (\Cal E_{\zeta}^{0, 1})\spcheck) =
3\zeta^2 - 3 \zeta \cdot K_X + K_X^2$. Therefore,
$$\align
S_0 &= \sum_{k = 0}^1~ s_2(\Cal E_{\zeta}^{1 - k, k} \oplus
(\Cal E_{-\zeta}^{k, 1 - k})\spcheck)\\
&= s_2(\Cal E_{\zeta}^{1, 0} \oplus (\Cal E_{-\zeta}^{0, 1})\spcheck)
+ s_2(\Cal E_{\zeta}^{0, 1} \oplus (\Cal E_{-\zeta}^{1, 0})\spcheck)\\
&= s_2(\Cal E_{\zeta}^{1, 0} \oplus (\Cal E_{-\zeta}^{0, 1})\spcheck)
+ s_2((\Cal E_{\zeta}^{0, 1})\spcheck \oplus \Cal E_{-\zeta}^{1, 0})\\
&= (3\zeta^2 + 3 \zeta \cdot K_X + K_X^2) +
(3\zeta^2 - 3 \zeta \cdot K_X + K_X^2)\\
&= 6 \zeta^2 + 2K_X^2.\qed\\
\endalign$$

Now we can compute the difference $\delta^X_{w, p}(\Cal C_-, \Cal C_+)$
when $\ell_\zeta = 1$.

\theorem{6.4} Let $\zeta$ define a wall of type $(w, p)$ with
$\ell_\zeta = 1$. Then,
$$[\mu_+(\alpha)]^d - [\mu_-(\alpha)]^d = (-1)^{h(\zeta) + 1} \cdot
\left \{ d(d - 1) \cdot a^{d - 2} \cdot \alpha^2 +
(2K_X^2 + 2d + 6) \cdot a^d \right \}$$
for $\alpha \in H_2(X; \Zee)$, where $a = (\zeta \cdot \alpha)/2$.
In other words, $\delta^X_{w, p}(\Cal C_-, \Cal C_+)$ is equal to
$$\delta(\Delta) \cdot (-1)^{h(\zeta) + 1} \cdot
\left \{ d(d - 1) \cdot \left(\zeta \over 2 \right)^{d - 2} \cdot q_X +
(2K_X^2 + 2d + 6) \cdot \left ({\zeta \over 2} \right )^d \right \}.$$
\endstatement
\par\noindent
{\it Proof.} From 5.4, 5.9, 5.12, and 6.3, we conclude that
$$\align
&\quad [\mu_+(\alpha)]^d - [\mu_-(\alpha)]^d\\
&= (-1)^{h(\zeta) + 1} \cdot d(d - 1) \cdot a^{d - 2} \cdot \alpha^2
+ (-1)^{h(\zeta) + 1} \cdot 8d \cdot a^d\\
&\quad\quad\quad\quad\quad\quad\quad\quad
 + (-1)^{h(\zeta) + 1} \cdot a^d \cdot (6 \zeta^2 + 2K_X^2)\\
&= (-1)^{h(\zeta) + 1} \cdot
\left \{ d(d - 1) \cdot a^{d - 2} \cdot \alpha^2 +
(2K_X^2 + 2d + 6) \cdot a^d \right \}. \qed\\
\endalign$$

For $[\mu_+(\alpha)]^{d - 2} \cdot \nu_+
 - [\mu_-(\alpha)]^{d - 2} \cdot \nu_-$, we have the following.

\theorem{6.5} Let $\zeta$ define a wall of type $(w, p)$ with
$\ell_\zeta = 1$, let $d = -p - 3$. For $\alpha \in H_2(X; \Zee)$,
let $a = (\zeta \cdot \alpha)/2$. Then,
$[\mu_+(\alpha)]^{d - 2} \cdot \nu_+
 - [\mu_-(\alpha)]^{d - 2} \cdot \nu_-$ is equal to
$${1 \over 4} \cdot (-1)^{h(\zeta)} \cdot \left [
(d - 2)(d - 3) \cdot a^{d - 4} \cdot \alpha^2
+ (2K_X^2 + 2d - 18) \cdot a^{d - 2} \right ].$$
\endstatement
\par\noindent
{\it Proof.} By Theorem 5.5,
$[\mu_+(\alpha)]^{d - 2} \cdot \nu_+
 - [\mu_-(\alpha)]^{d - 2} \cdot \nu_-$ is equal to
$$\align
&\quad {1 \over 4} \cdot \sum_{j = 0}^{2} {{d - 2} \choose j} \cdot
(-1)^{h(\zeta) + j} \cdot a^{d - 2 - j} \cdot S_j
 - (-1)^{h(\zeta)} \cdot a^{d - 2} \cdot \sum_{k = 0}^1
([{\Cal Z_{1 - k}}] + [{\Cal Z_{k}}])/x \\
&= {1 \over 4} \cdot \sum_{j = 0}^{2} {{d - 2} \choose j} \cdot
(-1)^{h(\zeta) + j} \cdot a^{d - 2 - j} \cdot S_j
 - (-1)^{h(\zeta)} \cdot 2a^{d - 2}.\\
\endalign$$
By Proposition 5.9, Proposition 5.12, and Lemma 6.3, we have
$$S_2 = 2 \alpha^2, S_1 = -8a, S_0 = 6 \zeta^2 + 2K_X^2.$$
Therefore, we conclude that $[\mu_+(\alpha)]^{d - 2} \cdot \nu_+
 - [\mu_-(\alpha)]^{d - 2} \cdot \nu_-$ is equal to
$${1 \over 4} \cdot (-1)^{h(\zeta)} \cdot \left [
(d - 2)(d - 3) \cdot a^{d - 4} \cdot \alpha^2
+ (2K_X^2 + 2d - 18) \cdot a^{d - 2} \right ]. \qed$$

In the rest of this section, we assume that $\ell_\zeta = 2$. The following
standard facts about double coverings can be found in \cite{2, 10}.

\lemma{6.6} Let $\phi: Y_1 \to Y_2$ be a double covering between
two smooth projective varieties with
$\phi_*\Cal O_{Y_1} = \Cal O_{Y_2} \oplus L^{-1}$
where $L$ is a line bundle on $Y_2$.

\roster
\item"{(i)}" $K_{Y_1} = \phi^*(K_{Y_2} \otimes L)$ and
$L^{\otimes 2} = \Cal O_{Y_2}(B)$
where $B$ is the branch locus in $Y_2$ and is the image of
the fixed set of the involution $\iota$ on $Y_1$;

\item"{(ii)}" If $D$ is a divisor on $Y_1$, then $\phi_*(\Cal O_{Y_1}(D))$
is a rank $2$ bundle on $Y_2$ with
$c_1(\phi_*(\Cal O_{Y_1}(D))) = \phi_*D - L$ and
$$c_2(\phi_*(\Cal O_{Y_1}(D))) = {1 \over 2} \cdot \left [(\phi_*D)^2
 - \phi_*(D^2) - \phi_*D \cdot L \right].$$

\endroster
\endstatement

Next, we recall some standard facts about the Hilbert scheme
$H_2 = \Hilb ^2(X)$. Let $\Delta_0 \subset X \times X$ be the diagonal,
and let $\iota$ be the obvious involution on
$\tilde H_2 = \operatorname{Bl}_{\Delta_0} (X\times X)$,
the blowup of $X\times X$ along $\Delta_0$.
Let $E$ be the exceptional divisor of the blowup in $\tilde H_2$.
Then, $H_2 = \tilde H_2/\iota$ and the branch locus lies under $E$.
Let $\tilde \Cal Z_2 \subset X \times \tilde H_2$ be the pullback of
the codimension $2$ cycle $\Cal Z_2 \subset X \times H_2$.
Then, $\tilde \Cal Z_2$ splits into a union of two cycles $\tilde H_{12}$
and $\tilde H_{13}$ in $X \times \tilde H_2$, which are the proper transforms
in $X \times \tilde H_2$ of the two morphisms of $X \times X$ into
$X \times (X \times X)$: the first maps the first factor
in $X \times X$ diagonally into $X \times X$ which is the product of
the first and second factors in $X \times (X \times X)$,
while the second maps the first factor
in $X \times X$ diagonally into $X \times X$ which is the product of
the first and third factors in $X \times (X \times X)$.
Thus each $\tilde H_{1j}$ is isomorphic to
$\operatorname{Bl}_{\Delta_0} (X\times X)$, and the projection of each
to $\tilde H_2$ is an isomorphism. If $\alpha \in H_2(X; \Zee)$, then
$$[\tilde \Cal Z_2]/\alpha = \alpha \otimes 1 + 1 \otimes \alpha
= \alpha \otimes 1 + \iota^*(\alpha \otimes 1) \eqno(6.7)$$
where $\alpha \otimes 1$ and $1 \otimes \alpha$ are the pull-backs
of $\alpha$ by the two projections of $\tilde H_2$ to $X$.
Fix $x \in X$. Let $\tilde X_x$ be the pull-back of
$X \times x \subset X \times X$ to $\tilde H_2$.
Then, $\tilde X_x$ is isomorphic to the blow-up of $X$ at $p$ with
the exceptional divisor $(\tilde X_x \cap E)$; moreover,
$$[\tilde \Cal Z_2]/x = \tilde X_x + \iota^* \tilde X_x. \eqno(6.8)$$
It is known (see p. 685 in \cite{9}) that $\Cal Z_2$ is smooth.
Let $B$ be the branch locus of the natural double covering from
$\Cal Z_2$ to $H_2$. Then, $B \sim 2L$ for some divisor $L$ on $H_2$,
and the pull-back of $B \subset H_2$ to $\tilde H_2$ is $2E$.
Let $i: \Cal Z_2 \to X \times H_2$ be the embedding,
and $\pi_1$ and $\pi_2$ be the natural projections of $X \times H_2$
to $X$ and $H_2$ respectively.

In the following, we compute the Chern and Segre classes of
$\Cal E_{\zeta}^{2 - k, k}$ for $k = 0, 1, 2$.
The method is to use Lemma 5.11 together with Lemma 6.6.
We start with $\Cal E_{\zeta}^{2, 0}$.

\lemma{6.9} $c_3(\Cal E_{\zeta}^{2, 0}) = c_4(\Cal E_{\zeta}^{2, 0}) = 0$,
$c_1(\Cal E_{\zeta}^{2, 0}) = [\Cal Z_2]/\zeta - L$,
and
$$c_2(\Cal E_{\zeta}^{2, 0}) = {1 \over 2} \left [([\Cal Z_2]/\zeta)^2
 - \zeta^2 \cdot X_x - [\Cal Z_2]/\zeta \cdot L \right]$$
where $x$ is any point on $X$, and $X_x$ stands for $[\Cal Z_2]/x$.
\endstatement
\proof Let notations be as in Lemma 5.11, and let $\ell_{\zeta} = 2$ and
$k = 0$. Then, $ Ext^1 = 0$. By Lemma 5.11 (i),
$\Cal E_{\zeta}^{2, 0}$ sits in an exact sequence
$$0 \to (\pi_2 \cdot i)_*(\pi_1 \cdot i)^*\Cal O_X(\zeta) \to
\Cal E_{\zeta}^{2, 0} \to [\Cal O_{H_2}]^{\oplus~ h(\zeta)} \to 0.$$
Since $(\pi_2 \cdot i)_*(\pi_1 \cdot i)^*\Cal O_X(\zeta)$ has rank $2$,
$c_3(\Cal E_{\zeta}^{2, 0}) = c_4(\Cal E_{\zeta}^{2, 0}) = 0$.
By Lemma 6.6 (ii),
$$c_1(\Cal E_{\zeta}^{2, 0}) = (\pi_2 \cdot i)_*(\pi_1 \cdot i)^*\zeta - L
= [\Cal Z_2]/\zeta - L$$
since $(\pi_2 \cdot i)_*(\pi_1 \cdot i)^*\zeta = [\Cal Z_2]/\zeta$;
moreover, we have
$$\align
c_2(\Cal E_{\zeta}^{2, 0}) &=
{1 \over 2} \left [((\pi_2 \cdot i)_*(\pi_1 \cdot i)^*\zeta)^2
 - (\pi_2 \cdot i)_*((\pi_1 \cdot i)^*\zeta)^2 -
(\pi_2 \cdot i)_*(\pi_1 \cdot i)^*\zeta \cdot L \right]\\
&= {1 \over 2} \left [([\Cal Z_2]/\zeta)^2 - \zeta^2 \cdot X_x
 - [\Cal Z_2]/\zeta \cdot L \right]\\
\endalign$$
since $(\pi_2 \cdot i)_*((\pi_1 \cdot i)^*\zeta)^2 = \zeta^2 \cdot
(\pi_2 \cdot i)_*(\pi_1 \cdot i)^*x = \zeta^2 \cdot [\Cal Z_2]/x
= \zeta^2 \cdot X_x$.
\endproof

The following follows from Lemma 6.9 and Remark 5.6.

\corollary{6.10} The Segre classes of the bundle $\Cal E_{\zeta}^{2, 0}$
are given by
$$\align
&s_1(\Cal E_{\zeta}^{2, 0}) = L - [\Cal Z_2]/\zeta\\
&s_2(\Cal E_{\zeta}^{2, 0}) = {1 \over 2} \left[[[\Cal Z_2]/\zeta]^2
 - 3 [\Cal Z_2]/\zeta \cdot L + 2L^2 + \zeta^2 \cdot X_x \right]\\
&s_3(\Cal E_{\zeta}^{2, 0}) = [\Cal Z_2]/\zeta]^2 \cdot L -
2 [\Cal Z_2]/\zeta \cdot L^2 + L^3 - \zeta^2 \cdot X_x \cdot [\Cal Z_2]/\zeta
+ \zeta^2 \cdot X_x \cdot L\\
&s_4(\Cal E_{\zeta}^{2, 0}) = {(\zeta^2)^2 \over 2} - 5\zeta^2
 - {5 \over 2} \zeta \cdot K_X + (6\chi(\Cal O_X) - K_X^2).\\
\endalign$$
Here we have identified degree $4$ classes with the corresponding integers.
\endstatement
\par\noindent
{\it Proof.} Since the computation is straightforward, we only calculate
$s_4(\Cal E_{\zeta}^{2, 0})$. For simplicity,
let $c_i$ denote the $i^{\text{th}}$ Chern class of $\Cal E_{\zeta}^{2, 0}$.
Note that $c_3 = c_4 = 0$ by Lemma 6.9. Thus,
$s_4(\Cal E_{\zeta}^{2, 0}) = c_1^4 - 3 c_1^2 c_2 + c_2^2$ by Remark 5.6.
Therefore,
$$\align
s_4(\Cal E_{\zeta}^{2, 0}) &= ([\Cal Z_2]/\zeta - L)^4 -
3 ([\Cal Z_2]/\zeta - L)^2 \cdot {1 \over 2} \left [([\Cal Z_2]/\zeta)^2
 -  \zeta^2 \cdot X_x - [\Cal Z_2]/\zeta \cdot L \right]\\
&\quad\quad\quad\quad\quad\quad
+ {1 \over 4} \left [([\Cal Z_2]/\zeta)^2
 -  \zeta^2 \cdot X_x - [\Cal Z_2]/\zeta \cdot L \right]^2\\
&= L^4 - {5 \over 2} \cdot [\Cal Z_2]/\zeta \cdot L^3 +
{7 \over 4} \cdot ([\Cal Z_2]/\zeta)^2 \cdot L^2 +
{3 \over 2}\zeta^2 \cdot X_x \cdot L^2\\
&\quad\quad\quad\quad\quad\quad
 -  {1 \over 4} ([\Cal Z_2]/\zeta)^4 + {1 \over 4}(\zeta^2)^2 \cdot X_x^2
+ \zeta^2 \cdot ([\Cal Z_2]/\zeta)^2 \cdot X_x\\
\endalign$$
since $([\Cal Z_2]/\zeta)^3 \cdot L = 0 =
[\Cal Z_2]/\zeta \cdot L \cdot X_x$. Now, we need a claim.

\claim{} Let $\alpha, \beta \in H_2(X; \Zee)$. Then, we have the following:
\roster

\item"{(i)}" $[\Cal Z_2]/\alpha \cdot [\Cal Z_2]/\beta \cdot X_x
= \alpha \cdot \beta$;

\item"{(ii)}" $X_x^2 = 1$;

\item"{(iii)}" $X_x \cdot L^2 = -1$;

\item"{(iv)}" $L^4 = 6\chi(\Cal O_X) - K_X^2$;

\item"{(v)}" $[\Cal Z_2]/\alpha \cdot L^3 = \alpha \cdot K_X$;

\item"{(vi)}" $[\Cal Z_2]/\alpha \cdot [\Cal Z_2]/\beta \cdot L^2
= - 2(\alpha \cdot \beta)$.

\endroster
\endstatement
\par\noindent
{\it Proof.} Let $\pi: \tilde H_2 \to H_2 = \tilde H_2/\iota$
be the quotient map. By (6.8), we have
$$\pi^*X_x = \pi^*([\Cal Z_2]/x) = [\tilde \Cal Z_2]/x =
(\tilde X_x + \iota^* \tilde X_x).$$

(i) Recall from (6.7) that
$\pi^*([\Cal Z_2]/\alpha) = [\tilde \Cal Z_2]/\alpha =
\alpha \otimes 1 + 1 \otimes \alpha$. Thus,
$$\align
[\Cal Z_2]/\alpha \cdot [\Cal Z_2]/\beta \cdot X_x
&= {1 \over 2} \cdot \pi^*([\Cal Z_2]/\alpha) \cdot \pi^*([\Cal Z_2]/\beta)
\cdot \pi^*X_x \\
&= {1 \over 2} \cdot (\alpha \otimes 1 + 1 \otimes \alpha) \cdot
(\beta \otimes 1 + 1 \otimes \beta) \cdot (\tilde X_x + \iota^* \tilde X_x)\\
&= \alpha \cdot \beta.\\
\endalign$$

(ii) Let $x_1 \in X$ be a point different from $x$. Then,
$$\align
X_x^2 &= X_x \cdot X_{x_1} =
{1 \over 2} \cdot \pi^*(X_x) \cdot \pi^*(X_{x_1}) \\
&= {1 \over 2} \cdot (\tilde X_x + \iota^* \tilde X_x) \cdot
(\tilde X_{x_1} + \iota^* \tilde X_{x_1}) \\
&= 1.\\
\endalign$$

(iii) Since $B \sim 2L$ and $\pi^*(B) = 2E$, $\pi^*(L) \sim E$. Thus,
$$X_x \cdot L^2 = {1 \over 2} \cdot (\tilde X_x + \iota^* \tilde X_x)
\cdot E^2 = \tilde X_x \cdot E^2 = (\tilde X_x \cdot E)^2 = -1.$$

(iv) Since $E = \Pee(N\spcheck)$ where $N$ is the normal bundle of $\Delta_0$
in $X \times X$, $-E|E = \xi$ is the tautological line bundle on $E$.
Since $N = T_{\Delta_0}$,
$$\xi^2 = -(\pi|E)^*c_1(N) \cdot \xi - c_2(N) =
(\pi|E)^*K_{\Delta_0} \cdot \xi + (K_X^2 - 12 \chi(\Cal O_X)).$$
It follows that $\xi^3 = (2K_X^2 - 12 \chi(\Cal O_X)) \cdot \xi$. Therefore,
$$L^4 = {1 \over 2} \cdot E^4 = -{1 \over 2} \cdot \xi^3
= 6\chi(\Cal O_X) - K_X^2.$$

(v) Note that $(\alpha \otimes 1)|E = (\pi|E)^*\alpha$
since ${\Delta_0} \cong X$. Thus,
$$[\Cal Z_2]/\alpha \cdot L^3 = {1 \over 2} \cdot
(\alpha \otimes 1 + 1 \otimes \alpha) \cdot E^3
= (\alpha \otimes 1) \cdot E^3
= (\pi|E)^*\alpha \cdot \xi^2 = \alpha \cdot K_X.$$

(vi) Again since $(\alpha \otimes 1)|E = (\pi|E)^*\alpha
= (1 \otimes \alpha)|E$, we have
$$\align
[\Cal Z_2]/\alpha \cdot [\Cal Z_2]/\beta \cdot L^2
&= {1 \over 2} \cdot (\alpha \otimes 1 + 1 \otimes \alpha) \cdot
(\beta \otimes 1 + 1 \otimes \beta) \cdot E^2\\
&= -2 \cdot (\pi|E)^*\alpha \cdot (\pi|E)^*\beta \cdot \xi \\
&= -2 (\alpha \cdot \beta).\qed \\
\endalign$$

We continue the calculation of $s_4(\Cal E_{\zeta}^{2, 0})$.
By Lemma 5.8 (i), $([\Cal Z_2]/\zeta)^4 = 3 (\zeta^2)^2$.
It follows from the above Claim with a straightforward computation that
$$s_4(\Cal E_{\zeta}^{2, 0}) = {(\zeta^2)^2 \over 2} - 5\zeta^2
 -  {5 \over 2} \zeta \cdot K_X + (6\chi(\Cal O_X) - K_X^2). \qed$$

Next, we compute the Chern and Segre classes of $\Cal E_{\zeta}^{0, 2}$
on $H_2$.

\lemma{6.11} $c_3(\Cal E_{\zeta}^{0, 2}) = c_4(\Cal E_{\zeta}^{0, 2}) = 0$,
$c_1(\Cal E_{\zeta}^{0, 2}) = [\Cal Z_2]/(\zeta - K_X) + L$,
and
$$c_2(\Cal E_{\zeta}^{0, 2}) = {1 \over 2} \left [
L \cdot [\Cal Z_2]/(\zeta - K_X) + [[\Cal Z_2]/(\zeta - K_X)]^2
 -  (\zeta - K_X)^2 \cdot X_x \right].$$
\endstatement
\proof Let $\ell_{\zeta} = 2$ and $k = 2$ in Lemma 5.11. By Lemma 6.6 (i),
$$(\det T_{\Cal Z_2})^{-1} = \Cal O_{\Cal Z_2}(K_{\Cal Z_2})
= (\pi_2 \cdot i)^*\Cal O_{H_2}(K_{H_2} + L).$$
Let $N_{\Cal Z_2}$ be the normal bundle of $\Cal Z_2$ in $X \times H_2$.
Since $\Cal Z_2$ is smooth and has codimension $2$ in $X \times H_2$,
$Ext^1 = Ext^1(I_{\Cal Z_2}, \Cal O_{X \times H_2})$
is isomorphic to
$$\det   N_{\Cal Z_2} = i^*\det   T_{X \times H_2} \otimes
(\det   T_{\Cal Z_2})^{-1} =
\Cal O_{\Cal Z_2}((\pi_2 \cdot i)^*L - (\pi_1 \cdot i)^*K_X).$$
By Lemma 5.11 (i), $\Cal E_{\zeta}^{0, 2}$ sits in an exact sequence
$$0 \to [\Cal O_{H_2}]^{\oplus~ h(\zeta)} \to \Cal E_{\zeta}^{0, 2} \to
(\pi_2 \cdot i)_*\Cal O_{\Cal Z_2}((\pi_2 \cdot i)^*L +
(\pi_1 \cdot i)^*(\zeta - K_X)) \to 0.$$
Note that $(\pi_2 \cdot i)_*(\pi_2 \cdot i)^*L = 2L$. Thus, by Lemma 6.6 (ii),
$$c_1(\Cal E_{\zeta}^{0, 2}) = (\pi_2 \cdot i)_*\left[(\pi_2 \cdot i)^*L +
(\pi_1 \cdot i)^*(\zeta - K_X)\right ] - L = [\Cal Z_2]/(\zeta - K_X) + L.$$
Also, Lemma 6.6 (ii) together with a straightforward calculation gives
$$c_2(\Cal E_{\zeta}^{0, 2}) = {1 \over 2} \left [
L \cdot [\Cal Z_2]/(\zeta - K_X) + [[\Cal Z_2]/(\zeta - K_X)]^2
 -  (\zeta - K_X)^2 \cdot X_x \right ]$$
where we have used the projection formula
$$(\pi_2 \cdot i)_*[(\pi_2 \cdot i)^*L \cdot (\pi_1 \cdot i)^*(\zeta - K_X)]
= L \cdot (\pi_2 \cdot i)_*(\pi_1 \cdot i)^*(\zeta - K_X)$$
and the fact that $(\pi_2 \cdot i)_*(\pi_2 \cdot i)^*L^2 = 2L^2$.
\endproof

The following follows from Lemma 6.11 and Remark 5.6.

\corollary{6.12} The Segre classes of $\Cal E_{\zeta}^{0, 2}$ are given by
$$\align
&s_1(\Cal E_{\zeta}^{0, 2}) = [\Cal Z_2]/(K_X - \zeta) - L\\
&s_2(\Cal E_{\zeta}^{0, 2}) = {1 \over 2} \left[[[\Cal Z_2]/(\zeta - K_X)]^2
+ 3 [\Cal Z_2]/(\zeta - K_X) \cdot L + 2L^2
+ (\zeta - K_X)^2 \cdot X_x \right]\\
&s_3(\Cal E_{\zeta}^{0, 2}) = -[[\Cal Z_2]/(\zeta - K_X)]^2 \cdot L -
2 [\Cal Z_2]/(\zeta - K_X) \cdot L^2 - L^3\\
&\quad\quad\quad\quad\quad\quad
 -  (\zeta - K_X)^2 \cdot X_x \cdot [\Cal Z_2]/(\zeta - K_X)
 -  (\zeta - K_X)^2 \cdot X_x \cdot L\\
&s_4(\Cal E_{\zeta}^{0, 2}) =
{((K_X - \zeta)^2)^2 \over 2} - 5(K_X - \zeta)^2
 -  {5 \over 2} (K_X - \zeta) \cdot K_X + (6\chi(\Cal O_X) - K_X^2).\\
\endalign$$
\endstatement
\proof The calculation of $s_4(\Cal E_{\zeta}^{0, 2})$ is similar
to that of $s_4(\Cal E_{\zeta}^{2, 0})$ in Corollary 6.10.
\endproof

Note that $s_4(\Cal E_{\zeta}^{0, 2})$ may be
obtained from $s_4(\Cal E_{\zeta}^{2, 0})$ by replacing $\zeta$
by $K_X - \zeta$, and indeed this holds more generally for $s_i$ when we add
the sign $(-1)^i$.

Now we compute the Chern and Segre classes of
$\Cal E_{\zeta}^{1, 1}$ on $X \times X$.

\lemma{6.13} Let $\tau_1$ and $\tau_2$ be the two natural projections
of $X \times X$ to $X$, let ${\Delta_0}$ be the diagonal in $X \times X$, and
let $j\: \Delta _0 \to X\times X$ be the inclusion. Then
$$\align
&c_1(\Cal E_{\zeta}^{1, 1}) = \tau_1^*\zeta + \tau_2^*(\zeta - K_X) \\
&c_2(\Cal E_{\zeta}^{1, 1}) = \tau_1^*\zeta \cdot \tau_2^*(\zeta - K_X)
+ {\Delta_0} \\
&c_3(\Cal E_{\zeta}^{1, 1}) = \tau_1^*\zeta \cdot {\Delta_0} -
\tau_2^*(\zeta - K_X) \cdot {\Delta_0} - j_*K_{\Delta_0} \\
&c_4(\Cal E_{\zeta}^{1, 1}) = -{{K_{\Delta_0}^2} \over 2}.\\
\endalign$$
\endstatement
\proof Let $\ell_{\zeta} = 2$ and $k = 1$ in Lemma 5.11. Recall that
$\pi_1$ and $\pi_2$ are the natural projections of $X \times (X \times X)$
to $X$ and $(X \times X)$ respectively.

\claim{1} $\pi_{2*}\left (\pi_1^*\Cal O_X(\zeta)
\otimes Ext^1 \right) \cong \tau_2^*\Cal O_X(\zeta - K_X)
\otimes I_{\Delta_0}$.
\endstatement
\par\noindent
{\it Proof.} Let $\Delta_{12}$ be the diagonal in $X \times X$ which is formed
by the first and second factors in $X \times (X \times X)$, and let
$\Delta_{13}$ be the diagonal in $X \times X$ which is formed
by the first and third factors in $X \times (X \times X)$. Then,
$\Delta_{12} \times X$ and $\Delta_{13} \times X$ are smooth codimension $2$
subvarieties in $X \times (X \times X)$. Here it is understood that
the factor $X$ in $\Delta_{13} \times X$ is embedded as the second
factor in $X \times (X \times X)$. Moreover,
$\Delta_{12} \times X$ and $\Delta_{13} \times X$ intersect properly
along the diagonal $\Delta_{123}$ in $X \times X \times X$.
Thus, from Lemma 5.10 (iii), we conclude that
$$Ext^1 = Ext^1(I_{\Delta_{13} \times X}, I_{\Delta_{12} \times X})
\cong I \otimes \det   N$$
where $N$ is the normal bundle $\Delta_{13} \times X$ in
$X \times (X \times X)$, and $I$ is the ideal sheaf of $\Delta_{123}$
in $\Delta_{13} \times X$. Now, the restriction of $\pi_2$ to
$\Delta_{13} \times X$ gives an isomorphism from $\Delta_{13} \times X$
to $X \times X$. Via this isomorphism,
$\Delta_{123}$ in $\Delta_{13} \times X$ is identified with
the diagonal ${\Delta_0}$ in $X \times X$, $\det   N$ is identified with
$\tau_2^*(-K_X)$, and the restriction
$\pi_1^*\Cal O_X(\zeta)|(\Delta_{13} \times X)$
is identified with $\tau_2^*(\zeta)$. Therefore,
$$\pi_{2*}\left (\pi_1^*\Cal O_X(\zeta) \otimes Ext^1 \right) \cong
\pi_{2*}\left (\pi_1^*\Cal O_X(\zeta) \otimes I \otimes \det   N \right)
= \tau_2^*\Cal O_X(\zeta - K_X) \otimes I_{\Delta_0}. \qed$$

Note that $\pi_{2*}(\pi_1^*\Cal O_X(\zeta) \otimes
\Cal O_{\Delta_{12} \times X}) = \tau_1^*(\zeta)$.
Thus by Lemma 5.11 (i) and Claim 1, we have a row exact sequence and
a column exact sequence
$$\matrix
          &0&\\
          &\downarrow&\\
          &\tau_1^*(\zeta) &\\
          &\downarrow&\\
0 \to & R^1\pi_{2*} \left (\pi_1^*\Cal O_X(\zeta) \otimes Hom \right )
& \to \Cal E_{\zeta}^{1, 1} \to
\tau_2^*\Cal O_X(\zeta - K_X) \otimes I_{\Delta_0} \to 0. \\
          &\downarrow&\\
          &[\Cal O_{X \times X}]^{\oplus~ h(\zeta)} &\\
          &\downarrow&\\
          &0&\\
\endmatrix \eqno (6.14)$$

In the next claim, we compute the Chern classes of $I_{\Delta_0}$.
Clearly, $c_0(I_{\Delta_0}) = 1$.

\claim{2} $c_1(I_{\Delta_0}) = 0$, $c_2(I_{\Delta_0}) = {\Delta_0}$,
$c_3(I_{\Delta_0}) = -j_*K_{\Delta_0}$,
$c_4(I_{\Delta_0}) = {K_{\Delta_0}^2}/2$.
\endstatement
\proof Note that $\Todd (N_{\Delta_0})^{-1} = 1 + K_{\Delta_0}/2 +
(K_{\Delta_0}^2/4 - \chi (\Cal O_{\Delta_0}))$.
By a formula on p.288 of \cite{12}
(a special case of the Grothendieck-Riemann-Roch Theorem),
$$\align
\ch  (j!\Cal O_{\Delta_0}) &= j_*(\Todd (N_{\Delta_0})^{-1} \cdot
\ch (\Cal O_{\Delta_0})) = j_*(\Todd (N_{\Delta_0})^{-1})\\
&= {\Delta_0} + {{j_*K_{\Delta_0}} \over 2} +
j_*\left({{K_{\Delta_0}^2} \over 4} - \chi (\Cal O_{\Delta_0}) \right).\\
\endalign$$
Since $\ch (j!\Cal O_{\Delta_0})$ is just equal to
$\ch (j_*\Cal O_{\Delta_0})$, we obtain
$$\ch (I_{\Delta_0}) = \ch (\Cal O_{X \times X}) -
\ch (j_*\Cal O_{\Delta_0})
= 1 - {\Delta_0} - {{j_*K_{\Delta_0}} \over 2} -
j_*\left({{K_{\Delta_0}^2} \over 4} - \chi (\Cal O_{\Delta_0}) \right).$$
 From this, the Chern classes of $I_{\Delta_0}$ follows immediately.
In particular,
$$c_4(I_{\Delta_0}) = {{{\Delta_0}^2} \over 2} +
j_*\left({{3K_{\Delta_0}^2} \over 2} - 6\chi (\Cal O_{\Delta_0}) \right)
= {K_{\Delta_0}^2 \over 2}$$
since ${\Delta_0}^2 = c_2(T_X) = 12 \chi (\Cal O_X) - K_X^2$
(see the Example 8.1.12 in \cite{12}).
\endproof

Now the calculation of the Chern classes of $\Cal E_{\zeta}^{1, 1}$ follows
from
(6.14)  and Claim 2. In particular,
$$\align
c_4(\Cal E_{\zeta}^{1, 1}) &= -\tau_1^*\zeta \cdot \tau_2^*(\zeta -
K_X) \cdot {\Delta_0} - \tau_1^*\zeta \cdot j_*K_{\Delta_0}
+ \tau_2^*(\zeta - K_X)^2 \cdot {\Delta_0}\\
&\quad\quad\quad\quad + 2 \tau_2^*(\zeta - K_X) \cdot j_*K_{\Delta_0}
+ {{j_*K_{\Delta_0}^2} \over 2}= -{{K_{\Delta_0}^2} \over 2}\\
\endalign$$
since $\tau_1^*\zeta \cdot \tau_2^*(\zeta - K_X)
\cdot {\Delta_0} = \zeta \cdot (\zeta - K_X)$ and
$\tau_1^*\zeta \cdot j_*K_{\Delta_0} = \zeta \cdot K_X$.
\endproof

The next result follows immediately from Lemma 6.13 and Remark 5.6.

\corollary{6.15} Let notations be the same as in Lemma 6.13. Then
$$\align
&s_1(\Cal E_{\zeta}^{1, 1}) = -\tau_1^*\zeta - \tau_2^*(\zeta - K_X) \\
&s_2(\Cal E_{\zeta}^{1, 1}) = \tau_1^*\zeta^2 +
\tau_1^*\zeta \cdot \tau_2^*(\zeta - K_X) + \tau_2^*(\zeta - K_X)^2
 -  {\Delta_0} \\
&s_3(\Cal E_{\zeta}^{1, 1}) = -\tau_1^*\zeta^2 \cdot \tau_2^*(\zeta - K_X)
 -  \tau_1^*\zeta \cdot \tau_2^*(\zeta - K_X)^2\\
&\quad\quad\quad\quad\quad + \tau_1^*\zeta \cdot {\Delta_0}
+ 3 \tau_2^*(\zeta - K_X) \cdot {\Delta_0} + j_*K_{\Delta_0}\\
&s_4(\Cal E_{\zeta}^{1, 1}) = (12 \zeta \cdot K_X - 12 \zeta^2 - 3K_X^2).\qed
\endalign$$
\endstatement

We can now work out (5.7) explicitly for $\ell_\zeta = 2$ and $j = 2, 1, 0$.
For simplicity, let
$$S_j = \sum_{k = 0}^2~ S_{j, k} = \sum_{k = 0}^2~
([{\Cal Z_{2 - k}}]/\alpha + [{\Cal Z_{k}}]/\alpha)^j
\cdot s_{4 - j}(\Cal E_{\zeta}^{2 - k, k} \oplus
(\Cal E_{-\zeta}^{k, 2 - k})\spcheck). \eqno (6.16)$$

\lemma{6.17} $S_2 = 64a^2 + (12 \zeta^2 + 4K_X^2 - 20)\alpha^2$
where $a = (\zeta \cdot \alpha)/2$.
\endstatement
\proof Note that $s_i(\Cal E_{-\zeta}^{k, 2 - k})$ (respectively, $S_{j, k}$)
can be obtained from $s_i(\Cal E_{\zeta}^{k, 2 - k})$
(respectively, $(-1)^j \cdot S_{j, 2 - k}$) by replacing $\zeta$ by $-\zeta$.
Also, $S_{2, 2}$ is equal to
$$([{\Cal Z_2}]/\alpha)^2 \cdot s_2(\Cal E_{\zeta}^{0, 2} \oplus
(\Cal E_{-\zeta}^{2, 0})\spcheck)
= ([{\Cal Z_2}]/\alpha)^2 \cdot \left [
s_2(\Cal E_{\zeta}^{0, 2}) - s_1(\Cal E_{\zeta}^{0, 2}) \cdot
s_1(\Cal E_{-\zeta}^{2, 0}) + s_2(\Cal E_{-\zeta}^{2, 0}) \right].$$
Therefore, by Corollary 6.10 and Corollary 6.12, we obtain
$$S_{2, 2} + S_{2, 0} = 32a^2 + (6 \zeta^2 + 2K_X^2 - 12)\alpha^2
+ 2(\alpha \cdot K_X)^2.$$
Let $\tau_1$ and $\tau_2$ be the projections of $X \times X$ to $X$.
Then, by Corollary 6.15,
$$\align
S_{2, 1} &= (\tau_1^*\alpha + \tau_2^*\alpha)^2
\cdot s_2(\Cal E_{\zeta}^{1, 1} \oplus (\Cal E_{-\zeta}^{1, 1})\spcheck)\\
&= (\tau_1^*\alpha + \tau_2^*\alpha)^2 \cdot \left[
s_2(\Cal E_{\zeta}^{1, 1}) - s_1(\Cal E_{\zeta}^{1, 1}) \cdot
s_1(\Cal E_{-\zeta}^{1, 1}) + s_2(\Cal E_{-\zeta}^{1, 1}) \right]\\
&= 32a^2 + (6 \zeta^2 + 2K_X^2 - 8)\alpha^2 - 2(\alpha \cdot K_X)^2.\\
\endalign$$
It follows that $S_2 = (S_{2, 2} + S_{2, 0}) + S_{2, 1} =
64a^2 + (12 \zeta^2 + 4K_X^2 - 20)\alpha^2$.
\endproof

Next, adopting the same method as in the proof of Lemma 6.17,
we compute the values of $S_1$ and $S_0$ in the next two lemmas respectively.

\lemma{6.18} $S_1 = -(48 \zeta^2 + 16K_X^2 - 120) a$
where $a = (\zeta \cdot \alpha)/2$.
\endstatement
\proof In view of (6.16), we have to compute $S_{1, 2}, S_{1, 1}$,
and $S_{1, 0}$. Note that $S_{1, 0}$ can be obtained from $-S_{1, 2}$
by replacing $\zeta$ by $-\zeta$. Using Corollary 6.10 and Corollary 6.12,
we see that $(S_{1, 2} + S_{1, 0}) = -(24 \zeta^2 + 8K_X^2 - 72) a
 -  6(\zeta \cdot K_X) (\alpha \cdot K_X)$.
Let $\tau_1$ and $\tau_2$ be the projections of $X \times X$ to $X$.
Then, by Corollary 6.15,
$$S_{1, 1} = (\tau_1^*\alpha + \tau_2^*\alpha)
\cdot s_3(\Cal E_{\zeta}^{1, 1} \oplus (\Cal E_{-\zeta}^{1, 1})\spcheck)
= -(24 \zeta^2 + 8K_X^2 - 48) a + 6(\zeta \cdot K_X) (\alpha \cdot K_X).$$
It follows that $S_1 = (S_{1, 2} + S_{1, 0}) + S_{1, 1} =
 - (48 \zeta^2 + 16K_X^2 - 120) a$.
\endproof

\lemma{6.19} $S_0 = 18 (\zeta^2)^2 + (14K_X^2 - 105) \zeta^2
+ [2(K_X^2)^2 - 50K_X^2 + 96]$.
\endstatement
\par\noindent
{\it Proof.} We need to compute $S_{0, 2}, S_{0, 1}$, and $S_{0, 0}$.
Again, $S_{0, 0}$ can be obtained from $S_{0, 2}$
by replacing $\zeta$ by $-\zeta$. Using Corollary 6.10 and Corollary 6.12,
we see that
$$(S_{0, 2} + S_{0, 0}) = 9 (\zeta^2)^2 + (8K_X^2 - 63) \zeta^2
+ [(K_X^2)^2 - 43K_X^2 + 60].$$
By Corollary 6.15, $S_{0, 1} = 9 (\zeta^2)^2 + (6K_X^2 - 42) \zeta^2
+ [(K_X^2)^2 - 7K_X^2 + 36]$. Therefore,
$$S_0 = (S_{0, 2} + S_{0, 0}) + S_{0, 1} = 18 (\zeta^2)^2 +
(14K_X^2 - 105) \zeta^2 + [2(K_X^2)^2 - 50K_X^2 + 96]. \qed$$

Now we can calculate the difference $\delta^X_{w, p}(\Cal C_-, \Cal C_+)$
when $\ell_\zeta = 2$.

\theorem{6.20} Let $\zeta$ define a wall of type $(w, p)$ with
$\ell_\zeta = 2$. Then
$$[\mu_+(\alpha)]^d - [\mu_-(\alpha)]^d = (-1)^{h(\zeta)} \cdot \left \{
g_0 \cdot a^d
+  g_1 \cdot a^{d - 2} \cdot \alpha^2
+  g_2 \cdot a^{d - 4} \cdot (\alpha^2)^2 \right \}$$
for $\alpha \in H_2(X; \Zee)$, where $a$ stands for $(\zeta \cdot \alpha)/2$
and
$$\align
&g_2 = {{d!} \over {2! \cdot (d - 4)!}}\\
&g_1 = {d \choose 2} \cdot (4K_X^2 + 4d + 8)\\
&g_0 = 2d^2 + 2d \cdot K_X^2 + 2 (K_X^2)^2 + 13d + 20K_X^2 + 21.\\
\endalign$$
In other words, the difference $\delta^X_{w, p}(\Cal C_-, \Cal C_+)$
is equal to
$$\delta(\Delta) \cdot (-1)^{h(\zeta)} \cdot \left \{
g_0 \cdot \left({\zeta \over 2} \right)^d
+ g_1 \cdot \left({\zeta \over 2} \right)^{d - 2} \cdot q_X
+ g_2 \cdot \left({\zeta \over 2} \right)^{d - 4} \cdot q_X^2 \right \}.$$
\endstatement
\par\noindent
{\it Proof.} In view of Theorem 5.4 and the notation (6.16), we have
$$[\mu_+(\alpha)]^d - [\mu_-(\alpha)]^d = \sum_{j = 0}^{4}~ {d \choose j}
\cdot (-1)^{h(\zeta) + j} \cdot a^{d - j} \cdot S_j.$$
Now, $S_4$ and $S_3$ are given by Proposition 5.9 and Proposition 5.12
respectively; $S_2, S_1$, and $S_0$ are computed in the previous three lemmas.
So it follows that the coefficient
of $(-1)^{h(\zeta)} \cdot a^{d - 4} \cdot (\alpha^2)^2$ is equal to
$$g_2 = {{d!} \over {2! \cdot (d - 4)!}}.$$
Similarly, also keeping in mind that $\zeta^2 = (p + 8)
= (5 - d)$, we have
$$\align
&g_1 = {d \choose 2} \cdot (12 \zeta^2 + 4K_X^2 + 16d - 52)
= {d \choose 2} \cdot (4K_X^2 + 4d + 8)\\
&g_0 = 64 \cdot {d \choose 2}
    + (48 \zeta^2 + 16K_X^2 - 120) \cdot d~ + \\
&\quad\quad\quad + \left [18(\zeta^2)^2 + 14 \cdot \zeta^2 \cdot K_X^2
+ 2 (K_X^2)^2 - 105 \zeta^2 - 50 K_X^2 + 96 \right ]\\
&\quad= 2d^2 + 2d \cdot K_X^2 + 2 (K_X^2)^2 + 13d + 20K_X^2 + 21. \qed \\
\endalign$$

\corollary{6.21} Let $\zeta$ define a wall of type $(w, p)$ with
$\ell_\zeta \le 2$. Then, the difference $\delta^X_{w, p}(\Cal C_-, \Cal C_+)$
of Donaldson polynomial invariants is
a polynomial in $\zeta$ and $q_X$ with coefficients involving only
$\zeta^2$, homotopy invariants of $X$, and universal constants.
\endstatement
\proof Follows from Theorems 6.1, 6.4, and 6.20.
\endproof

Finally, we compute the difference $[\mu_+(\alpha)]^{d - 2} \cdot \nu_+
 -  [\mu_-(\alpha)]^{d - 2} \cdot \nu_-$ for $\ell_\zeta = 2$.

\theorem{6.22} Let $\zeta$ define a wall of type $(w, p)$ with
$\ell_\zeta = 2$, and let $d = -p - 3$. Then,
$[\mu_+(\alpha)]^{d - 2} \cdot \nu_+
 -  [\mu_-(\alpha)]^{d - 2} \cdot \nu_-$ is equal to
$${1 \over 4} \cdot (-1)^{h(\zeta) + 1} \cdot \left \{
{\tilde g}_0 \cdot a^{d - 2}
+ {\tilde g}_1 \cdot a^{d - 4} \cdot \alpha^2
+ {\tilde g}_2 \cdot a^{d - 6} \cdot (\alpha^2)^2 \right \}$$
for $\alpha \in H_2(X; \Zee)$, where $a$ stands for $(\zeta \cdot \alpha)/2$
and
$$\align
&{\tilde g}_2 = {{(d - 2)!} \over {2! \cdot (d - 6)!}}\\
&{\tilde g}_1 = {{d - 2} \choose 2} \cdot (4K_X^2 + 4d - 40)\\
&{\tilde g}_0 = 2d^2 + 2d \cdot K_X^2 + 2 (K_X^2)^2 - 35d - 28K_X^2 - 99.\\
\endalign$$
\endstatement
\proof By Theorem 5.5, $[\mu_+(\alpha)]^{d - 2} \cdot \nu_+
 -  [\mu_-(\alpha)]^{d - 2} \cdot \nu_-$ is equal to
$${1 \over 4} \cdot \sum_{j = 0}^4 {{d - 2} \choose j} \cdot
(-1)^{h(\zeta) + 1 + j} \cdot a^{d - 2 - j} \cdot S_j
 -  \sum_{j = 0}^2 {{d - 2} \choose j} \cdot
(-1)^{h(\zeta) + 1 + j} \cdot a^{d - 2 - j} \cdot T_j$$
where for simplicity we have defined $T_j$ as
$$T_j = \sum_{k = 0}^2 T_{j, k} = \sum_{k = 0}^2
([{\Cal Z_{2 - k}}]/\alpha + [{\Cal Z_{k}}]/\alpha)^j \cdot
([{\Cal Z_{2 - k}}] + [{\Cal Z_{k}}])/x \cdot
s_{2 - j}(\Cal E_{\zeta}^{2 - k, k}
\oplus (\Cal E_{-\zeta}^{k, 2 -k})\spcheck).$$

Next, we compute $T_0$. Using Corollary 6.10 and Corollary 6.12, we obtain
$$\align
T_{0, 0} &= X_x \cdot s_2(\Cal E_{\zeta}^{2, 0}
\oplus (\Cal E_{-\zeta}^{0, 2})\spcheck) \\
&= X_x \cdot \left [s_2(\Cal E_{\zeta}^{2, 0}) - s_1(\Cal E_{\zeta}^{2, 0})
\cdot s_1(\Cal E_{-\zeta}^{0, 2}) + s_2(\Cal E_{-\zeta}^{0, 2}) \right ]\\
&= (3\zeta^2 + 3 \zeta \cdot K_X + K_X^2 - 3).\\
\endalign$$
Note that $T_{0, 2}$ can be obtained from $T_{0, 0}$
by replacing $\zeta$ by $-\zeta$. Thus,
$$T_{0, 2} = (3\zeta^2 - 3 \zeta \cdot K_X + K_X^2 - 3).$$
Similarly, using Corollary 6.15, we get $T_{0, 1} = (6\zeta^2 + 2K_X^2 - 4)$.
Therefore,
$$T_0 = \sum_{k = 0}^2 T_{0, k} = (12\zeta^2 + 4K_X^2 - 10).$$
By similar but much simpler arguments, we conclude that
$T_1 = -16a$ and $T_2 = 4 \alpha^2$.

 From (5.9), (5.12), (6.17), (6.18), and (6.19), we have
$$\align
&S_4 = 12 (\alpha^2)^2, \\
&S_3 = -48 a \cdot \alpha^2,\\
&S_2 = 64a^2 + (12 \zeta^2 + 4K_X^2 - 20)\alpha^2, \\
&S_1 = -(48 \zeta^2 + 16K_X^2 - 120) a, \\
&S_0 = 18 (\zeta^2)^2 + (14K_X^2 - 105) \zeta^2
+ [2(K_X^2)^2 - 50K_X^2 + 96].\\
\endalign$$
Putting all these together, we see that
$[\mu_+(\alpha)]^{d - 2} \cdot \nu_+
 -  [\mu_-(\alpha)]^{d - 2} \cdot \nu_-$ is equal to
$${1 \over 4} \cdot (-1)^{h(\zeta) + 1} \cdot \left \{
{\tilde g}_0 \cdot a^{d - 2}
+ {\tilde g}_1 \cdot a^{d - 4} \cdot \alpha^2
+ {\tilde g}_2 \cdot a^{d - 6} \cdot (\alpha^2)^2 \right \}$$
where ${\tilde g}_0, {\tilde g}_1$, and ${\tilde g}_2$ are as defined in the
statement of Theorem 6.22 above.
\endproof

\enddocument